\documentclass[review]{elsarticle}

\usepackage{hyperref}

\journal{Computers$\And$Mathematics with Applications}

\usepackage{color}
\usepackage{amsmath}
\usepackage{graphicx}
\usepackage{subfigure}
\usepackage{rotating}
\usepackage{multirow, makecell}
\usepackage{geometry}
\usepackage{appendix}

\biboptions{sort&compress}              

\begin{document}

\sloppy

	
\begin{frontmatter}
		
\title{Numerical stability analysis of shock-capturing methods for strong shocks I: second-order MUSCL schemes}

\author[mymainaddress]{Weijie Ren}

\author[mymainaddress]{Wenjia Xie\corref{mycorrespondingauthor}}
\cortext[mycorrespondingauthor]{Corresponding author}
\ead{xiewenjia@nudt.edu.cn}

\author[mymainaddress]{Ye Zhang}

\author[mymainaddress]{Hang Yu}

\author[mymainaddress]{Zhengyu Tian}

\address[mymainaddress]{College of Aerospace Science and Engineering, National University of Defense Technology, Hunan 410073, China}

\begin{abstract}
Modern shock-capturing schemes often suffer from numerical shock anomalies if the flow field contains strong shocks, which may limit their further application in hypersonic flow computations. In the current study, we devote our efforts to exploring the primary numerical characteristics and the underlying mechanism of shock instability for second-order finite-volume schemes. To this end, we, for the first time, develop the matrix stability analysis method for the finite-volume MUSCL approach. Such a linearized analysis method allows to investigate the shock instability problem of the finite-volume shock-capturing schemes in a quantitative and efficient manner. Results of the stability analysis demonstrate that the shock stability of second-order scheme is strongly related to the Riemann solver, Mach number, limiter function, numerical shock structure, and computational grid. Unique stability characteristics associated with these factors for second-order methods are revealed quantitatively with the established method. Source location of instability is also clarified by the matrix stability analysis method. Results show that the shock instability originates from the numerical shock structure. Such conclusions pave the way to better understand the shock instability problem and may shed new light on developing more reliable shock-capturing methods for compressible flows with high Mach number.
\end{abstract}

\begin{keyword}
Finite-volume\sep Shock-capturing\sep Hypersonic\sep Shock instability\sep Matrix stability analysis\sep MUSCL
\end{keyword}
		
\end{frontmatter}
	
\section{Introduction}
Numerical shock prediction is one of the most important issues in computational fluid dynamics, especially for hypersonic flow computations. Shock-capturing methods, based on the mathematical theory of weak solutions, are relatively simple to code and able to simulate any type of flow field, regardless of the presence of shocks \cite{Bonfiglioli2016}. As a result, they are widely used in practical fluid flow simulations. Unfortunately, the computed flowfields involving strong shock waves by modern shock-capturing schemes are often characterized by the appearance of numerical anomalies, such as the carbuncle phenomenon and post-shock oscillations. They frequently arise in hypersonic flow simulations and considerably deteriorate the reliability of the hypersonic heating prediction.\par

As one of the most famous kinds of shock anomalies, the carbuncle phenomenon was first observed by Peery and Imlay \cite{Perry1988} when they computed the supersonic flow field around the blunt-body with the Roe solver. Although the carbuncle phenomenon refers to the spurious solution of blunt-body calculations in which a protuberance grows ahead of the bow shock along the stagnation line \cite{Quirk1994}, it has been generally considered to be a concrete manifestation of numerical shock instabilities. Previous studies have shown that these numerical shock instabilities are more prone to occur in cases where shock-capturing methods own minimal dissipation on discontinuities. Moreover, recent studies \cite{Kitamura2010,Ohwada2013,Ohwada2018} even show that dissipative solvers such as the HLLE scheme \cite{Harten1983,Einfeldt1988} and the Rusanov scheme \cite{Rusanov1961} (local Lax-Friedrichs), the robustness of which is highly appreciated among experts, are also plagued by the shock instability problem. Thus, although significant progress has been made in numerical methods for compressible flows, it is no exaggeration to say that the modern shock capturing is far away from excellence. During the past three decades, extensive studies have been conducted to clarify the numerical characteristics of the carbuncle phenomenon and propose possible cures for these instabilities. Readers are referred to references \cite{Ismail2006,Shen2014,Rodionov2017,Simon2019b,Qu2022} and references therein for detailed literature reviews of this subject. Here, we only give a brief account of the persistent efforts to reveal the underlying mechanism of the carbuncle phenomenon, which is the core issue of the current study.\par

Quirk \cite{Quirk1994} first performs a comprehensive study of the shock instability problem of approximate Riemann solvers. With the linearized perturbation analysis, he points out that low-dissipative Riemann solvers, such as Roe, add no dissipation on the contact and shear waves in the direction parallel to the shock. As a result, the constant interaction of perturbed pressure and density fields triggers shock instability. Pandolfi and D'Ambrosio \cite{Pandolfi2001} extend Quirk's linearized perturbation analysis to more schemes. They also focus on the relationship between perturbed pressure and density fields and argue that schemes are carbuncle-free if they are able to damp out the density perturbations, such as HLL \cite{Harten1983} and van Leer \cite{VanLeer1982}. Gressier and Moschetta \cite{Gressier2005} employ Quirk's linearized perturbation analysis and find that the strict stability and exact resolution of contact discontinuities are incompatible. And based on the analysis of the hybrid schemes, Shen et al.\cite{Shen2014} find that the third flux component (transverse directional momentum flux component) in the y-direction plays an important role in triggering shock instability. Furthermore, based on the linearized analysis, they reveal that the inconsistent normal momentum distribution along the shock is the source of the shock instability. Simon and Mandal \cite{Simon2018a,Simon2019d} believe that both the mass flux component and interface-normal momentum flux component on the interfaces transverse to the shock front are critical to the shock instability. They investigate the stability of HLLE \cite{Harten1983,Einfeldt1988} and HLLC \cite{Toro1994} schemes by linearized analysis and numerical experiments and find that the shock instability in the HLLC scheme could be triggered by the unwanted activation of the anti-diffusive term appearing in mass and interface-normal momentum flux component on the interfaces transverse to the shock front. Xie et al.\cite{Xie2017} conduct a series of numerical tests with various approximate Riemann solvers and find that the scheme that can preserve the mass flux across the normal shock is shock-stable. Combining the linearized analysis with numerical experiments, they argue that the shock instability can be alleviated if the dissipation corresponding to the momentum is added behind the shock.\par

In the researches mentioned above, the shock instability is attributed to the insufficient multidimensional dissipation on contact or shear waves. This is consistent with the conclusions of theoretical analysis \cite{Einfeldt2010} and physical considerations \cite{Xu2010}. However, recent studies by Chen et al.\cite{Chen2018,Chen2018a} reveal that the inappropriate pressure dissipation of the Riemann solver at the vertical transverse face of a shock is deeply connected to the carbuncle phenomenon and shock instability. Such a useful conclusion is obtained by the stability analysis of a novel shock instability analysis model, which keeps the essential characteristics of numerical shock instability. It also gives support for the partial reasonability of Liou's conjecture \cite{Liou2000}, which states that schemes with the mass flux that is independent of the pressure term are not affected by the carbuncle phenomenon. By a linear perturbation analysis, it is found that the pressure dissipation term is equivalent to reducing mass flux perturbations behind shocks, which is found to trigger the shock instability in the view of perturbations \cite{Xie2019towards,Xie2021,Xie2022}. It is noteworthy that Fleischmann et al.\cite{Fleischmann2020,Fleischmann2020shock} attribute the carbuncle phenomenon to the low Mach number problem of Godunov-type schemes. They demonstrate that an excessive acoustic contribution to dissipation in the flux calculation in the transverse direction to the shock front propagation is the prime reason for the numerical instability. The wrong amplitude of pressure fluctuations induced by the low Mach number of the transverse direction will trigger the shock instability. Thus, a fundamentally different approach to stabilize shock-capturing schemes can be achieved by decreasing the viscosity on the acoustic waves for low Mach numbers.\par

Numerical shock structure is demonstrated to be an important factor that will influence shock instability. Based on the matrix stability analysis, Dumbser et al.\cite{Dumbser2004} and Chauvat et al.\cite{Chauvat2005} find that the shock instability is strongly related to the numerical shock structure. When the internal shock conditions are closer to the upstream states, numerical schemes are more prone to be unstable. These conclusions are validated by the numerical experiments on the 1D and 2D normal shocks from Kitamura et al.\cite{Kitamura2013,Kitamura2009,Kitamura2013a} and Xie et al.\cite{Xie2017}. Moreover, based on the matrix stability analysis, Liu et al.\cite{Liu2020} find that there exists a threshold independent of shock strength to trigger instability. According to their work, if the density ratio between both sides of a specific cell interface in the numerical shock structure exceeds a threshold, the computation is unstable. Xie et al.\cite{Xie2021} link the shock instability to the inappropriate entropy production inside the numerical shock structure. By dissipation analysis and numerical experiments, they find that if enough entropy production is guaranteed inside shocks, the instability problem can be successfully eliminated. Zaide and Roe \cite{Zaide2011,Zaide2012} investigate the mechanism of numerical shock structure affecting the shock instability and attribute the shock instability to the nonlinearity of the Hugoniot curve. They argue that the intermediate states in the numerical shock structure are always assumed to satisfy the local thermodynamic equilibrium, while this does not hold inside the physical shock, which is the cause of shock instability. Based on their conclusions, Zaide and Roe\cite{zaide2013,Zaide2013a} develop two novel flux functions which can permit stationary shocks with one intermediate state and an unambiguous numerical shock structure and demonstrate their stability.\par

The computational grid is known as another critical factor that will control the shock instability. Quirk \cite{Quirk1994} finds that the carbuncle phenomenon will be more pronounced if the grid is aligned to the shock. By the numerical tests, Ohwada et al.\cite{Ohwada2013} reveal that even the robust schemes, such as HLLE \cite{Harten1983,Einfeldt1988} and AUSM$ ^+ $ \cite{Liou1996}, can also yield unstable flow field when the grid is aligned to the shock. Based on a series of numerical experiments, Henderson and Menart \cite{Henderson2007} conduct a comprehensive study on the effect of the computational grid on the shock anomalies. They find that the aspect ratio of cells near the shock is a factor influencing the magnitude of the carbuncle phenomenon. A larger aspect ratio will yield a more stable flow field. The same conclusion can be found in \cite{Ohwada2013}. Chen et al.\cite{Chen2018a} also investigate the grid dependence of the shock instability through the matrix stability analysis method. They find that the instability may arise from the transverse face when it is perpendicular to the shock wave, and the shock anomalies may be cured by increasing the angle between the transverse face and the shock front.\par

It should be noted that Dumbser et al.\cite{Dumbser2004} propose the matrix stability analysis method, coupling the capacity of schemes to stably capture shocks with the eigenvalues of the stability matrix. Based on the matrix stability analysis method, the evolution of perturbation errors can be quantitatively analyzed and easily implemented by programming. Moreover, such an analysis method allows to incorporate the features of the numerical shock, effects of the computational grid, and boundary conditions \cite{Simon2018a}. As a result, the matrix stability analysis method becomes a powerful tool for judging whether a scheme is stable or not. When the stability analysis yields an unstable result for a scheme, the instability in Euler computations can always be found, although the unstable phenomenon may only appear in certain conditions\cite {Shen2014}. Consequently, the matrix stability analysis method is widely used to investigate the mechanism of the shock instabilities \cite{Dumbser2004,Chauvat2005,Chen2018a,Liu2020} and validate the shock stability of the novel shock-stable schemes \cite{Chen2018,Xie2019b,Hu2022a,Chen2018c,Chen2023,Sun2022a}. However, the matrix stability analysis method proposed in \cite{Dumbser2004} is only applicable to the first-order scheme. Since second or even high-order schemes are actually applied in practical hypersonic flow computations, the shock stability of such methods needs to be investigated in a comprehensive manner. Results from numerical experiments \cite{Liou2000,Zhang2017b} reveal that the stability of the second-order scheme may differ from the first-order case and the limiter used in the second-order scheme plays an important role in triggering the shock instability. Tu et al.\cite{Tu2014} and Jiang et al.\cite{Jiang2017} even investigate the stability of high-order schemes by various numerical experiments and reveal that high-order schemes are at a higher risk of shock instability. However, Kemm\cite{Kemm_Heuristical_2018} thinks that increasing the order offers an alternative way to stabilize the shock position by introducing more degrees of freedom to remodel the Rankine–Hugoniot condition at a captured shock. From the researches above, it can be found that there is no consensus about the influence of spatial order and despite these examples of progress, few efforts have been devoted to analyzing the shock stability of second or even high-order shock-capturing methods quantitatively. Such a situation is mainly due to the lack of an effective analytical method. In the current study, we develop the matrix stability analysis method for the finite-volume MUSCL approach, providing a quantitative analysis tool for the stability of second-order schemes. A follow-up study will present the matrix stability analysis method for high-order finite-volume schemes. The primary numerical characteristics and the underlying mechanism of shock instability for finite-volume schemes can be investigated in detail.\par

The outline of the rest of this paper is as follows. In section \ref{section 2}, governing equations of compressible flows and their related finite volume discretization are presented. In section \ref{section 3}, the planar steady shock problem is introduced and its associated instability behaviors are presented. In section \ref{section 4}, we establish the matrix stability analysis method for the second-order finite-volume method. Numerical characteristics and the underlying mechanism of shock instability for second-order finite-volume schemes are explored with the developed analysis method in section \ref{section 5}. Section \ref{section 6} contains conclusions and an outlook on future developments.\par

\section{Governing equations and finite-volume discretization}\label{section 2}

\subsection{The Euler equations}\label{subsection 2.1}
We consider a compressible flow governed by the two-dimensional Euler equations written in integral form as
\begin{equation}\label{eq Euler equations}
    \frac{\partial}{\partial t} \int_{\Omega} \mathbf{U} {\rm{d}}{\Omega} + \oint_{\partial \Omega} {\mathbf{F}} {\rm{d}}S = 0,
\end{equation}
where $\partial \Omega$ denotes boundaries of the control volume $\Omega$ and ${\mathbf{F}}$ is the flux component normal to the boundary $\partial \Omega$. The vectors of conservative variables and flux are given respectively by
\begin{equation}
	\mathbf{U} = \left[ {\begin{array}{{c}}
	\rho \\
	{\rho u}\\
	{\rho v}\\
	{\rho e}
	\end{array}} \right],
	\quad
	{\mathbf{F}} = \left[ {\begin{array}{{c}}
	\rho q \\
	{\rho q u + p n_x}\\
	{\rho q v + p n_y}\\	
	{\left(\rho e+p\right)q}
	\end{array}} \right],
\end{equation}
where $\rho$, $e$, and $p$ represent density, specific total energy, and pressure respectively, and ${\bf{u}} = \left( {u,v} \right)$ is the flow velocity. The directed velocity, $q=un_x+vn_y$, is the component of velocity acting in the $\mathbf{n}$ direction, where $\mathbf{n}={\left[ {{n_x},\;{n_y}} \right]^T}$ is the outward unit vector normal to the surface element ${\rm{d}}S$. The equation of state is in the form
\begin{equation}
	p = \left(\gamma-1\right) \rho \left[e-\frac{1}{2}\left(u^2+v^2\right)\right],
\end{equation}
where $\gamma$ is the specific heat ratio.\par

\subsection{The finite-volume MUSCL schemes}\label{subsection 2.2}

\begin{figure}[htbp]
	\centering
	\includegraphics[width=0.50\textwidth]{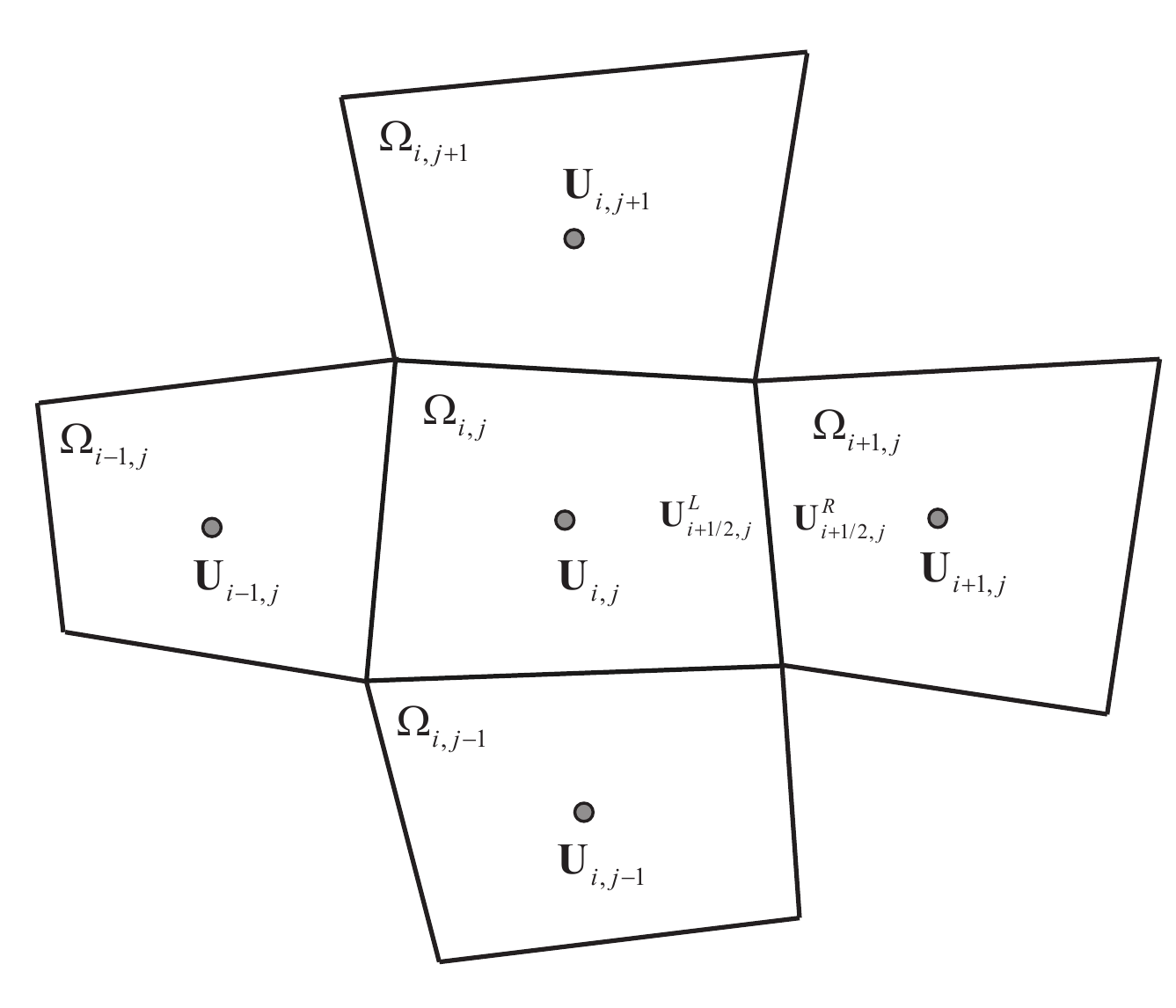}
	\caption{Possible grid configuration for two-dimensional domain in x-y space.}
	\label{fig control volume}
\end{figure}

We consider discretizing the system (\ref{eq Euler equations}) with a cell-centered finite-volume method over a two-dimensional domain subdivided into some structured quadrilateral cells. A particular control volume $ \Omega_{i,j} $ is shown in Fig.\ref{fig control volume}, and the semi-discrete finite volume scheme over it can be written as
\begin{equation}\label{eq discrete Euler equations}
    \frac{\mathrm{d} {\mathbf{U}}_{i,j}}{\mathrm{~d} t}=-\frac{1}{\left|\Omega_{i,j}\right|} \sum_{k=1}^{4} \mathcal{L}_{k} \mathbf{F}_{k},
\end{equation}
where $ {\mathbf{U}}_{i,j} $ denotes the average of \textbf{U} on $ \Omega _{i,j} $. $ \left|\Omega_{i,j}\right| $ is the volume of $ \Omega _{i,j} $ and $\mathcal{L}_{k}$ stands for the length of the cell interface. The flux $\mathbf{F}_{k}$ is the calculated numerical flux that is supposed to be constant along the individual cell interface. In the current study, approximate Riemann solvers are implemented to determine the numerical flux $\mathbf{F}_{k}$ at each cell interface.

\begin{figure}[htbp]
	\centering
	\includegraphics[width=0.7\textwidth]{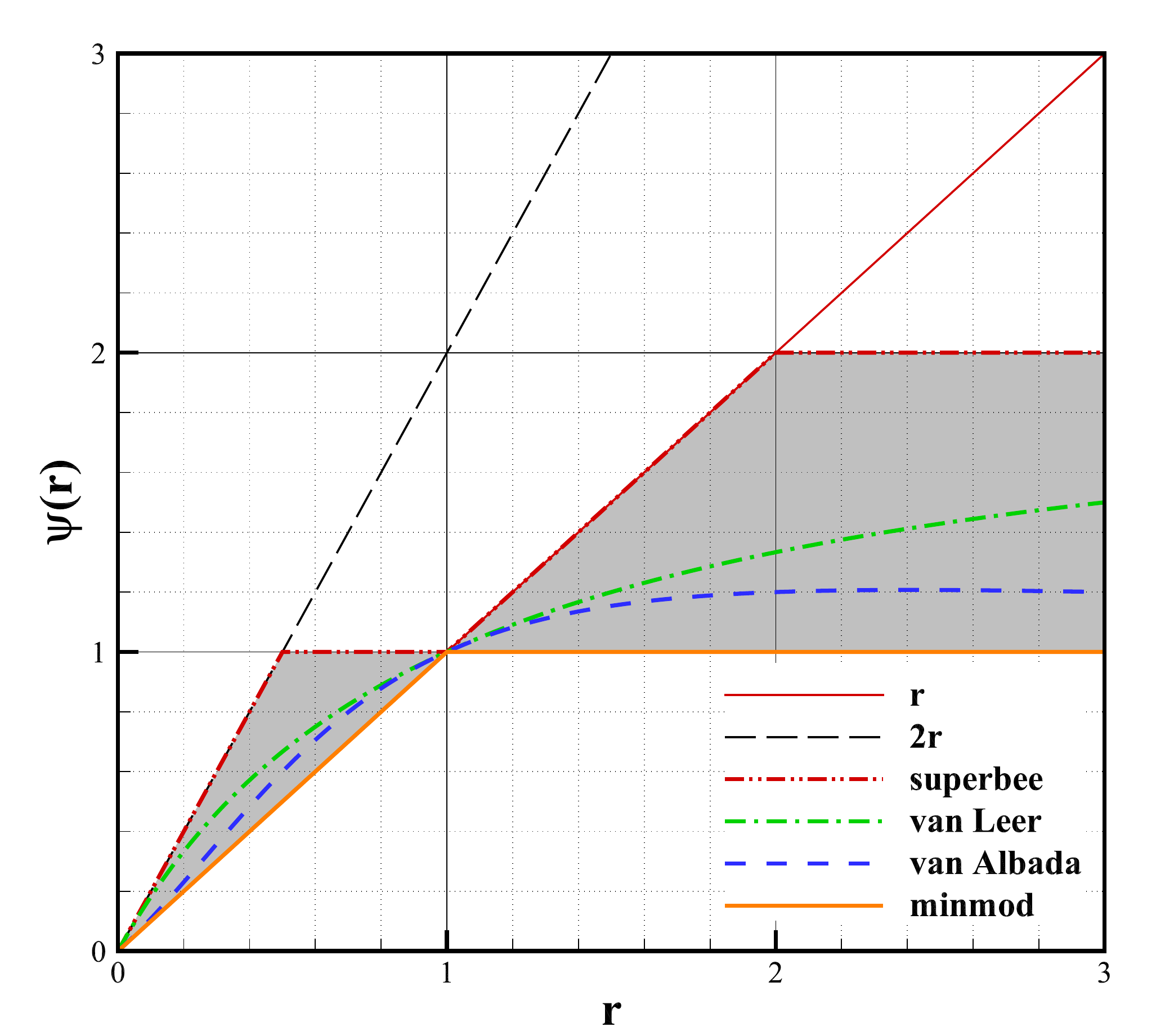}
	\caption{Second-order TVD region.}
	\label{fig limiters comparison}
\end{figure}

In our study, we devote our efforts to investigating the shock instability of second-order finite-volume schemes, particularly the MUSCL approach \cite{VanLeer1979,VanLeer2021}. As seen in Fig.\ref{fig control volume}, considering the interface between $ \Omega_{i,j} $ and $ \Omega_{i+1,j} $, we can get \cite{Blazek2015,Hemker1987}
\begin{equation}\label{eq MUSCL method}
    \begin{aligned}
        & \mathbf{U}_{i+1 / 2,j}^{L}=\mathbf{U}_{ i,j}+\frac{1}{2} \Psi_{i+1 / 2,j}^{L}\left(\mathbf{U}_{i,j}-\mathbf{U}_{i-1,j}\right) \\
        & \mathbf{U}_{i+1 / 2,j}^{R}=\mathbf{U}_{i+1,j}-\frac{1}{2} \Psi_{i+1 / 2,j}^{R}\left(\mathbf{U}_{i+2,j}-\mathbf{U}_{i+1,j}\right)
        \end{aligned},
\end{equation}
where $ \mathbf{U}_{i+1 / 2,j}^{L} $ and $ \mathbf{U}_{i+1 / 2,j}^{R} $ are the variables on the left and right sides of the interface and $ \Psi_{i+1/2,j}^{L/R} $ is called limiter. The limiter is a function of the consecutive gradients, which can be written as
\begin{equation}
	\Psi_{i+1/2,j}^{L/R}=\Psi(r_{i+1/2,j}^{L/R}),
\end{equation}
where
\begin{equation}
    \begin{aligned}
        & r_{i+1 / 2,j}^{L}=\frac{\mathbf{U}_{i+1,j}-\mathbf{U}_{i,j}}{\mathbf{U}_{i,j}-\mathbf{U}_{i-1,j}} \equiv \frac{\Delta_{+}}{\Delta_{-}} \mathbf{U}_{i,j} \\
        & r_{i+1 / 2,j}^{R}=\frac{\mathbf{U}_{i+1,j}-\mathbf{U}_{i,j}}{\mathbf{U}_{i+2,j}-\mathbf{U}_{i+1,j}} \equiv \frac{\Delta_{-}}{\Delta_{+}} \mathbf{U}_{i+1,j}
    \end{aligned},
\end{equation}
where $ \Delta_{+} $ and $ \Delta_{-} $ are the forward and backward difference operators. When $ \Delta_{-}\mathbf{U}_{i,j}=0 $ and $ \Delta_{+}\mathbf{U}_{i+1,j}=0 $, $ \Psi_{i+1/2,j}^{L}$ and $ \Psi_{i+1/2,j}^{R}$ are set to be zero in this paper, respectively. And note that the corresponding scheme becomes first-order accurate if $ \Psi_{i+1/2,j}^{L/R}=0 $. In the current work, four popular limiters including superbee \cite{Roe1985}, van Leer \cite{BramVanLeer1974}, van Albada \cite{VanAlbada1982}, and minmod \cite{Roe1986a} are employed. Their functions can be written as
\begin{equation}\label{eq limiter functions}
    \begin{aligned}
        & \Psi_{superbee} (r)=\max [0,\min (2 r, 1), \min (r, 2)]\\
        & \Psi_{van Leer} (r)=\frac{r+|r|}{1+r}\\
        & \Psi_{van Albada} (r)=\frac{r^{2}+r}{1+r^{2}}\\
        & \Psi_{minmod} (r)=\max [0, \min (r, 1)]
	\end{aligned}.
\end{equation}
As shown in Fig.\ref{fig limiters comparison}, four limiter functions all lie in the shaded region, which is the second-order TVD region \cite{Sweby1984}. Generally, from the upper boundary of the shadow region to the lower, the limiters become more dissipative \cite{Tang2020}.\par

\section{The planar steady shock instability}\label{section 3}
The carbuncle phenomenon is conventionally referred to as the distorted shock ahead of the blunt-body in the supersonic or hypersonic flow \cite{Quirk1994}. It is essential to choose a simpler test case that shares the fundamental characteristics of the blunt-body carbuncle and can be analyzed simply. In the present work, we consider the 2D steady normal shock problem. It has been well demonstrated that if a scheme yields unstable solutions for the steady normal shock problem, it will also suffer from the blunt-body carbuncle \cite{Dumbser2004,Ismail2006,Kitamura2009}.

\begin{figure}[htbp]
	\centering
	\includegraphics[width=0.4\textwidth]{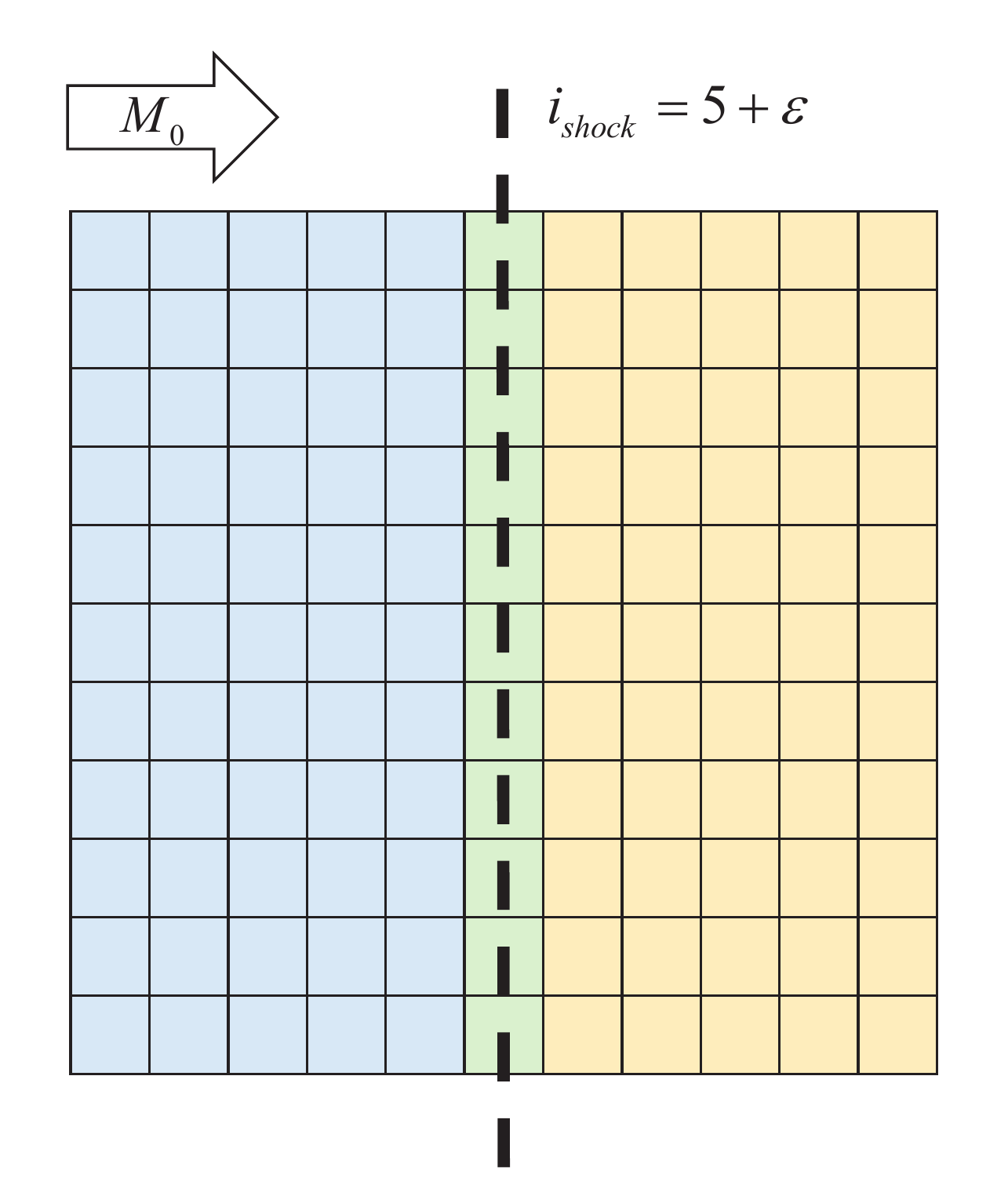}
	\caption{Computational grid for 2D steady normal shock problem.}
	\label{fig 2D grid}
\end{figure}

\subsection{Definition of the problem}\label{subsection 3.1}
Fig.\ref{fig 2D grid} shows the computational grid for the steady shock problem. The initial conditions are prescribed for left ($ L: i\leq 5 $) and right ($ R: i\geq 7 $) following the Rankine-Hugoniot conditions across the normal shock as
\begin{equation}\label{eq initial conditions}
	\mathbf{U}_{L}=\begin{pmatrix}1\\ 1\\ 0\\ \frac{1}{\gamma (\gamma-1)M_{0}^{2}}+\frac{1}{2}\end{pmatrix},\quad \mathbf{U}_{R}=\begin{pmatrix}f(M_{0})\\ 1\\ 0\\ \frac{g(M_{0})}{\gamma(\gamma-1)M_{0}^{2}}+\frac{1}{2f(M_{0})}\end{pmatrix},
\end{equation}
with
\begin{equation}
	f(M_{0})=\left(\frac{2}{(\gamma+1)M_{0}^{2}}+\frac{\gamma-1}{\gamma+1}\right)^{-1},\quad g(M_{0})=\frac{2\gamma M_{0}^{2}}{\gamma+1}-\frac{\gamma-1}{\gamma+1},
\end{equation}
where $ M_0 $ is the upstream Mach number. The intermediate states within the shock ($ M: i=6 $) are set according to the Hugoniot curve \cite{Chauvat2005}
\begin{equation}
	\begin{aligned}\rho_M&=(1-\alpha_\rho)\rho_L+\alpha_\rho\rho_R\\ u_M&=(1-\alpha_u)u_L+\alpha_u u_R\\ p_M&=(1-\alpha_p)p_L+\alpha_p p_R\end{aligned}
\end{equation}
with
\begin{equation}  \label{eq numerical shock structure}
	\begin{aligned}
		&a_{\rho}=\varepsilon\\
		&\alpha_u=1-(1-\varepsilon)\left(1+\varepsilon\frac{M_0^2-1}{1+(\gamma-1)M_0^2/2}\right)^{-1/2}\left(1+\varepsilon\frac{M_0^2-1}{1-2\gamma M_0^2/(\gamma-1)}\right)^{-1/2}\\
		&\alpha_p=\varepsilon\left[1+(1-\varepsilon)\frac{\gamma+1}{\gamma-1}\frac{M_0^2-1}{M_0^2}\right]^{-1/2}
	\end{aligned},
\end{equation}
where $ \varepsilon \in [0,1] $ is a weighting average that describes the initial state of the internal cell and is called shock position here. The boundary conditions in all directions are imposed through ghost cells. The inflow boundary conditions are set to freestream values, while in order to fix the shock at the same position, the outflow boundary conditions are determined by setting the mass flux at the ghost cells as \cite{Kitamura2009,Ismail2009}
\begin{equation}
	(\rho u)_{imax+1,j}=(\rho u)_{imax+2,j}=(\rho u)_{0}=1.
\end{equation}
Meanwhile, other values are simply extrapolated. The upper and lower boundaries are set as periodic conditions. In order to explore the stability of numerical methods for strong shock capturing, a convergent and stable solution of the 2D planar steady shock problem should be obtained for the linearized stability analysis. Here, the numerical experiment setup presented above is used to carry out numerical tests for the convergent and stable solutions. However, such a convergent and stable flow field is challenging to be obtained, considering that the numerical shock instability may occur in certain cases. In the following section, a simple and effective strategy to enforce a convergent and stable flow field is introduced.\par

\subsection{Initialization of the two-dimensional flow field for stability analysis}\label{subsection 3.2}

In the current study, the 2D flow field for stability analysis is initialized by projecting the steady flow field from 1D computation onto the 2D domain. This initialization method is also employed in \cite{Sanders1998,Dumbser2004}. In this way, when a small initial perturbation is introduced into the 2D flow field, it will be damped if the scheme is stable. Otherwise, the initial perturbation will increase, leading to the instability. As a result, we can investigate the stability of different shock-capturing methods by tracking the evolution of the perturbation error. In this paper, the small initial perturbation is set as $ \delta=10^{-7} $.\par

Unfortunately, it has been demonstrated that the shock instability or the carbuncle phenomenon also occurs in the 1D computation. When it happens, the computation fails to converge and a stable solution cannot be obtained. As a result, it is hard to initialize the 2D flow field. The 1D shock instability or carbuncle usually occurs when the Riemann solver that produces stationary discrete shocks with a single interior point is used, such as Roe under the shock position $ \varepsilon={0.1,0.2,0.3} $ \cite{Kitamura2009,Zaide2011,Zaide2012,Xie2021}. Xie et al.\cite{Xie2017} find that the shock can be stabilized by enforcing the consistency of the mass flux across the normal shock, as a result of which, the stable solution can be obtained. In this paper, this method proposed in \cite{Xie2017} is used to obtain a stable result in the 1D domain, which can be written as
\begin{equation}\label{eq modification}
	(\rho u)_{R}=(\rho u)_{L},
\end{equation}
where the subscripts $ L $ and $ R $ denote the cells just in front of and after the shock structure, respectively. The modification (\ref{eq modification}) is used in this paper for Roe Riemann solver when $ \varepsilon=0.1,0.2,0.3 $, if not mentioned specifically.\par

\begin{figure}[htbp]
	\centering
	\includegraphics[width=0.8\textwidth]{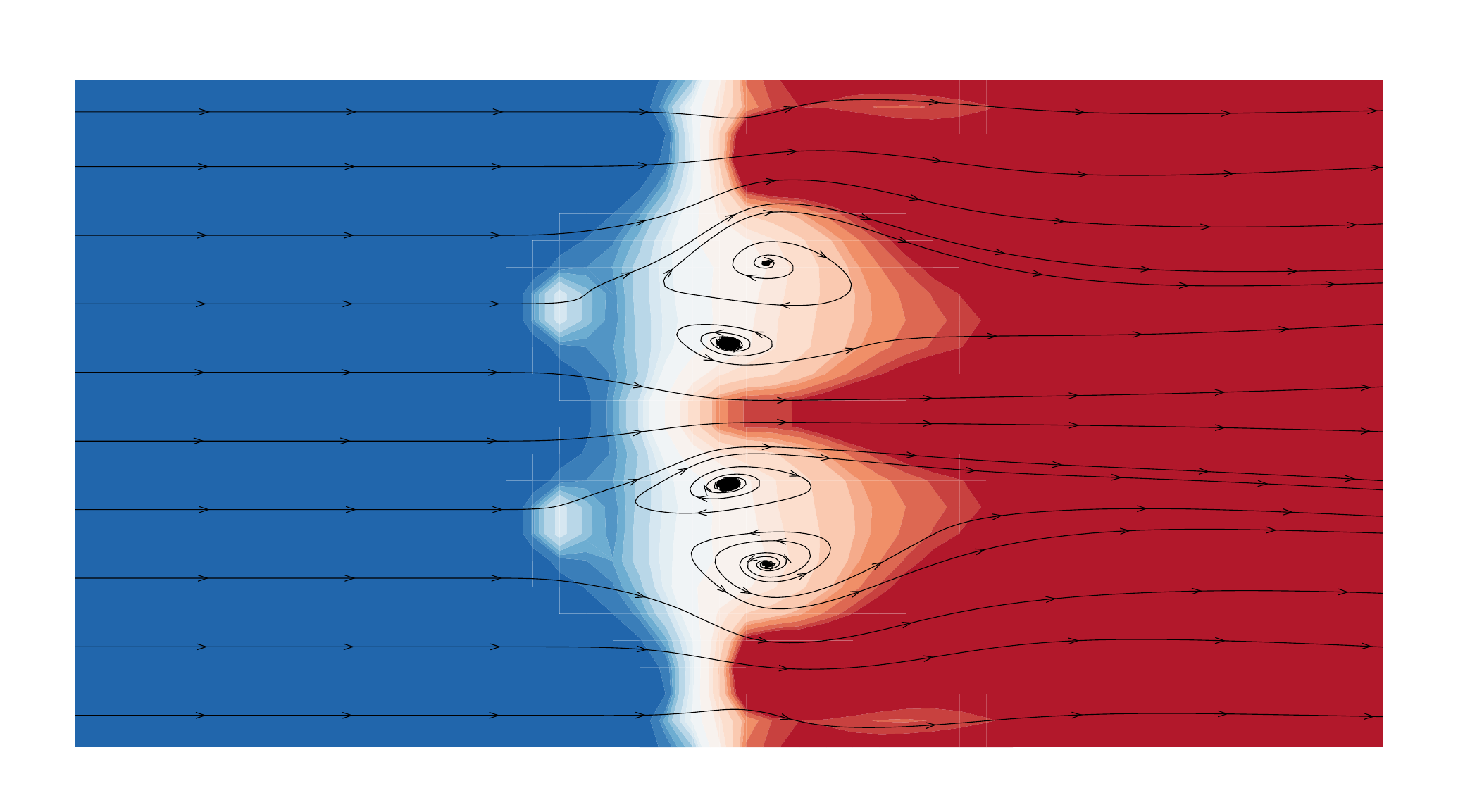}
	\caption{Density contour of the typical carbuncle phenomenon in 2D steady normal shock problem.(Grid with $50 \times 25$ cells, first-order scheme with Roe solver, $ M_0=20$, $ \varepsilon =0.1 $, and t=150.)}
	\label{fig typical carbuncle}
\end{figure}

\subsection{Carbuncle in 2D steady normal shock problem}\label{subsection 3.3}
The typical carbuncle phenomenon in 2D steady normal shock problem is shown in Fig.\ref{fig typical carbuncle}. As shown, the carbuncle phenomenon is characterized by several convex or concave wedge-shaped shock profiles distributed along the transversal direction in a staggered manner, and vortices appear behind the shock profile \cite{Xie2017}.\par

An initial perturbation of size $ 10^{-7} $ is added on each cell at the beginning of the computation, and Fig.\ref{fig error evolution} shows the evolution of the maximal perturbation error. In the current work, we take the norm $ \|v\|_{\infty}(t) $ of the transverse velocity as the indicator of perturbation error since it is known that the exact solution for mean flows is $ v = 0 $ \cite{Dumbser2004,Henderson2007}. Fig.\ref{fig error evolution} and Fig.\ref{fig typical examples} are computed with the condition of $ M_0=20 $ and $ \varepsilon=0.1 $, where the computation is more prone to be unstable. The second-order scheme with Roe solver and van Albada limiter is employed. The computational grid is 11$ \times $11 Cartesian grid. Third-order TVD Runge Kutta discretization \cite{Gottlieb1998} with $ \text{CFL} = 0.1 $ is used throughout the paper, if not mentioned specifically. As shown in Fig.\ref{fig error evolution}, according to the evolution of the perturbation error, there are three stages in the development of the shock instability \cite{Dumbser2004}:

\begin{figure}[htbp]
	\centering
	\includegraphics[width=0.7\textwidth]{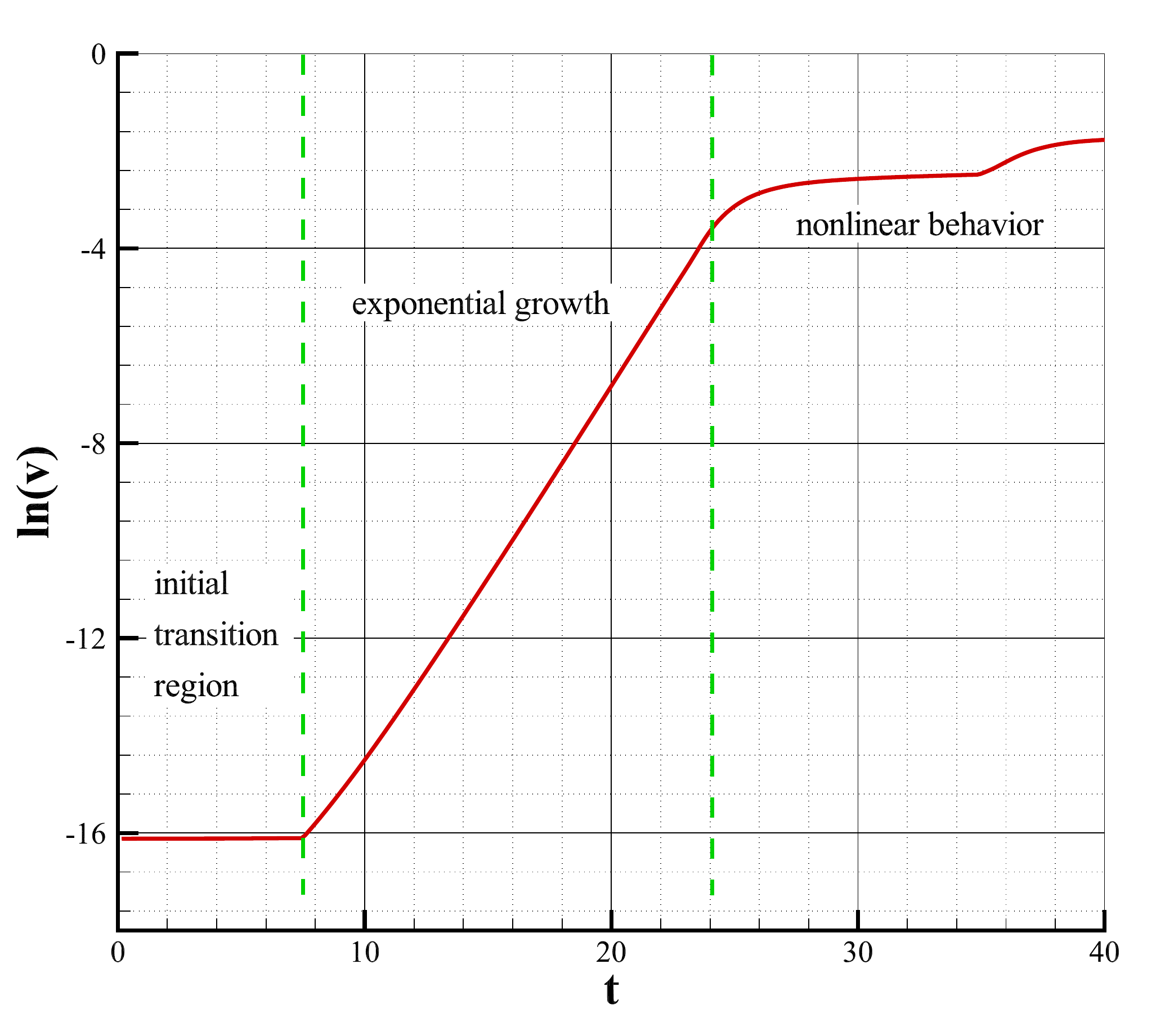}
	\caption{Evolution of the perturbation error $ \delta v $.(Grid with 11$ \times $11 cells, second-order scheme with the Roe solver and van Albada limiter, $ M_0=20 $ and $ \varepsilon =0.1 $.)}
	\label{fig error evolution}
\end{figure}

\begin{itemize}
	\item[1.]In the time interval from 0 to 7.5, which can be called the initial transition region, the perturbation error almost remains static. This might be due to the damping of stable modes \cite{Dumbser2004}. The flow field in this stage is consistent with the initial one and remains steady, as seen in Fig.\ref{fig typical examples}(a).
	\item[2.]The evolution of $ \|v\|_{\infty}(t) $ exhibits an exponential increase in the time interval from 7.5 to 24. Thus, this region is called the exponential growth region. The evolution of the perturbation error at this stage satisfies the exponential law
	\begin{equation}\label{eq exponential law}
		\|v\|_{\infty}(t)=v_{0} \mathrm{e}^{\lambda_{num}\left(t-t_{0}\right)},
	\end{equation}
	where $ v_{0} $, $ \lambda_{num} $, and $ t_{0} $ can be easily obtained from Fig.\ref{fig error evolution}, and $ \lambda_{num} $ is called the temporal error growth rate in the present work. The flow field becomes unstable during the exponential growth stage, as shown in Fig.\ref{fig typical examples}(b) and Fig.\ref{fig typical examples}(c), and the shock profile starts to deform.
	\item[3.]After time 24, the perturbation error is sufficiently large, nonlinear effects become significant, and the exponential growth is minimized. As seen in Fig.\ref{fig error evolution}, the error may fluctuate at a high level during this stage. The flow field, meanwhile, is unphysical. The instability is more severe, as shown in Fig.\ref{fig typical examples}(d), and the shock is severely distorted.
\end{itemize}

\begin{figure}[htbp]
	\centering
	\subfigure[t=7]{
	\begin{minipage}[t]{0.46\linewidth}
	\centering
	\includegraphics[width=0.95\textwidth]{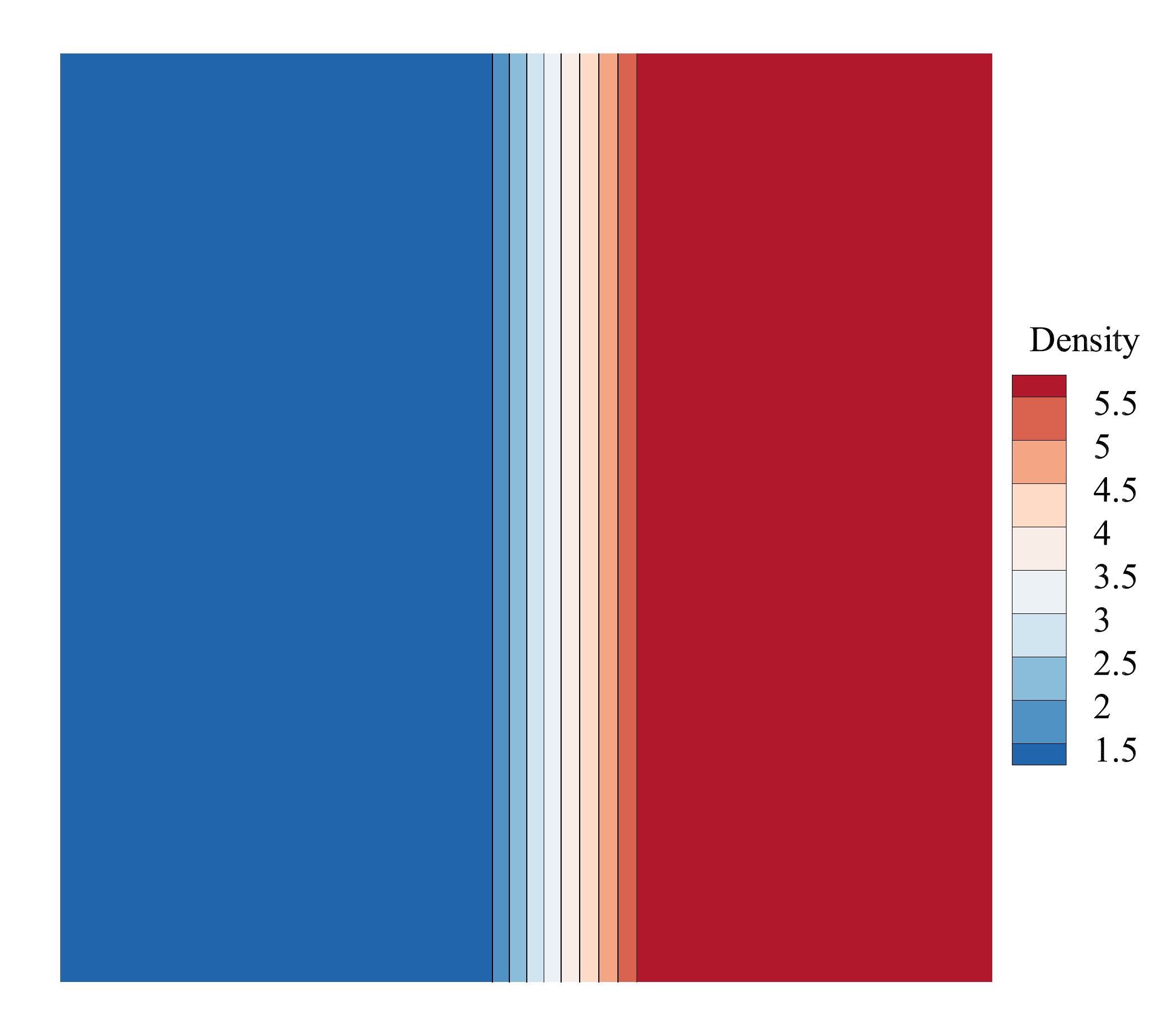}
	\end{minipage}
	}
	\subfigure[t=15]{
	\begin{minipage}[t]{0.46\linewidth}
	\centering
	\includegraphics[width=0.95\textwidth]{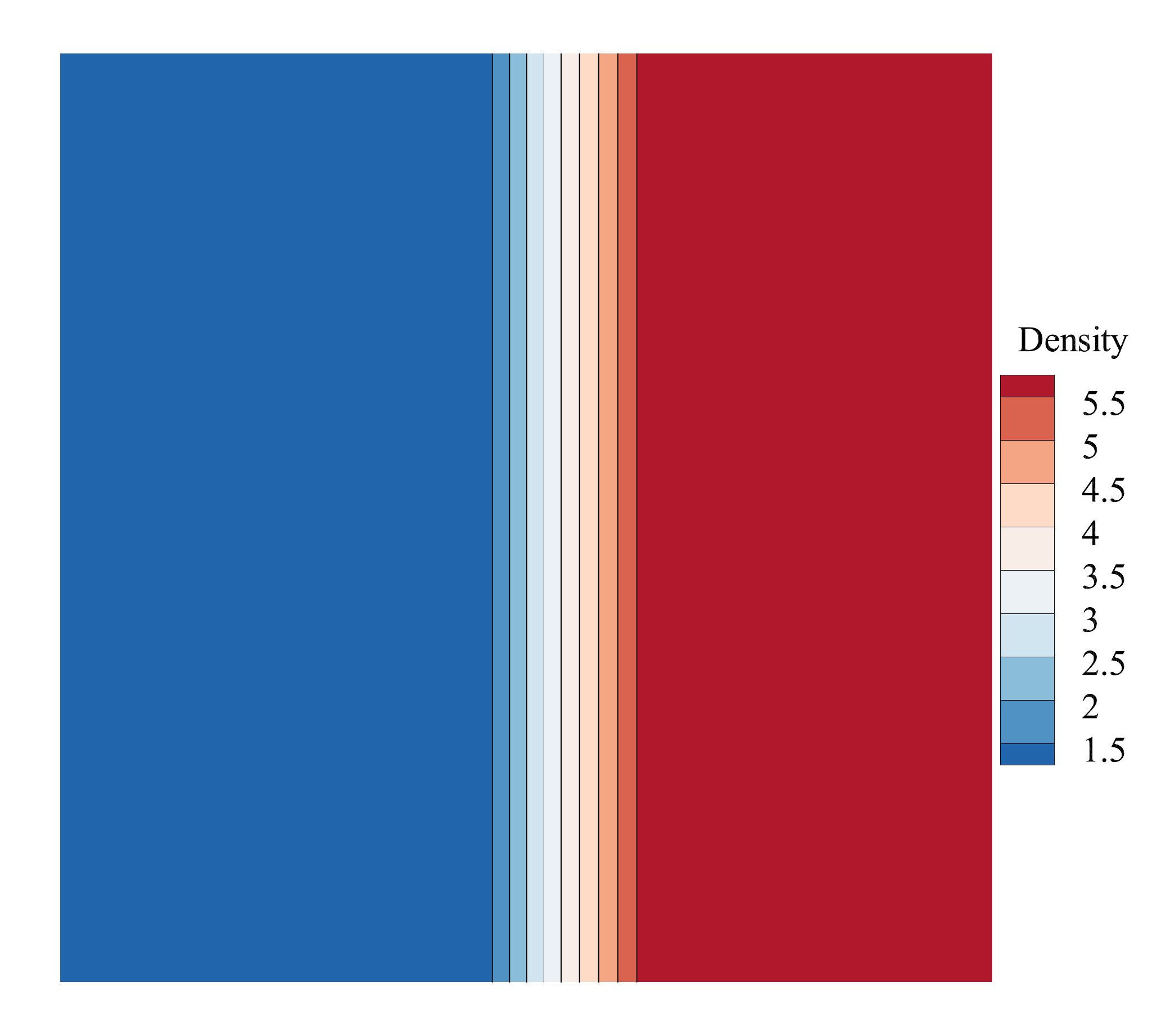}
	\end{minipage}
	}

	\subfigure[t=25]{
	\begin{minipage}[t]{0.46\linewidth}
	\centering
	\includegraphics[width=0.95\textwidth]{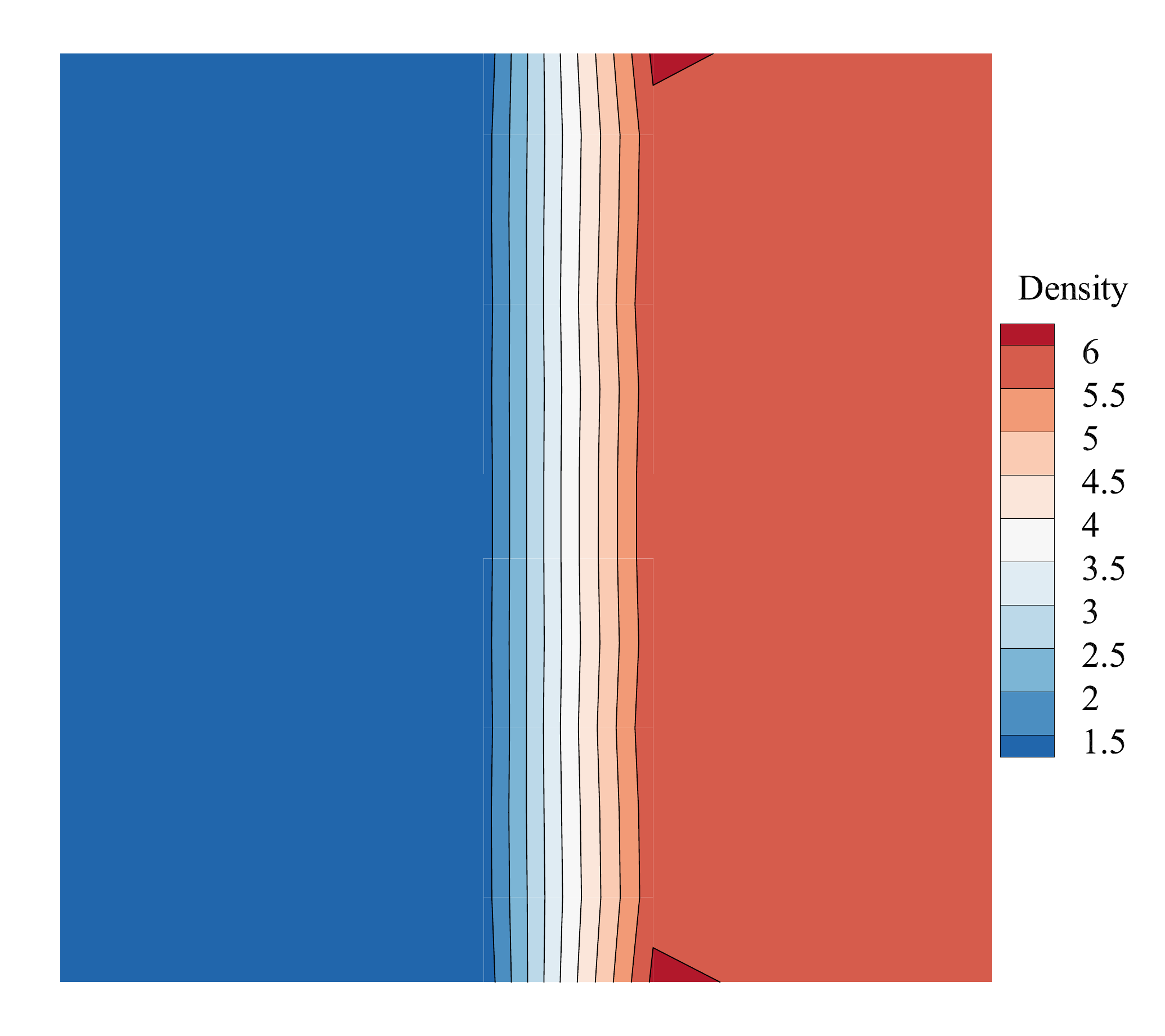}
	\end{minipage}
	}
	\subfigure[t=80]{
	\begin{minipage}[t]{0.46\linewidth}
	\centering
	\includegraphics[width=0.95\textwidth]{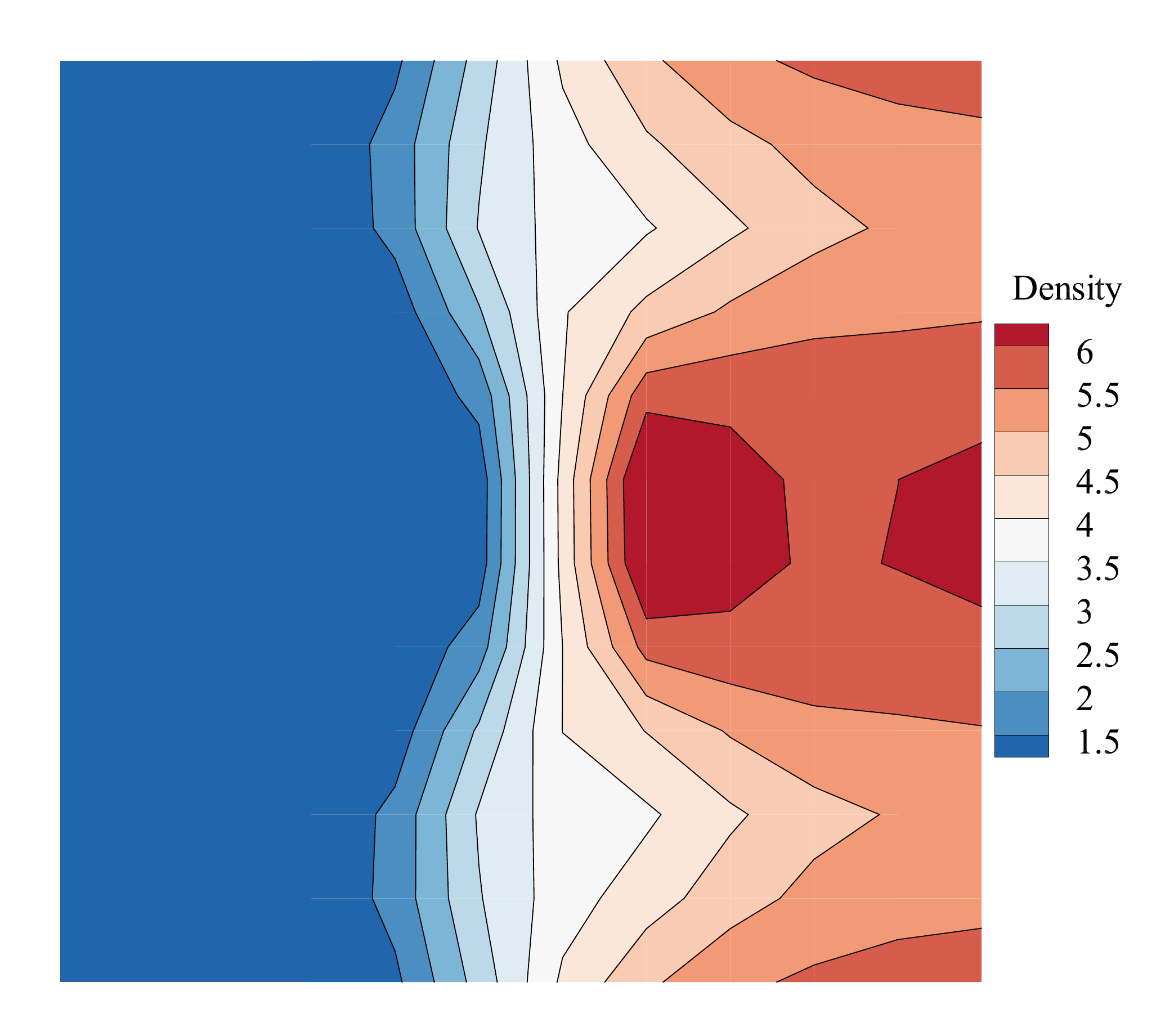}
	\end{minipage}
	}

	\centering
	\caption{Evolution of the flow fields.(Grid with 11$ \times $11 cells, second-order scheme with the Roe solver and van Albada limiter, $ M_0=20 $ and $ \varepsilon =0.1 $.)}\label{fig typical examples}
\end{figure}

As shown in Fig.\ref{fig error evolution}, it can be found that the exponential growth stage plays a vital role in the evolution of the perturbation error, determining whether the error will increase or decrease, further stable or unstable, and how quickly it will develop towards instability. Therefore, in the current study we concentrate on the exponential growth stage and relate the shock instability to the temporal error growth rate $ \lambda_{num} $.\par

\begin{figure}[htbp]
	\centering
	\subfigure[first-order scheme]{
	\begin{minipage}[t]{0.46\linewidth}
	\centering
	\includegraphics[width=0.95\textwidth]{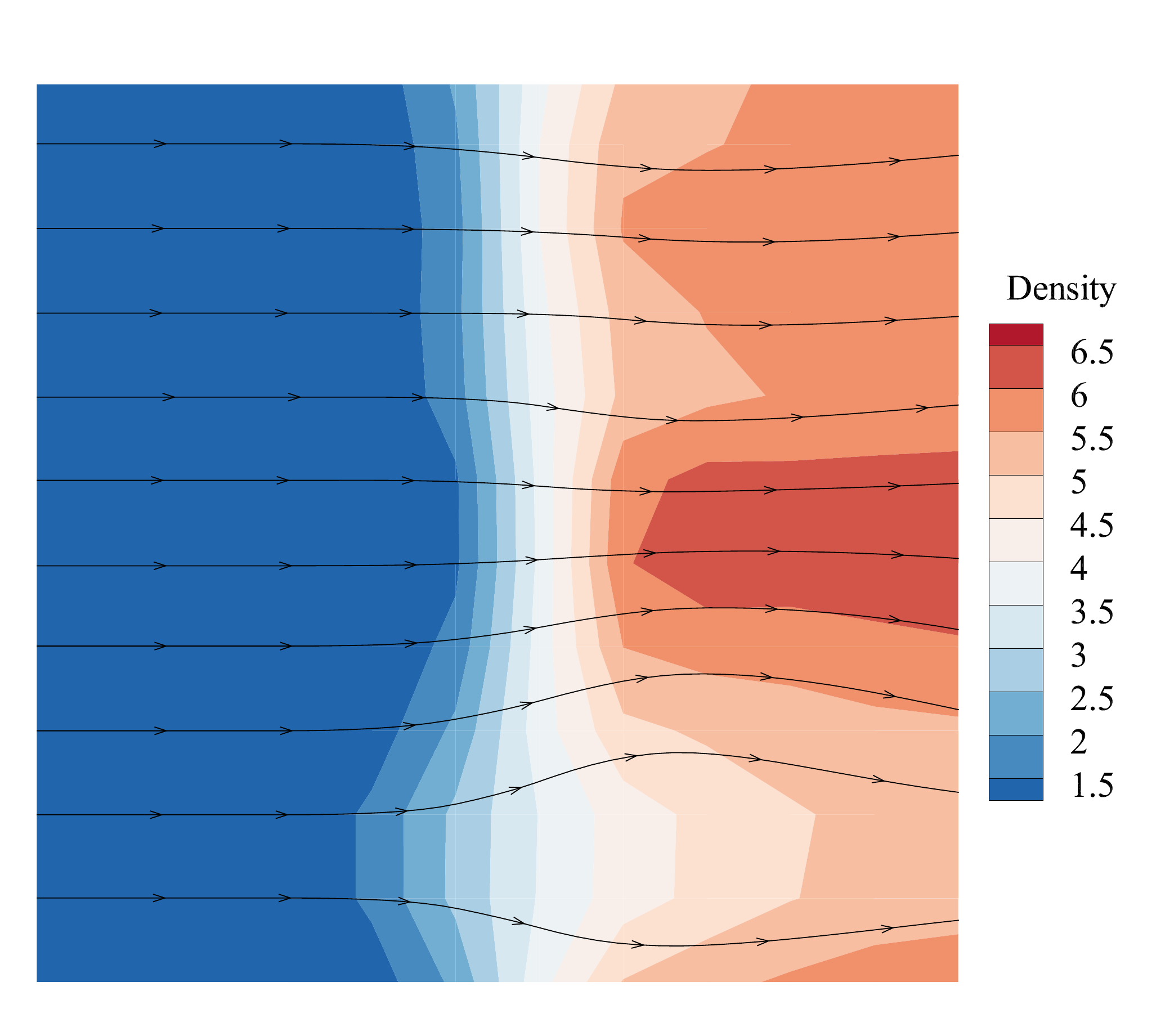}
	\end{minipage}
	}
	\subfigure[second-order scheme]{
	\begin{minipage}[t]{0.46\linewidth}
	\centering
	\includegraphics[width=0.95\textwidth]{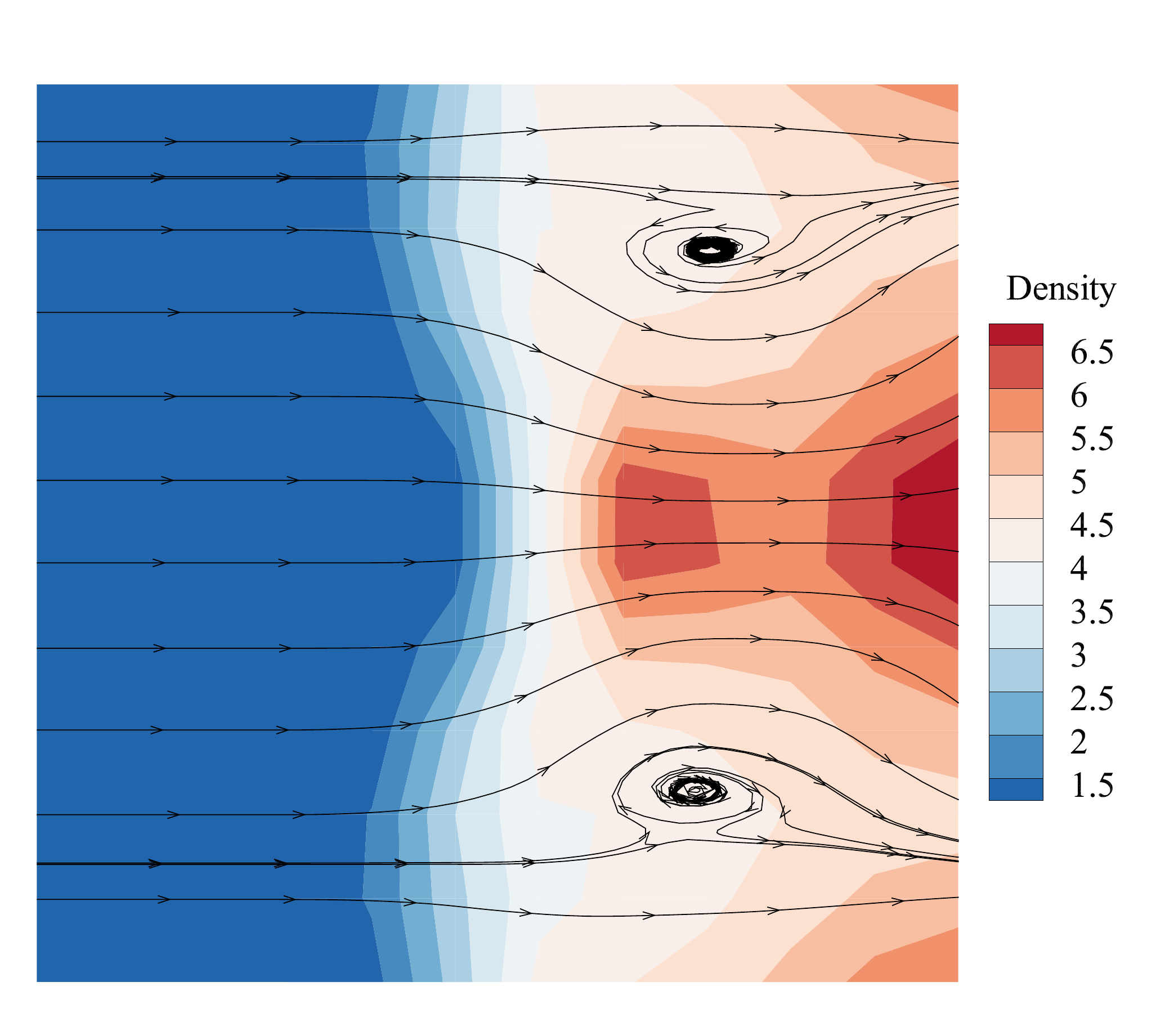}
	\end{minipage}
	}
	\centering
	\caption{Density contours of first and second-order schemes.(Grid with 11$ \times $11 cells, Roe solver, van Albada limiter used in the second-order scheme, $ M_0=20 $, $ \varepsilon =0.1 $, and t=100.)}\label{fig comparison of first- and second-order schemes}
\end{figure}

\begin{figure}[htbp]
	\centering
	\includegraphics[width=0.69\textwidth]{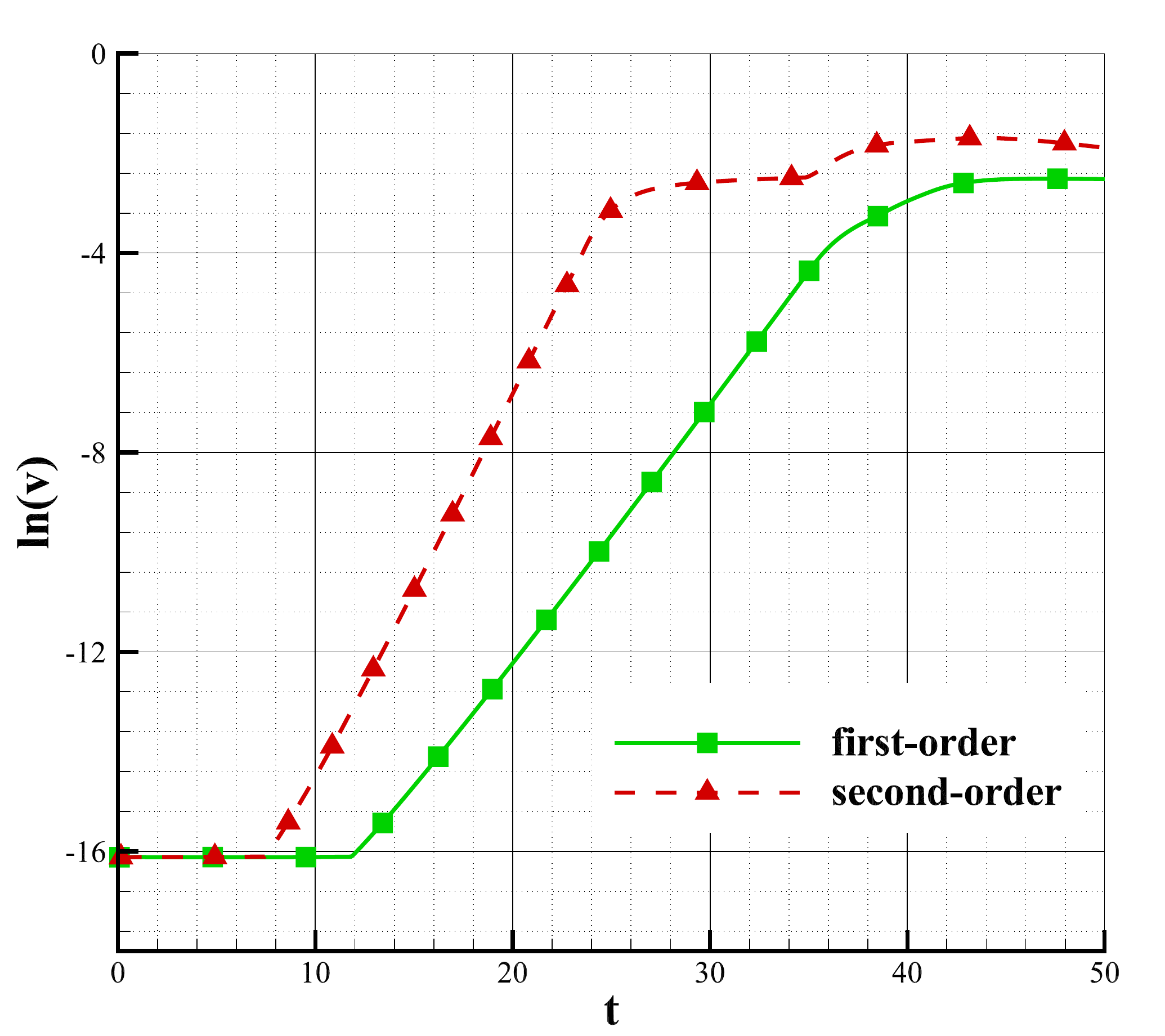}
	\caption{Evolution of the perturbation errors for first and second-order schemes.(Grid with 11$ \times $11 cells, Roe solver, van Albada limiter used in the second-order scheme, $ M_0=20$, $ \varepsilon =0.1 $, and t=100 (only shows $ t<50 $ for clarity).)}
	\label{fig comparison of the errors of first- and second-order schemes}
\end{figure}

Fig.\ref{fig comparison of first- and second-order schemes} shows the density contours computed by first and second-order schemes at the same computing time. As shown, the shock profile computed by the second-order scheme distorts more severely than the first-order case. And there are vortices appearing behind the shock in the flow field computed by the second-order scheme. Fig.\ref{fig comparison of the errors of first- and second-order schemes} is the comparison of the perturbation error, which also shows that the perturbation error produced by the second-order scheme steps into the exponential growth region earlier and increases more quickly. It is also larger at the nonlinear region than that of the first-order scheme. Both Fig.\ref{fig comparison of first- and second-order schemes} and Fig.\ref{fig comparison of the errors of first- and second-order schemes} suggest that the shock instability of the second-order scheme is more severe than the first-order scheme. And in the current study, the shock instability of the second-order scheme will be studied in more detail.\par

\section{The matrix stability analysis method for the MUSCL schemes}\label{section 4}

\subsection{The eatablishment of the matrix stability analysis method}\label{subsection 4.1}
The matrix stability analysis method proposed by Dumbser et al.\cite{Dumbser2004} has been well demonstrated to be a useful tool for assessing the shock instability problem of shock-capturing methods. However, such a stability analysis method is only applicable to first-order schemes, which limits its application in analyzing the stability of second and high-order schemes. In the current study, we extend this method for second-order finite-volume schemes, which is presented as follows.\par

For the stability analysis of a steady field, we assume that
\begin{equation}\label{eq variables decomposition}
	\mathbf{U}_{i,j}=\mathbf{U}_{i,j}^{0}+\delta \mathbf{U}_{i,j},
\end{equation}
where $ \mathbf{U}_{i,j}^{0} $ represents the steady mean value and $ \delta \mathbf{U}_{i,j} $ is small numerical random perturbation. Assume that the random perturbation has no effect on the limiter function. Thus, the limiter is the function of the steady mean value
\begin{equation}\label{eq assumption of limiter}
	\Psi_{i+1 / 2, j}^{L / R}=\Psi\left(r_{i+1 / 2, j}^{L / R, 0}\right), \quad
	\text { with }
	\left\{\begin{aligned}
		r_{i+1 / 2, j}^{L,0}&=\frac{\mathbf{U}_{i+1, j}-\mathbf{U}_{i, j}}{\mathbf{U}_{i, j}-\mathbf{U}_{i-1, j}}=\frac{\mathbf{U}_{i+1, j}^{0}-\mathbf{U}_{i, j}^{0}}{\mathbf{U}_{i, j}^{0}-\mathbf{U}_{i-1, j}^{0}} \\
		r_{i+1 / 2, j}^{R,0}&=\frac{\mathbf{U}_{i+1, j}-\mathbf{U}_{i, j}}{\mathbf{U}_{i+2, j}-\mathbf{U}_{i+1, j}}=\frac{\mathbf{U}_{i+1, j}^{0}-\mathbf{U}_{i, j}^{0}}{\mathbf{U}_{i+2, j}^{0}-\mathbf{U}_{i+1, j}^{0}}
	\end{aligned}\right..
\end{equation}
As a consequence, the matrix stability analysis method can also be applicable for non-differentiable limiters, such as mimod. And the results in section \ref{section 4} and \ref{section 5} show that such assumption is appropriate. Substituting (\ref{eq variables decomposition}) into (\ref{eq MUSCL method}), we can get
\begin{equation}\label{eq detailed second-order variables decomposition}
	\begin{aligned}
		\mathbf{U}_{i+1/2,j}^{L}&=\left(\mathbf{E}+\frac{1}{2} \Psi_{i+1/2,j}^{L}\right) \mathbf{U}_{i,j}^{0}  -\frac{1}{2} \Psi_{i+1/2,j}^{L} \mathbf{U}_{i-1,j}^{0}\\
		&+\left(\mathbf{E}+\frac{1}{2} \Psi_{i+1/2,j}^{L}\right) \delta \mathbf{U}_{i,j}  -\frac{1}{2} \Psi_{i+1/2,j}^{L} \delta \mathbf{U}_{i-1,j}\\
		\mathbf{U}_{i+1/2,j}^{R}&=\left(\mathbf{E}+\frac{1}{2} \Psi_{i+1/2,j}^{R}\right) \mathbf{U}_{i+1,j}^{0}-\frac{1}{2} \Psi_{i+1/2,j}^{R} \mathbf{U}_{i+2,j}^{0}\\
		&+\left(\mathbf{E}+\frac{1}{2} \Psi_{i+1/2,j}^{R}\right) \delta \mathbf{U}_{i+1,j}-\frac{1}{2} \Psi_{i+1/2,j}^{R} \delta \mathbf{U}_{i+2,j}
	\end{aligned},
\end{equation}
where $ \mathbf{U}_{i+1 / 2,j}^{L} $ and $ \mathbf{U}_{i+1 / 2,j}^{R} $ are the variables on the left and right sides of the interface between $ \Omega_{i,j} $ and $ \Omega_{i+1,j} $. \textbf{E} denotes the $ 4 \times 4 $ unit matrix. Since the numerical flux $ \mathbf{F}_{i+1/2,j} $ is the function of $ \mathbf{U}_{i+1/2,j}^{L} $ and $ \mathbf{U}_{i+1/2,j}^{R} $, $ \mathbf{F}_{i+1/2,j} $ can be linearized around the steady mean value as
\begin{equation}
	\begin{aligned}\label{eq detailed second-order flux linearized}
		\mathbf{F}_{i+1/2,j}&=\mathbf{F}_{i+1/2,j}\left( \mathbf{U}^{L}_{i+1/2,j},\mathbf{U}^{R}_{i+1/2,j}\right)\\
		& =\mathbf{F}_{i+1/2,j}\left( \mathbf{U}^{L,0}_{i+1/2,j},\mathbf{U}^{R,0}_{i+1/2,j}\right)+\frac{\partial \mathbf{F}_{i+1/2,j}}{\partial \mathbf{U}_{i+1/2,j}^{L,0}}\delta\mathbf{U}_{i+1/2,j}^{L}+\frac{\partial \mathbf{F}_{i+1/2,j}}{\partial \mathbf{U}_{i+1/2,j}^{R,0}}\delta\mathbf{U}_{i+1/2,j}^{R}\\
		& = \mathbf{F}_{i+1/2,j}\left( \mathbf{U}^{L,0}_{i+1/2,j},\mathbf{U}^{R,0}_{i+1/2,j}\right)\\
		&-\frac{1}{2} \frac{\partial \mathbf{F}_{i+1/2,j}}{\partial \mathbf{U}_{i+1/2,j}^{L,0}} \Psi_{i+1/2,j}^{L} \delta \mathbf{U}_{i-1,j}+\frac{\partial \mathbf{F}_{i+1/2,j}}{\partial \mathbf{U}_{i+1/2,j}^{L,0}} \left(\mathbf{E}+\frac{1}{2} \Psi_{i+1/2,j}^{L}\right) \delta \mathbf{U}_{i,j} \\
		&+\frac{\partial \mathbf{F}_{i+1/2,j}}{\partial \mathbf{U}_{i+1/2,j}^{R,0}} \left(\mathbf{E}+\frac{1}{2} \Psi_{i+1/2,j}^{R}\right) \delta \mathbf{U}_{i+1,j}-\frac{1}{2} \frac{\partial \mathbf{F}_{i+1/2,j}}{\partial \mathbf{U}_{i+1/2,j}^{R,0}} \Psi_{i+1/2,j}^{R} \delta \mathbf{U}_{i+2,j}
	\end{aligned}.
\end{equation}
To make the expression brief, the following variables are introduced
\begin{equation}
	\begin{aligned}
		&\eta_{i+1/2,j}^{L/R}=\frac{1}{2} \frac{\partial \mathbf{F}_{i+1/2,j}}{\partial \mathbf{U}_{i+1/2,j}^{L/R,0}}\Psi_{i+1/2,j}^{L/R}\\
		&\beta_{i+1/2,j}^{L/R}=\frac{\partial \mathbf{F}_{i+1/2,j}}{\partial \mathbf{U}_{i+1/2,j}^{L/R,0}}\left(\mathbf{E}+\frac{1}{2} \Psi_{i+1/2,j}^{L/R}\right)
	\end{aligned}.
\end{equation}
Then (\ref{eq detailed second-order flux linearized}) can be written as
\begin{equation}
	\begin{aligned}\label{eq simplify second-order flux linearized}
		\mathbf{F}_{i+1/2,j}& = \mathbf{F}_{i+1/2,j}\left( \mathbf{U}^{L,0}_{i+1/2,j},\mathbf{U}^{R,0}_{i+1/2,j}\right)\\
		&-\eta_{i+1/2,j}^{L} \delta \mathbf{U}_{i-1,j}+\beta_{i+1/2,j}^{L} \delta \mathbf{U}_{i,j}+\beta_{i+1/2,j}^{R} \delta \mathbf{U}_{i+1,j}-\eta_{i+1/2,j}^{R}  \delta \mathbf{U}_{i+2,j}
	\end{aligned}.
\end{equation}
Substituting (\ref{eq variables decomposition}) and (\ref{eq simplify second-order flux linearized}) into (\ref{eq discrete Euler equations}) and taking into account that the mean-field is steady, finally, we finally get the linear error evolution mode
\begin{equation}
	\begin{aligned}\label{eq second-order linear error evolution model for conservative variables}
		\frac{\mathrm{d} \delta \mathbf{U}_{i, j}}{\mathrm{dt}}=&-\left(\xi _{i+1 / 2, j}^{L} +\xi _{i, j+1 / 2}^{L} +\xi _{i-1 / 2, j}^{R}+\xi _{i, j-1 / 2}^{R}\right) \delta \mathbf{U}_{i, j} \\
		&-\left(\xi _{i+1 / 2, j}^{R}-\mu  _{i-1 / 2, j}^{R}\right) \delta \mathbf{U}_{i+1, j} -\left(\xi _{i, j+1 / 2}^{R}-\mu  _{i, j-1 / 2}^{R}\right) \delta \mathbf{U}_{i, j+1} \\
		&-\left(\xi _{i-1 / 2, j}^{L}-\mu  _{i+1 / 2, j}^{L}\right) \delta \mathbf{U}_{i-1, j} -\left(\xi _{i, j-1 / 2}^{L}-\mu  _{i, j+1 / 2}^{L}\right) \delta \mathbf{U}_{i, j-1} \\
		&+\mu  _{i+1 / 2, j}^{R} \delta \mathbf{U}_{i+2, j}+\mu  _{i, j+1 / 2}^{R} \delta \mathbf{U}_{i, j+2}+\mu  _{i-1 / 2, j}^{L} \delta \mathbf{U}_{i-2, j}+\mu  _{i, j-1 / 2}^{L} \delta \mathbf{U}_{i, j-2}
	\end{aligned},
\end{equation}
where
\begin{equation}\label{eq temporary variable}
	\xi _{i+1 / 2, j}^{L/R}=\frac{\mathcal{L}_{i+1/2,j}}{|\Omega _{i,j}|}\beta_{i+1/2,j}^{L/R}\quad , \quad \mu _{i+1 / 2, j}^{L/R}=\frac{\mathcal{L}_{i+1/2,j}}{|\Omega _{i,j}|}\eta _{i+1/2,j}^{L/R}.
\end{equation}
It can be found from (\ref{eq second-order linear error evolution model for conservative variables}) that the evolution of perturbation error in $ \Omega_{i,j} $ is influenced by the errors in cell $ \Omega_{i,j} $ itself, four neighbors, and four sub-adjacent cells. However, for the first-order scheme, only the errors in $ \Omega_{i,j} $ and four neighbors will affect the evolution of the error in $ \Omega_{i,j} $ \cite{Dumbser2004}. The difference in the linear error evolution model suggests that the shock stability of the first and second-order schemes may be significantly different.\par

Equation (\ref{eq second-order linear error evolution model for conservative variables}) holds for all cells in the computational domain and we finally get the error evolution of all cells in the computational domain
\begin{equation}\label{eq linear error evolution model for all cells}
	\frac{\mathrm{d}}{\mathrm{d} t}\left(\begin{array}{c}
		\delta \mathbf{U}_{1,1} \\
		\vdots \\
		\delta \mathbf{U}_{imax,jmax}
		\end{array}\right)=\mathbf{S} \cdot\left(\begin{array}{c}
		\delta \mathbf{U}_{1,1} \\
		\vdots \\
		\delta \mathbf{U}_{imax,jmax}
		\end{array}\right),
\end{equation}
where \textbf{S} is called the stability matrix in the present study. When considering the evolution of the initial error only, (\ref{eq linear error evolution model for all cells}) can be solved analytically. The solution is
\begin{equation}\label{eq solution of linear error evolution model}
	\left(\begin{array}{c}
		\delta \mathbf{U}_{1,1} \\
		\vdots \\
		\delta \mathbf{U}_{imax,jmax}
		\end{array}\right)(t)=\mathrm{e}^{\mathbf{S} t} \cdot\left(\begin{array}{c}
		\delta \mathbf{U}_{1,1} \\
		\vdots \\
		\delta \mathbf{U}_{imax,jmax}
		\end{array}\right)_{t=0},
\end{equation}
and will remain bounded if the maximal real part of the eigenvalues of \textbf{S} is negative. So the stability criterion is
\begin{equation}\label{eq stability criterion}
	\max (\operatorname{Re}(\lambda(\mathbf{S}))) \leq 0.
\end{equation}

One should note that equation (\ref{eq detailed second-order flux linearized}) is accurate only when the numerical flux is differentiable at the mean value, which is not always holding, for example, when the Roe solver is employed and the shock is exactly between two cells ($ \varepsilon =0 \enspace \text{or} \enspace 1 $). In the current work, the cases that have numerical shock structure ($ 0< \varepsilon< 1 $) are analyzed. As a result, the numerical flux calculated by the solvers used in this paper is differentiable at the mean value. The gradients of the numerical flux functions {such as $ \dfrac{\partial \mathbf{F}_{i+1/2,j}}{\partial\mathbf{U}_{i+1/2,j}^{L,0}} $} can be calculated by the centered difference approximation \cite{Shen2014}
\begin{equation}
	\dfrac{\partial \mathbf{F}_{i+1/2,j}}{\left(\partial \mathbf{U}_{i+1/2,j}^{L,0}\right)_k}=\dfrac{\mathbf{F}_{i+1/2,j}(\mathbf{U}_{i+1/2,j}^{L,0}+\delta\mathbf{I}_k,\mathbf{U}_{i+1/2,j}^{R,0})-\mathbf{F}_{i+1/2,j}(\mathbf{U}_{i+1/2,j}^{L,0}-\delta\mathbf{I}_k,\mathbf{U}_{i+1/2,j}^{R,0})}{2\delta},
\end{equation}
where $ \mathbf{I}_k $ is unit vector of which the $ kth $ component is 1, and $ \delta =10^{-7} $ \cite{Dumbser2004,Shen2014}.

\subsection{The primitive variables form of the matrix stability analysis method}\label{subsection 4.2}
It is known that the second-order MUSCL approach (\ref{eq MUSCL method}) can be performed not only in conservative variables but also in primitive variables. The reconstruction of primitive variables is widely used in flow simulation, especially in hypersonic flow simulation, due to its lower cost and broad applicability \cite{Rider1993,Zanotti2016}. Therefore, the matrix stability analysis method in the form of primitive variables is pursued in the current study. The numerical flux $ \mathbf{F}_{i+1/2,j} $ is also the function of the primitive variables
\begin{equation}
	\mathbf{F}_{i+1/2,j}=\mathbf{F}_{i+1/2,j}\left( \mathbf{W}_{i+1/2,j}^{L},\mathbf{W}_{i+1/2,j}^{R} \right),
\end{equation}
where $ \mathbf{W}=(\rho,u,v,p)^T $ is the vector of primitive variables. Then the flux function $\mathbf{F}_{i+1/2,j}$ can be linearized as
\begin{equation}
	\begin{aligned}\label{eq second-order flux linearized for primitive variables}
		\mathbf{F}_{i+1/2,j}& = \mathbf{F}_{i+1/2,j}\left( \mathbf{W}^{L,0}_{i+1/2,j},\mathbf{W}^{R,0}_{i+1/2,j}\right)\\
		&-\eta_{i+1/2,j}^{L} \delta \mathbf{W}_{i-1,j}+\beta_{i+1/2,j}^{L} \delta \mathbf{W}_{i,j}+\beta_{i+1/2,j}^{R} \delta \mathbf{W}_{i+1,j}-\eta_{i+1/2,j}^{R}  \delta \mathbf{W}_{i+2,j}
	\end{aligned},
\end{equation}
where the variables $ \eta_{i+1/2,j}^{L/R} $ and $ \beta_{i+1/2,j}^{L/R} $ in (\ref{eq second-order flux linearized for primitive variables}) are as follows
\begin{equation}
	\begin{aligned}
		&\eta_{i+1/2,j}^{L/R}=\frac{1}{2} \frac{\partial \mathbf{F}_{i+1/2,j}}{\partial \mathbf{W}_{i+1/2,j}^{L/R,0}}\Psi_{i+1/2,j}^{L/R}\\
		&\beta_{i+1/2,j}^{L/R}=\frac{\partial \mathbf{F}_{i+1/2,j}}{\partial \mathbf{W}_{i+1/2,j}^{L/R,0}}\left(\mathbf{E}+\frac{1}{2} \Psi_{i+1/2,j}^{L/R}\right)
	\end{aligned}.
\end{equation}
Substituting (\ref{eq variables decomposition}) and (\ref{eq second-order flux linearized for primitive variables}) into (\ref{eq discrete Euler equations}), we can get
\begin{equation}
	\begin{aligned}\label{eq second-order linear error evolution model for primitive variables}
		\frac{\mathrm{d} \delta \mathbf{W}_{i, j}}{\mathrm{dt}}=&\left(\frac{\mathrm{d} \mathbf{U}_{i,j}}{\mathrm{d} {\mathbf{W}_{i,j}}}\right)^{-1}\left[-\left(\xi _{i+1 / 2, j}^{L} +\xi _{i, j+1 / 2}^{L} +\xi _{i-1 / 2, j}^{R}+\xi _{i, j-1 / 2}^{R}\right) \delta \mathbf{W}_{i, j}\right. \\
		&-\left(\xi _{i+1 / 2, j}^{R}-\mu  _{i-1 / 2, j}^{R}\right) \delta \mathbf{W}_{i+1, j} -\left(\xi _{i, j+1 / 2}^{R}-\mu  _{i, j-1 / 2}^{R}\right) \delta \mathbf{W}_{i, j+1} \\
		&-\left(\xi _{i-1 / 2, j}^{L}-\mu  _{i+1 / 2, j}^{L}\right) \delta \mathbf{W}_{i-1, j} -\left(\xi _{i, j-1 / 2}^{L}-\mu  _{i, j+1 / 2}^{L}\right) \delta \mathbf{W}_{i, j-1} \\
		&\left.+\mu  _{i+1 / 2, j}^{R} \delta \mathbf{W}_{i+2, j}+\mu  _{i, j+1 / 2}^{R} \delta \mathbf{W}_{i, j+2}+\mu  _{i-1 / 2, j}^{L} \delta \mathbf{W}_{i-2, j}+\mu  _{i, j-1 / 2}^{L} \delta \mathbf{W}_{i, j-2} \right]
	\end{aligned},
\end{equation}
where $ \xi $ and $ \mu $ are computed by (\ref{eq temporary variable}). And $\frac{\mathrm{d} \mathbf{U}_{i,j}}{\mathrm{d} {\mathbf{W}_{i,j}}}$ is the transformation matrix between the conservative variables and the primitive variables. It can be written as
\begin{equation}\label{eq the transfer matrix}
	\frac{\mathrm{d} \mathbf{U}_{i,j}}{\mathrm{d} \mathbf{W}_{i,j}}=\left[\begin{array}{cccc}
		1 & 0 & 0 & 0 \\
		u & \rho & 0 & 0 \\
		v & 0 & \rho & 0 \\
		(u^{2}+v^{2}) / 2 & \rho u & \rho v & 1 /(\gamma-1)
	\end{array}\right]_{i,j}.
\end{equation}
Equation (\ref{eq second-order linear error evolution model for primitive variables}) is the primitive variables form of linear error evolution model in $ \Omega _{i,j} $. (\ref{eq linear error evolution model for all cells}) can also be written in primitive variables form as follows
\begin{equation}\label{eq linear error evolution model of primitive variables for all cells}
	\frac{\mathrm{d}}{\mathrm{d} t}\left(\begin{array}{c}
		\delta \mathbf{W}_{1,1} \\
		\vdots \\
		\delta \mathbf{W}_{imax,jmax}
		\end{array}\right)=\mathbf{S} \cdot\left(\begin{array}{c}
		\delta \mathbf{W}_{1,1} \\
		\vdots \\
		\delta \mathbf{W}_{imax,jmax}
		\end{array}\right).
\end{equation}
Considering the evolution of initial errors only, the solution of (\ref{eq linear error evolution model of primitive variables for all cells}) is
\begin{equation}
	\left(\begin{array}{c}
		\delta \mathbf{W}_{1,1} \\
		\vdots \\
		\delta \mathbf{W}_{imax,jmax}
		\end{array}\right)(t)=\mathrm{e}^{\mathbf{S} t} \cdot\left(\begin{array}{c}
		\delta \mathbf{W}_{1,1} \\
		\vdots \\
		\delta \mathbf{W}_{imax,jmax}
		\end{array}\right)_{t=0}.
\end{equation}
And the stability criterion can be seen in (\ref{eq stability criterion}).\par

\begin{figure}[htbp]
	\centering
	\includegraphics[width=0.69\textwidth]{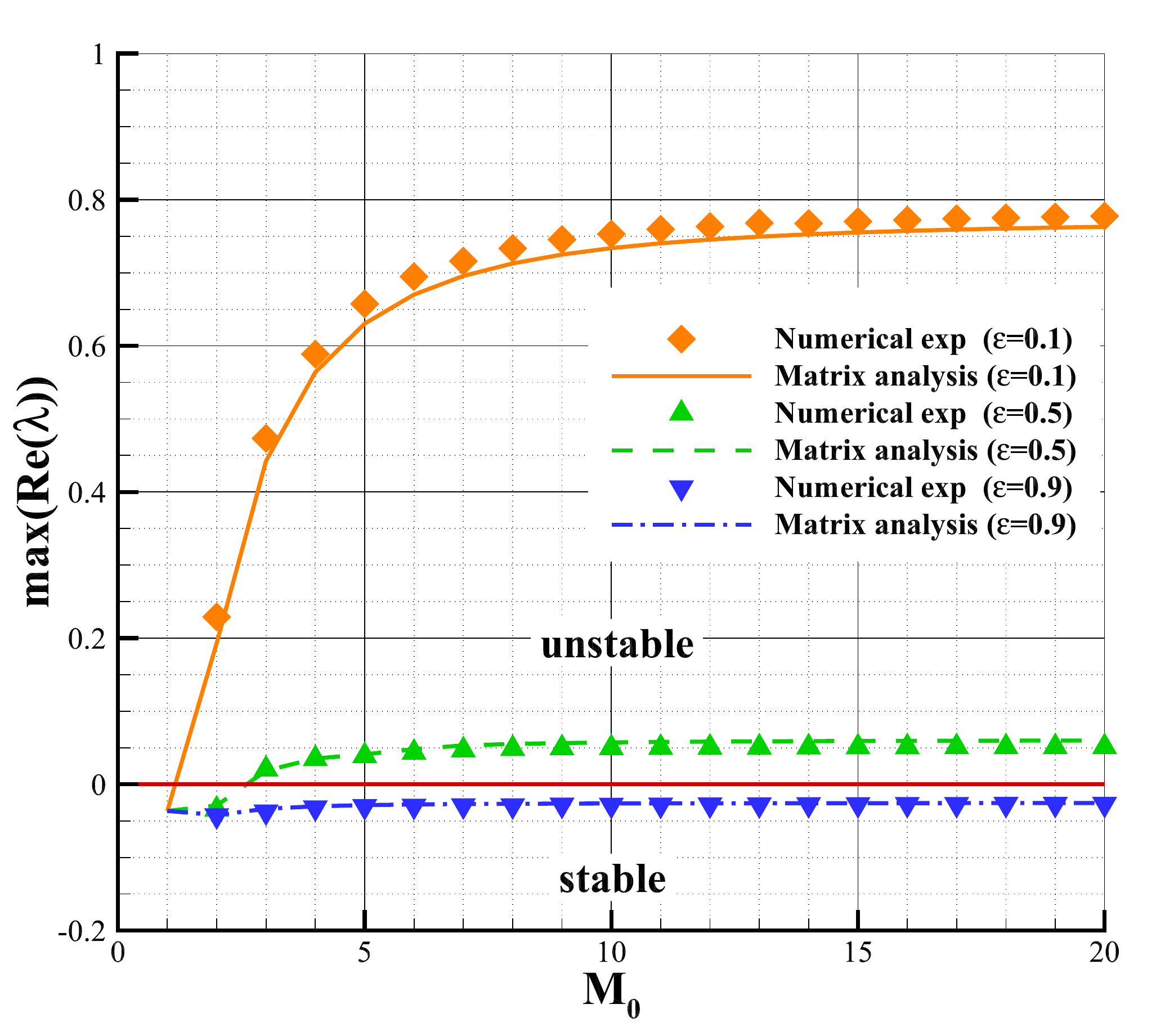}
	\caption{Quantitative validation of the stability theory proposed in this paper.(Grid with 11$ \times $11 cells, Roe solver and van Albada limiter, $ M_0={1,2,\cdots,20} $ and $ \varepsilon =0.1,0.5,0.9 $.)}
	\label{fig quantitative validation}
\end{figure}

\subsection{Quantitative validation of the stability theory}\label{subsection 4.3}
In section \ref{subsection 4.1} and \ref{subsection 4.2}, we introduce the matrix stability analysis method for the second-order MUSCL scheme. By analyzing (\ref{eq solution of linear error evolution model}), we can find the maximal real part of the eigenvalues can describe the exponential growth of the perturbation error, which is $ \lambda_{num} $ in (\ref{eq exponential law}). Since $ \lambda_{num} $ can be obtained easily by numerical experiments, the reliability of the matrix stability analysis method can be quantitatively verified by comparing $ \text{max}(\text{Re}(\lambda)) $ and $ \lambda_{num} $.\par

The second-order scheme with Roe solver and van Albada limiter is used to perform the quantitative validation. The conditions are $ M_0 = 1,2,\cdots ,20 $ and $ \varepsilon =0.1,0.5,0.9 $. The computational grid is $ 11 \times 11 $ Cartesian grid. Fig.\ref{fig quantitative validation} provides the comparison between $ \text{max}(\text{Re}(\lambda)) $ and $ \lambda_{num} $, which shows a good agreement between them in both unstable region and stable region. One should note that good agreements can also be obtained by other second-order schemes with different Riemann solvers, limiter functions, and computational conditions. That confirms the reliability of the matrix stability analysis method developed in the current study.\par

\section{Stability analysis of the MUSCL schemes}\label{section 5}

\begin{figure}[htbp]
	\centering
	\subfigure[first-order Roe]{
	\begin{minipage}[t]{0.46\linewidth}
	\centering
	\includegraphics[width=0.9\textwidth]{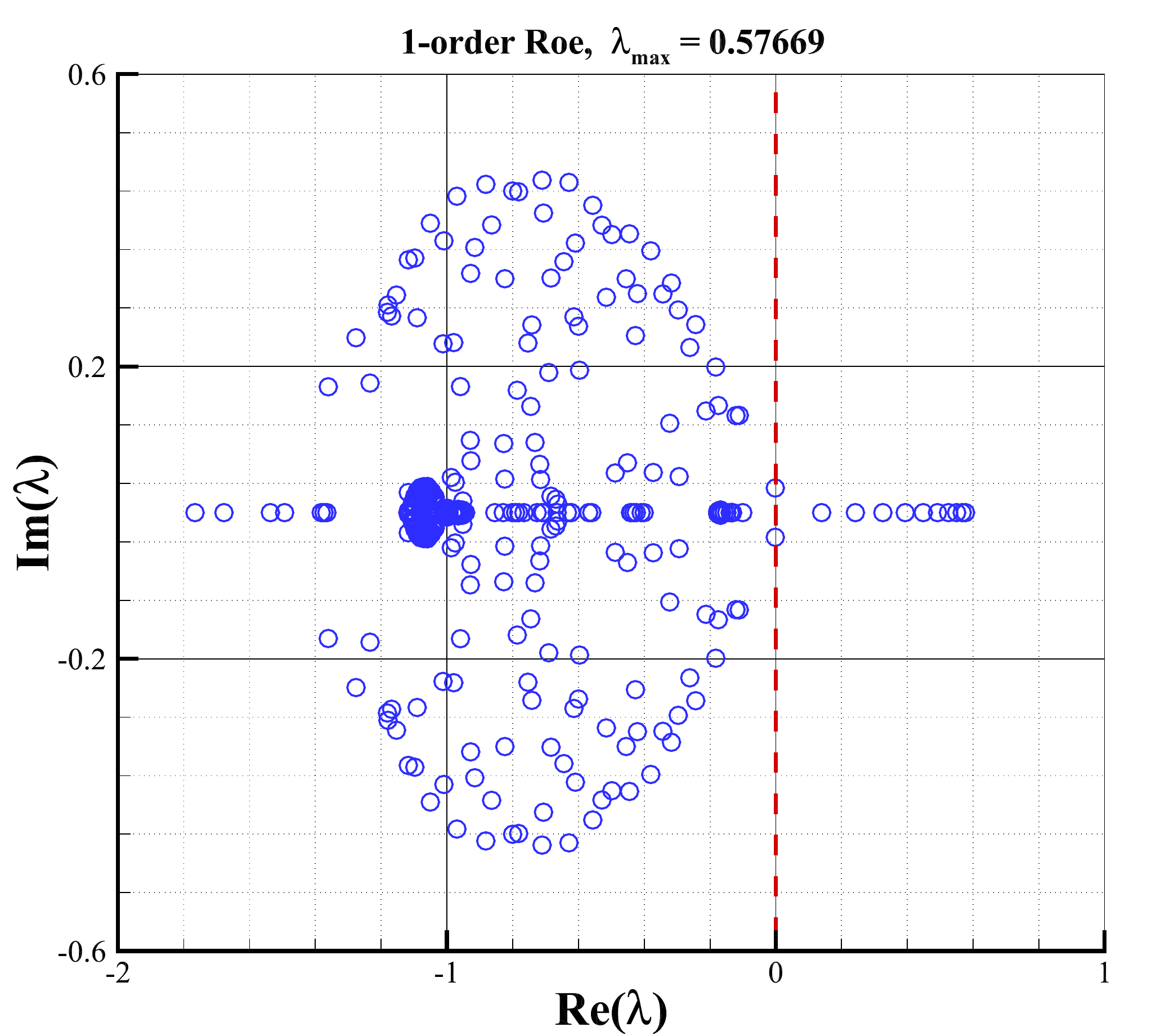}
	\end{minipage}
	}
	\subfigure[second-order Roe]{
	\begin{minipage}[t]{0.46\linewidth}
	\centering
	\includegraphics[width=0.9\textwidth]{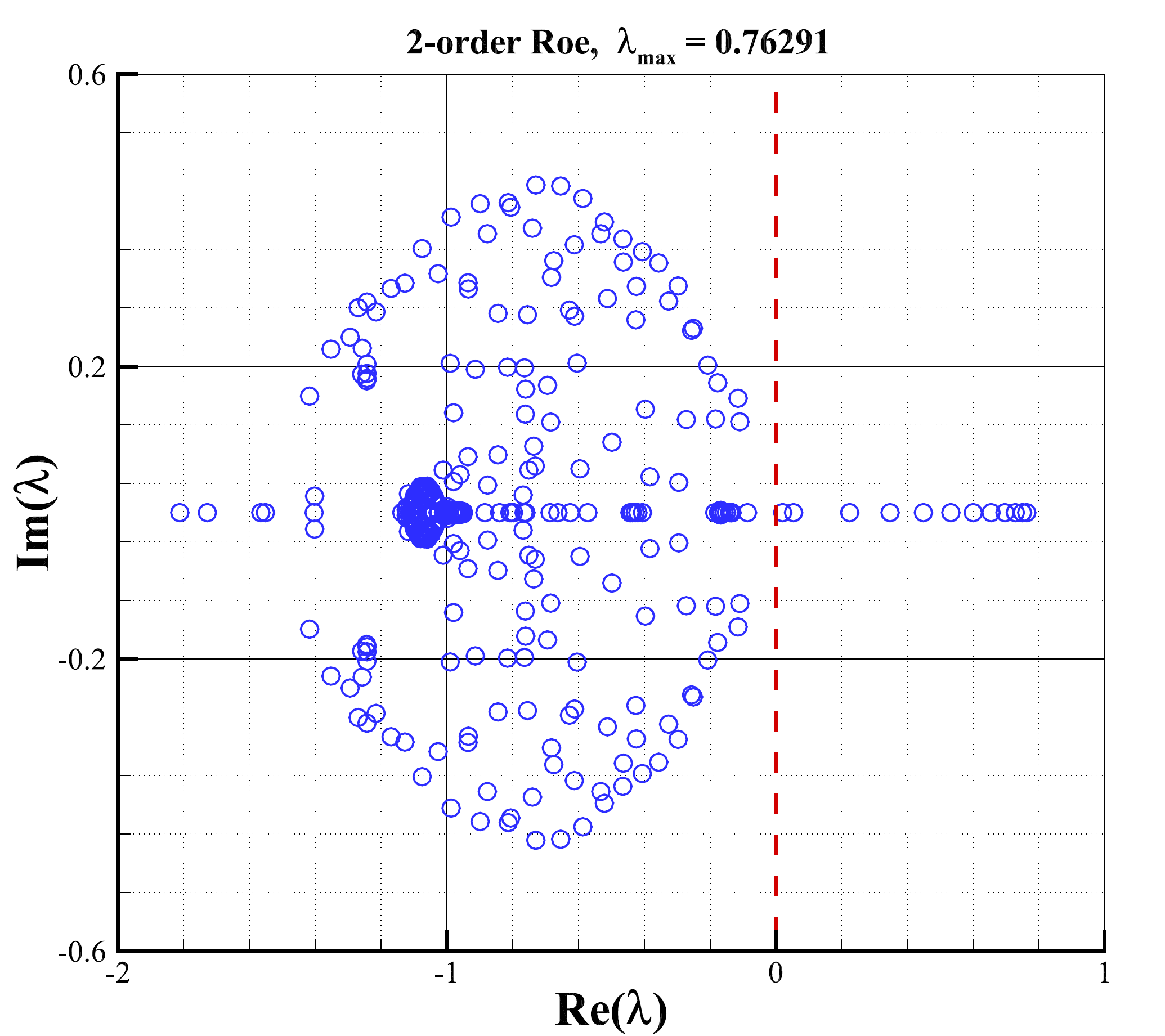}
	\end{minipage}
	}

	\subfigure[first-order HLLC]{
	\begin{minipage}[t]{0.46\linewidth}
	\centering
	\includegraphics[width=0.9\textwidth]{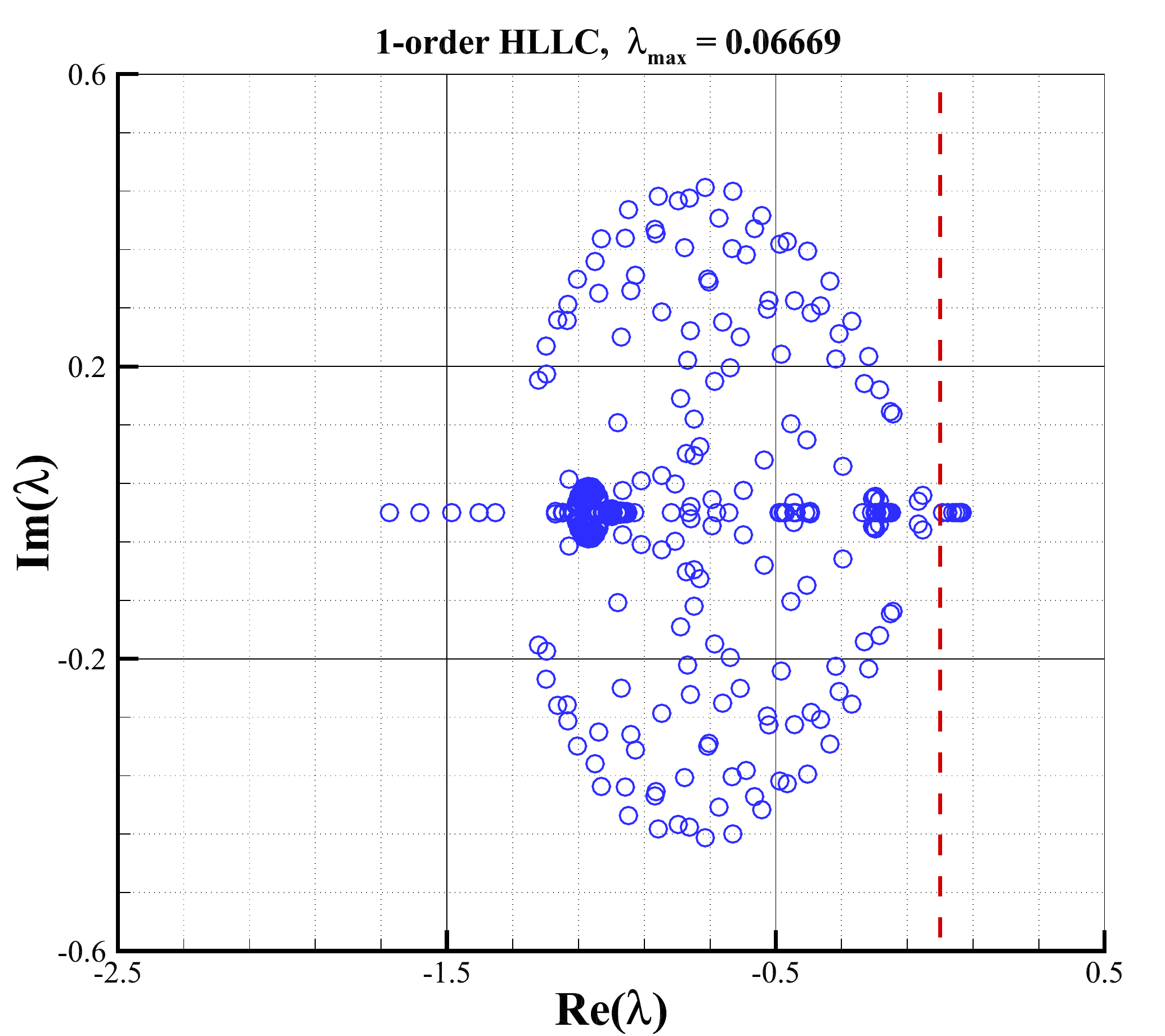}
	\end{minipage}
	}
	\subfigure[second-order HLLC]{
	\begin{minipage}[t]{0.46\linewidth}
	\centering
	\includegraphics[width=0.9\textwidth]{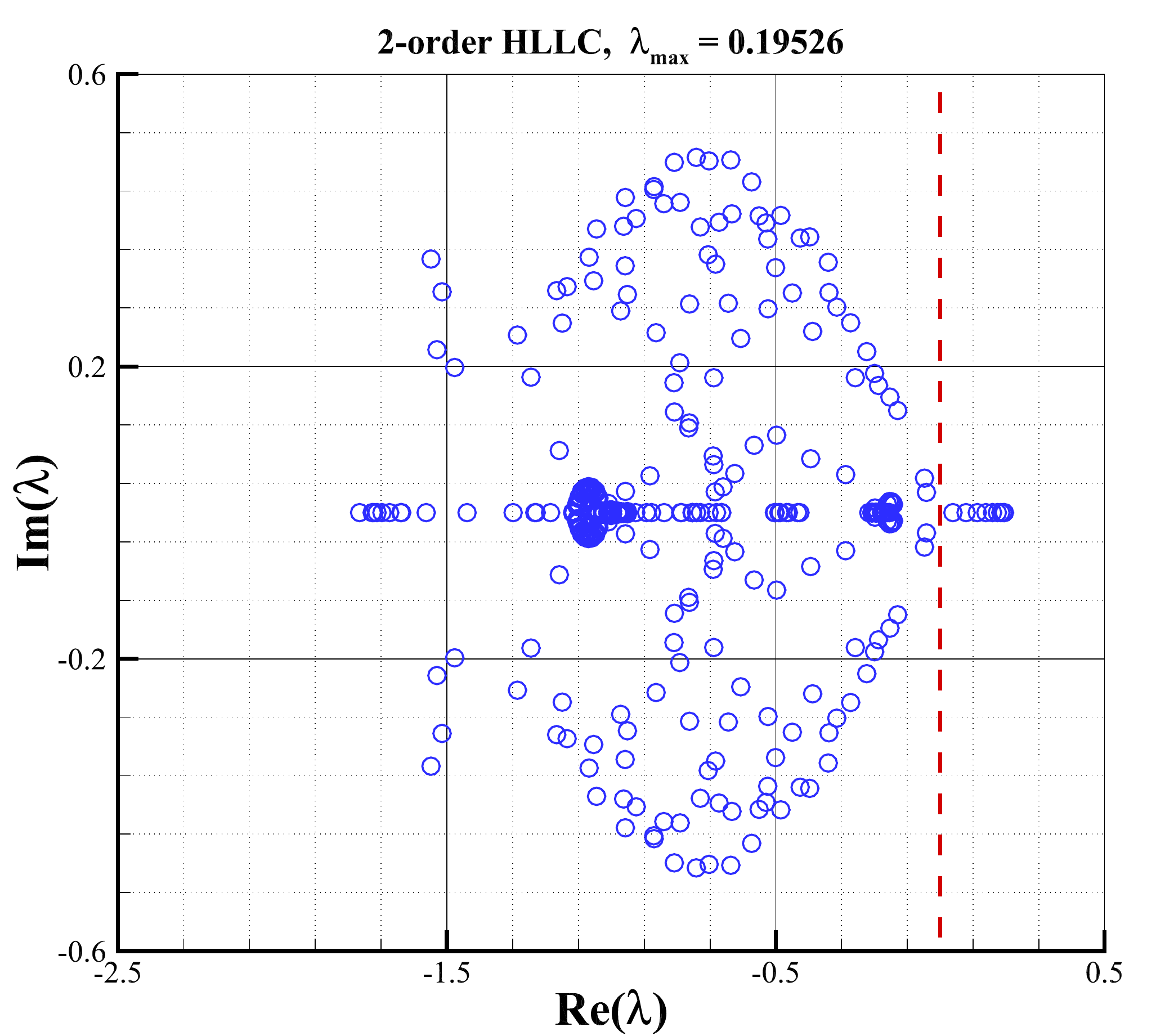}
	\end{minipage}
	}

	\subfigure[first-order HLL]{
	\begin{minipage}[t]{0.46\linewidth}
	\centering
	\includegraphics[width=0.9\textwidth]{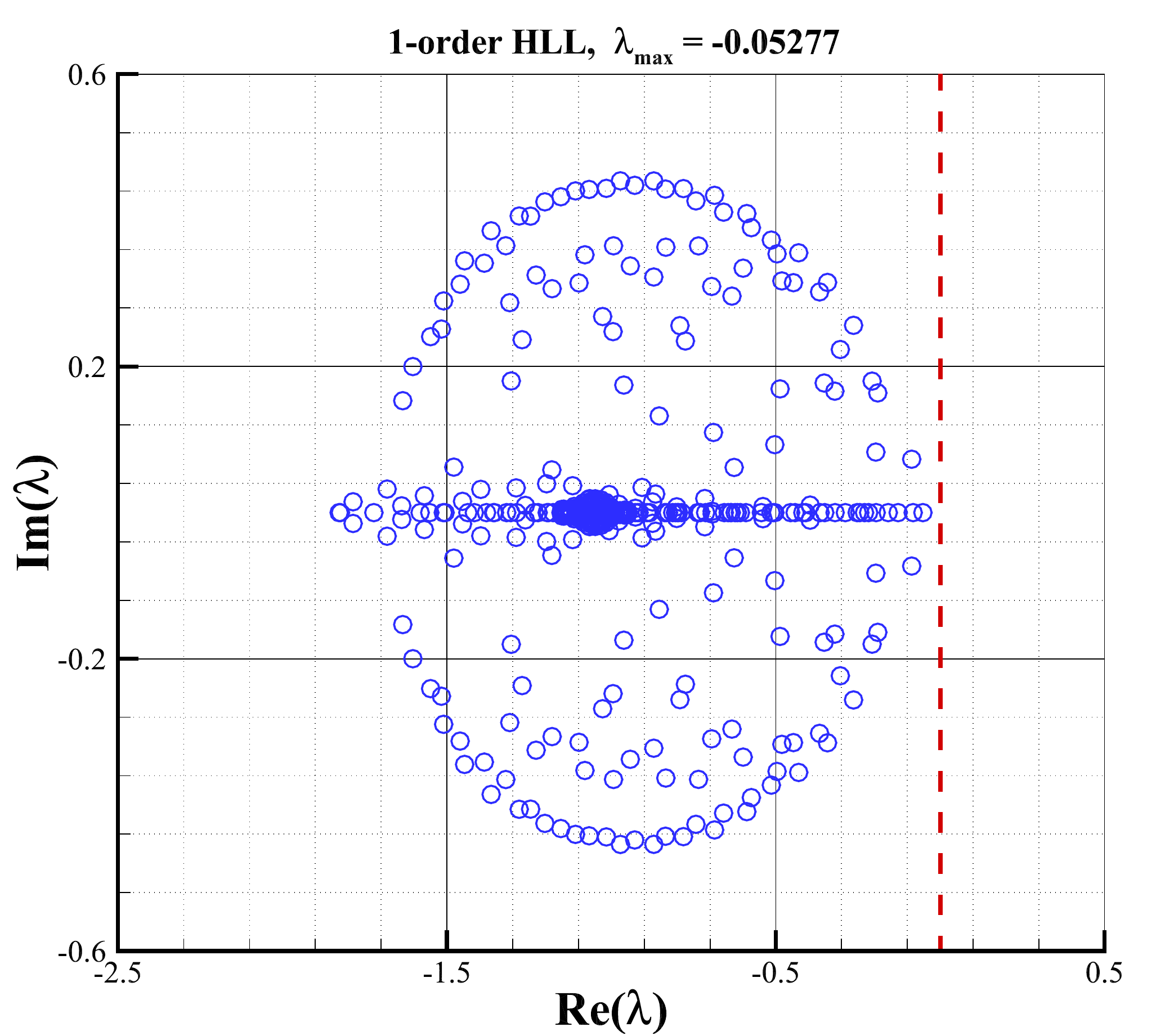}
	\end{minipage}
	}
	\subfigure[second-order HLL]{
	\begin{minipage}[t]{0.46\linewidth}
	\centering
	\includegraphics[width=0.9\textwidth]{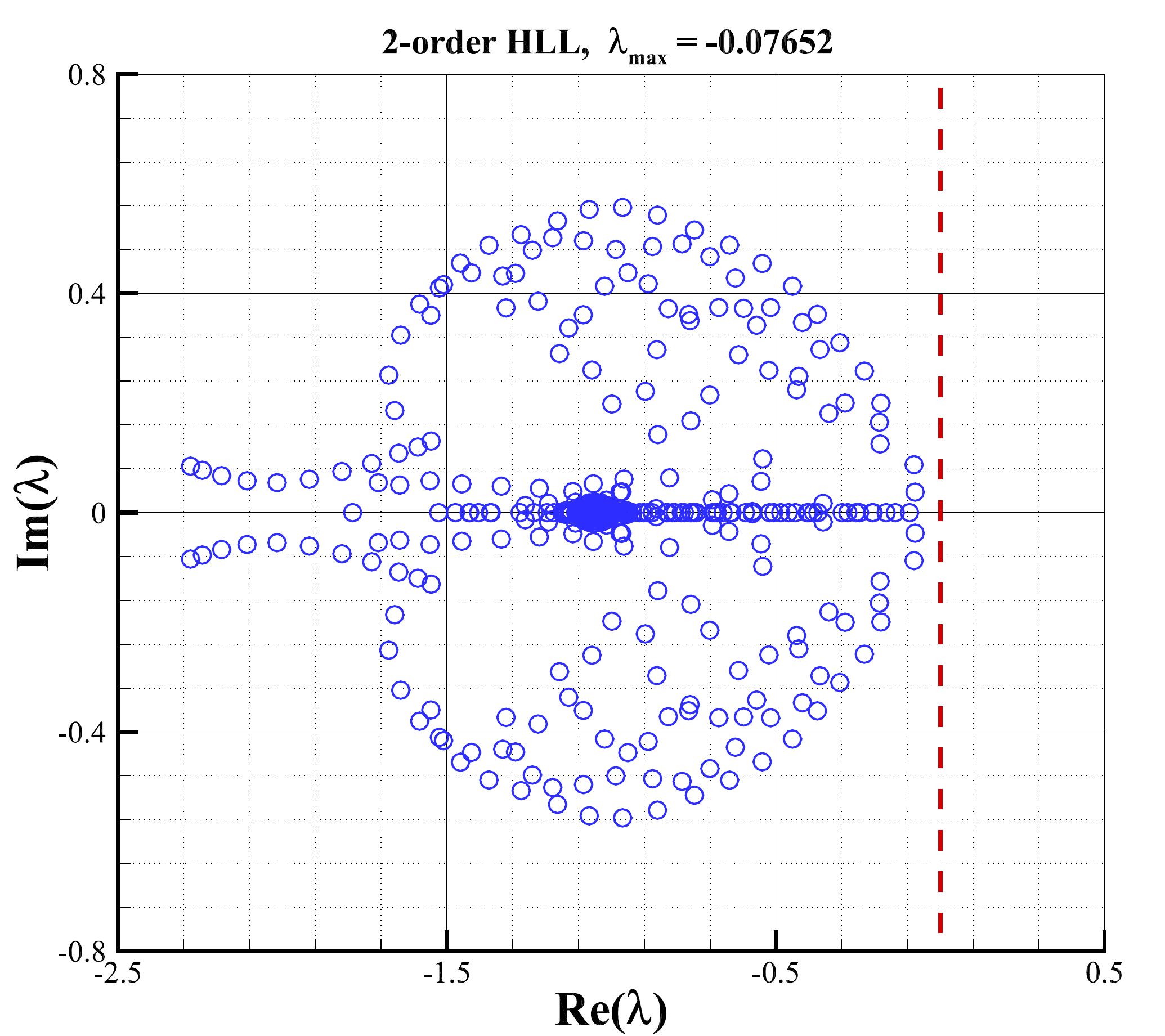}
	\end{minipage}
	}
	\centering
	\caption{Distribution of the eigenvalues of \textbf{S} in the complex plane for different Riemann solvers.(Grid with 11$ \times $11 cells, $ M_0=20 $ and $ \varepsilon =0.1 $. Left column: first-order schemes; right column: second-order schemes with van Albada limiter.)}\label{fig scatters a}
\end{figure}
\addtocounter{figure}{-1}
\begin{figure}[htbp]
	\centering
	\addtocounter{subfigure}{6}

	\subfigure[first-order van Leer]{
	\begin{minipage}[t]{0.46\linewidth}
	\centering
	\includegraphics[width=0.9\textwidth]{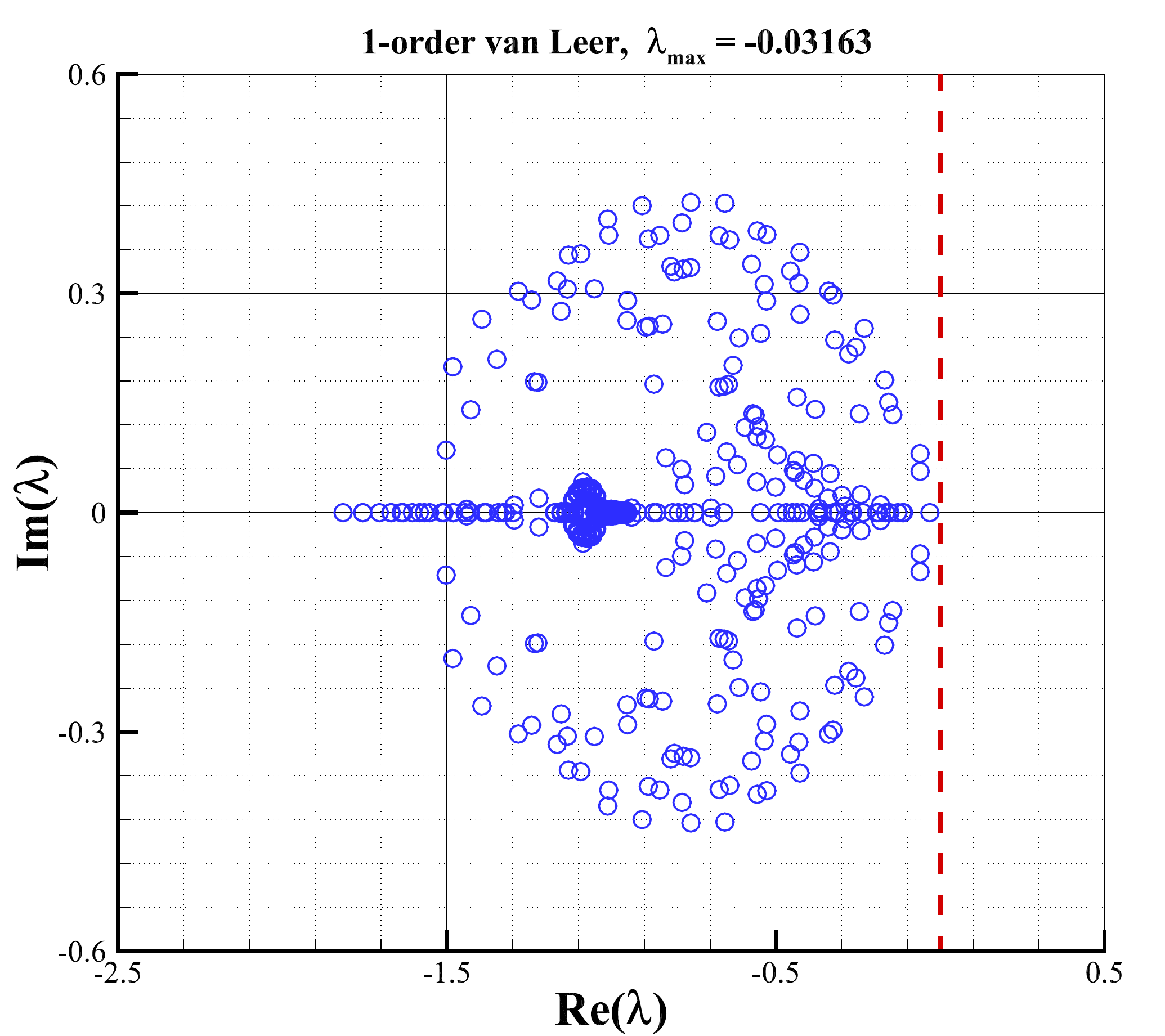}
	\end{minipage}
	}
	\subfigure[second-order van Leer]{
	\begin{minipage}[t]{0.46\linewidth}
	\centering
	\includegraphics[width=0.9\textwidth]{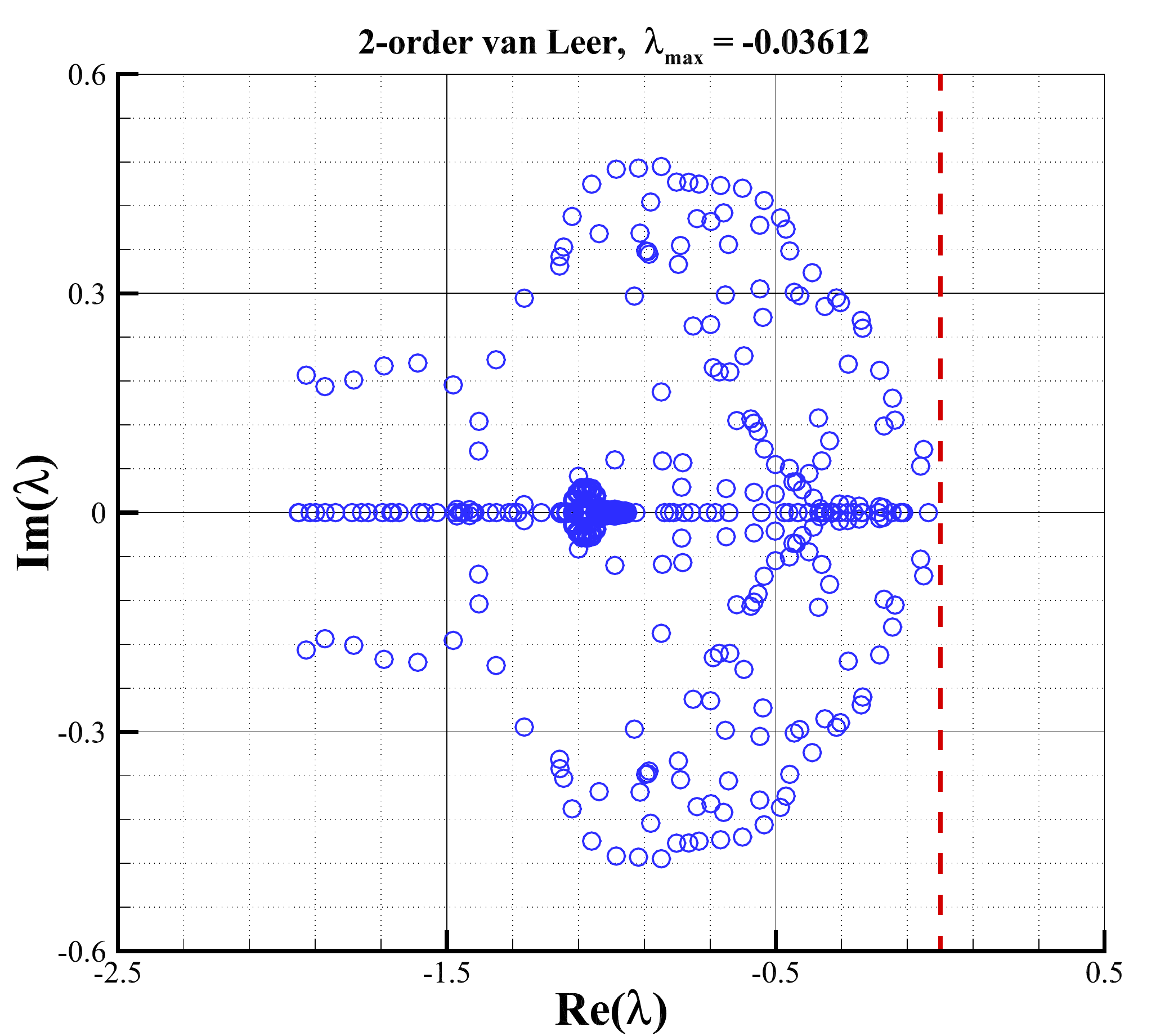}
	\end{minipage}
	}

	\subfigure[first-order AUSM$ ^+ $]{
	\begin{minipage}[t]{0.46\linewidth}
	\centering
	\includegraphics[width=0.9\textwidth]{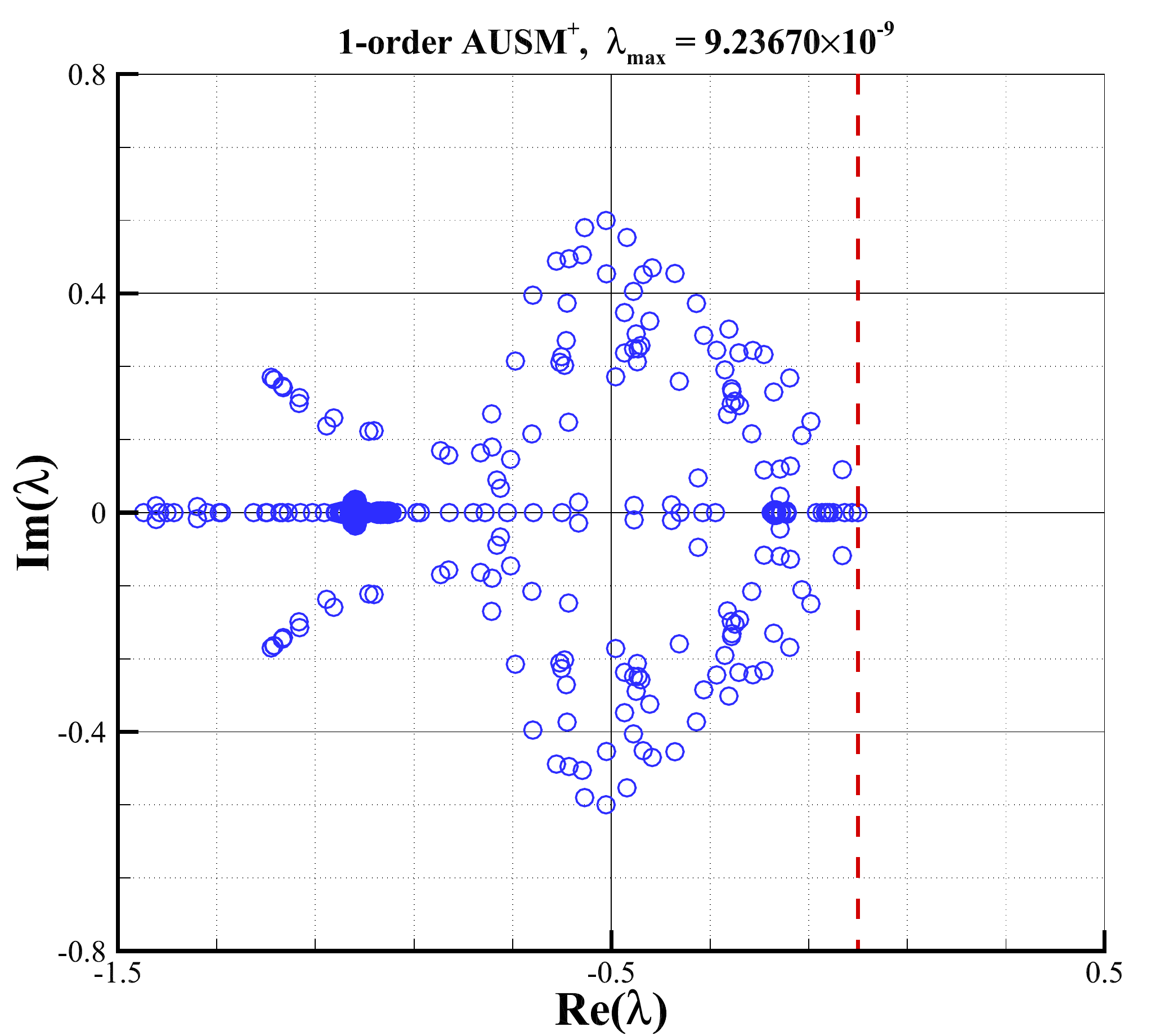}
	\end{minipage}
	}
	\subfigure[second-order AUSM$ ^+ $]{
	\begin{minipage}[t]{0.46\linewidth}
	\centering
	\includegraphics[width=0.9\textwidth]{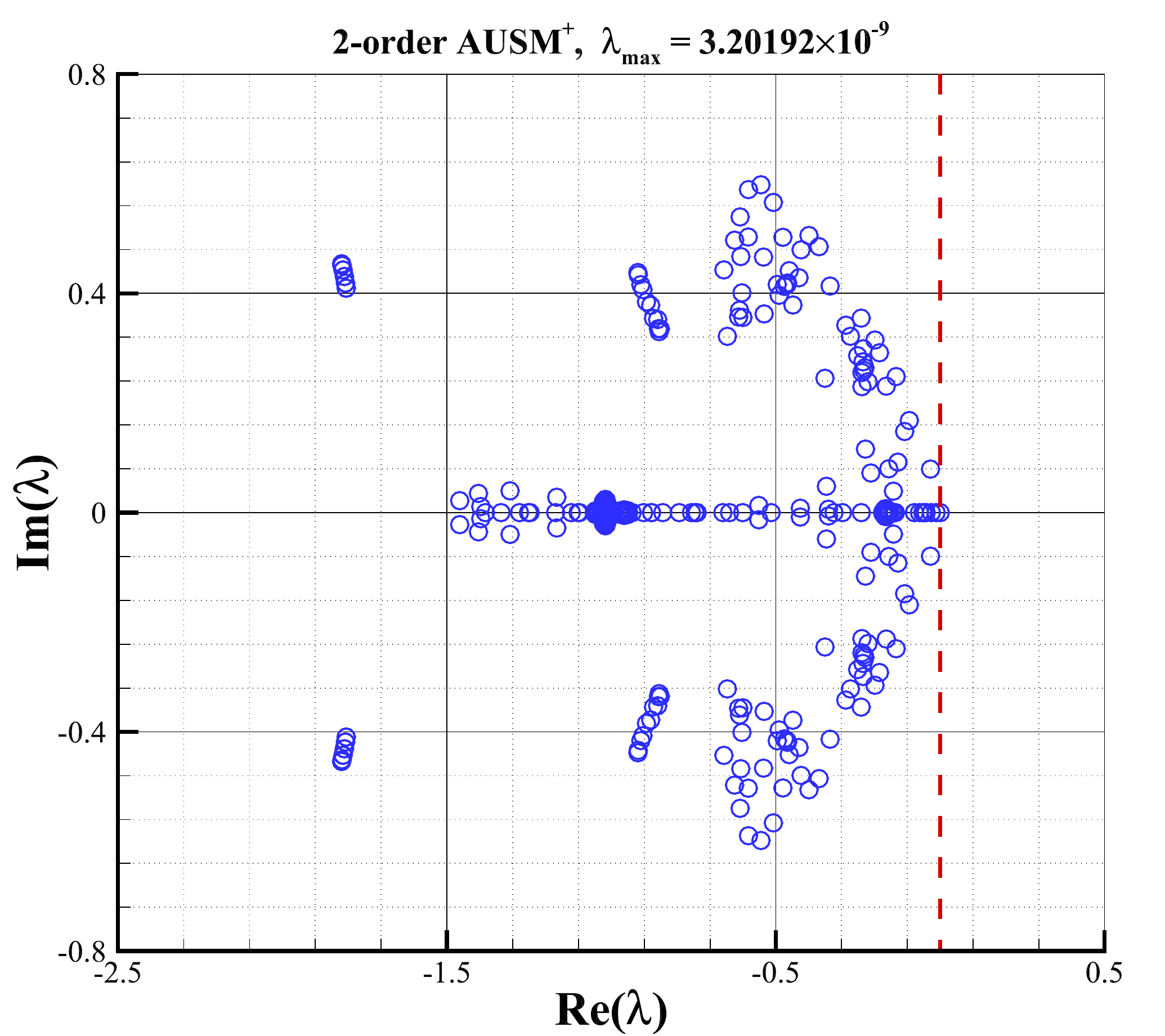}
	\end{minipage}
	}

	\centering
	\caption{Distribution of the eigenvalues of \textbf{S} in the complex plane for different Riemann solvers.(Grid with 11$ \times $11 cells, $ M_0=20 $ and $ \varepsilon =0.1 $. Left column: first-order schemes; right column: second-order schemes with van Albada limiter.)}\label{fig scatters b}
\end{figure}

It is well demonstrated that the shock instability is influenced by a variety of factors, such as Riemann solvers, shock intensity, grid distribution, grid aspect ratio, numerical shock structure, and so on \cite{Kitamura2009,Tu2014,Xie2019b,Henderson2007}. Numerous numerical experiments are used to investigate how these factors relate to shock instability \cite{Kitamura2009,Xie2017,Henderson2007,Tu2014,Kitamura2019,Ohwada2013}. The matrix stability analysis method connects the shock stability with the eigenvalues of the stability matrix, making it a quantitative tool for investigating shock stability. With this method, Dumbser et al.\cite{Dumbser2004} and Chen et al.\cite{Chen2018,Chen2018a} investigate the stability characteristics associated with these factors for the first-order scheme. How will these factors affect the stability of the second-order scheme? In this section, we will answer this question using the matrix stability analysis method established in section \ref{section 4}.

\subsection{Stability analysis for different Riemann solvers}\label{subsection 5.1}
It has been demonstrated that the Riemann solver plays an important role in triggering the shock instability \cite{Dumbser2004,Kitamura2009,Xie2017,Kitamura2019,Quirk1994}. In this section, we further investigate the influence of the Riemann solver on shock instability for the second-order scheme. Here, several approximate Riemann solvers with different dissipative characteristics are considered, including Roe \cite{Roe1981,Roe1986}, HLL \cite{Harten1983}, HLLC \cite{Toro1994}, AUSM$ ^+ $ \cite{Liou1996}, and van Leer \cite{VanLeer1982}. Note that in this paper, the HLL and HLLC Riemann solvers employ the wave speeds estimate proposed by Davis \cite{Davis1988}.\par

The stability analysis is carried out with the condition of $ M_0  = 20 $ and $ \varepsilon =0.1 $. All eigenvalues of the stability matrix for the first and second-order schemes with these Riemann solvers are shown in Fig.\ref{fig scatters a}. The right column shows the results of second-order schemes with the van Albada limiter, whereas the left column shows the results of first-order ones. According to Fig.\ref{fig scatters a}, it can be found that the sign of the maximal real parts of eigenvalues for second-order schemes is closely related to the Riemann solver and is consistent with that of the first-order schemes. It indicates that the Riemann solver still plays a key role in determining the shock stability of second-order schemes as it does for first-order cases. Moreover, As shown, eigenvalues with positive real parts are presented for Roe and HLLC Riemann solvers, which obviously, result in unstable behaviors. It is well-known that both solvers own minimal dissipation on contact and shear waves. On the contrary, the HLL and van Leer Riemann solvers always have eigenvalues with negative real parts, indicating their stability. And both solvers are known as dissipative schemes. What is more interesting is that the maximal real part of the AUSM$ ^+ $ Riemann solver is close to zero, indicating that AUSM$ ^+ $ Riemann solver is marginally stable.\par

\begin{figure}[htbp]
	\centering
	\subfigure[Roe Riemann solver]{
	\begin{minipage}[t]{0.8\linewidth}
	\centering
	\includegraphics[width=0.85\textwidth]{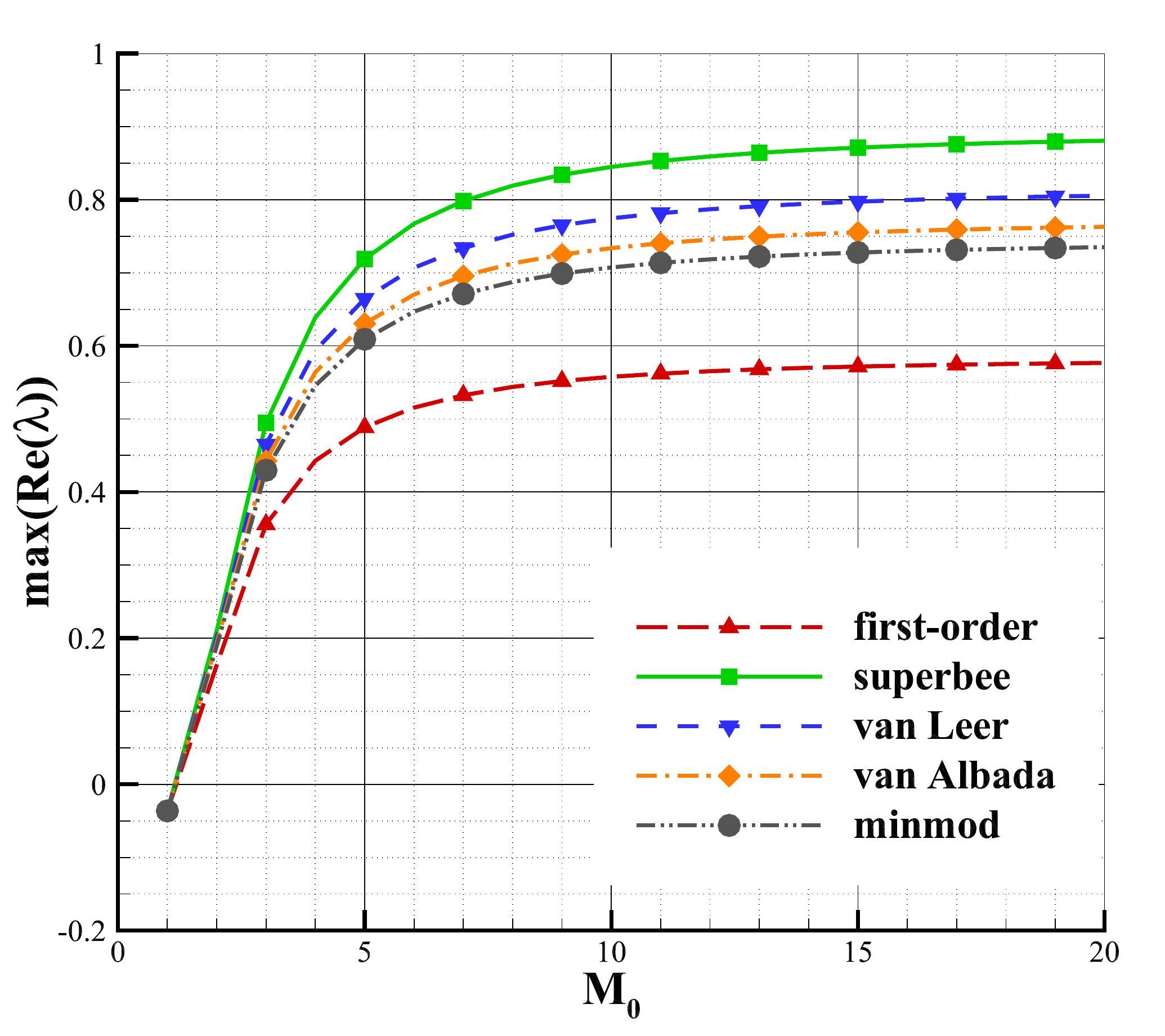}
	\end{minipage}
	}

	\subfigure[HLL Riemann solver]{
	\begin{minipage}[t]{0.8\linewidth}
	\centering
	\includegraphics[width=0.85\textwidth]{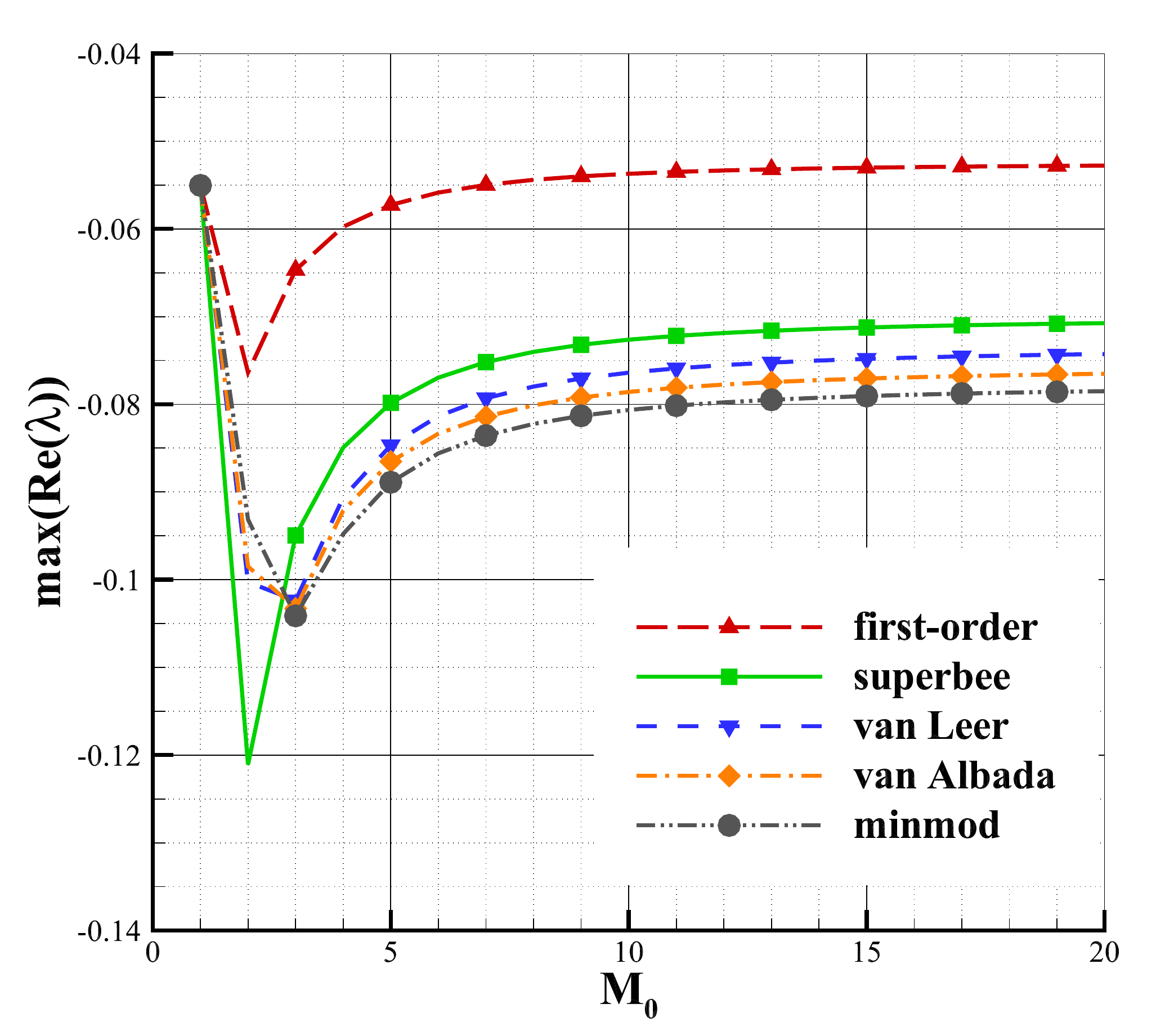}
	\end{minipage}
	}

	\centering
	\caption{Influence of special accuracy, limiter, and Mach number.(Grid with 11$ \times $11 cells, $ \varepsilon =0.1 $.)}\label{fig accuracy, limiters and Mach number}
\end{figure}

\subsection{Influence of spatial accuracy, limiter, and Mach number}\label{subsection 5.2}
To study the influence of spatial accuracy, limiter, and Mach number on shock instability, two typical Riemann solvers, Roe and HLL, are employed in this section. The same investigation can also be conducted on other solvers, such as HLLC, van Leer, and AUSM$ ^+ $, which is not shown for brevity. As seen in Fig.\ref{fig accuracy, limiters and Mach number}, the analysis is carried out for the second-order scheme, and in order to investigate the effect of spatial accuracy, results of the first-order scheme are also presented for comparison. Four limiters (superbee, van Leer, van Albada, and minmod) are employed when considering second-order accuracy. By analyzing the results in Fig.\ref{fig accuracy, limiters and Mach number}, the following conclusions can be obtained:\par

\begin{itemize}
	\item[1.]As shown in Fig.\ref{fig accuracy, limiters and Mach number}, we can find that the spatial accuracy will significantly affect the shock instability. The second-order scheme with Roe solver has a larger maximal real part of the eigenvalues than the first-order scheme, which suggests that the second-order scheme with Roe solver is more prone to be unstable than the first-order scheme. However, the maximal real part of the eigenvalues of the second-order scheme with HLL solver is less than that of the first-order case, indicating that the second-order scheme with the HLL solver is more stable than the first-order scheme. The same conclusion as HLL solver can also be obtained from AUSM$ ^+ $ solver, which is marginally stable as shown in Section \ref{subsection 5.1}. As a result, we can infer that the stability of the stable and marginally stable schemes will be better as the spatial accuracy is enhanced to second-order, while the stability of the unstable schemes will be poorer.
	\item[2.]For the second-order scheme, the limiter is a critical factor affecting the shock instability. The second-order scheme has the largest maximal real part of the eigenvalues when employing the superbee limiter with the lowest dissipation, while has the smallest maximal real part when using the minmod limiter with the highest dissipation. This indicates that the computation is more prone to be unstable when employing the limiter with lower dissipation. One exception is the HLL solver when $ M_0=2 $. As shown in Fig.\ref{fig accuracy, limiters and Mach number}(b), the second-order scheme with HLL solver is more unstable when using minmod limiter, but it is still stable. Anyway, if there is a need to alleviate the shock instability problem in supersonic and hypersonic flow simulations, the limiter with higher dissipation, such as minmod, is preferred.
	\item[3.]It can be seen from Fig.\ref{fig accuracy, limiters and Mach number} that the Mach number has a significant impact on the shock instability. When the Mach number is modest, the maximal real part of the eigenvalues increases quickly with the increasing Mach number, but when the Mach number is greater than a threshold, the growth will slow down and the maximal real part of the eigenvalues even remains unchanged. The threshold is determined by the Riemann solver, spatial accuracy, and limiter. It implies that when the Mach number is modest, its influence is evident, and the shock instability becomes more severe as the Mach number increases. However, the influence of the Mach number is limited. If the Mach number exceeds a threshold, it will do little to the shock stability. It should also be noted that the scheme with HLL solver is most stable when the Mach number is 2 or 3 (for the first-order scheme and second-order scheme with superbee limiter, the Mach number is 2, and for the others is 3.).
\end{itemize}

\begin{figure}[htbp]
	\centering
	\includegraphics[width=0.7\textwidth]{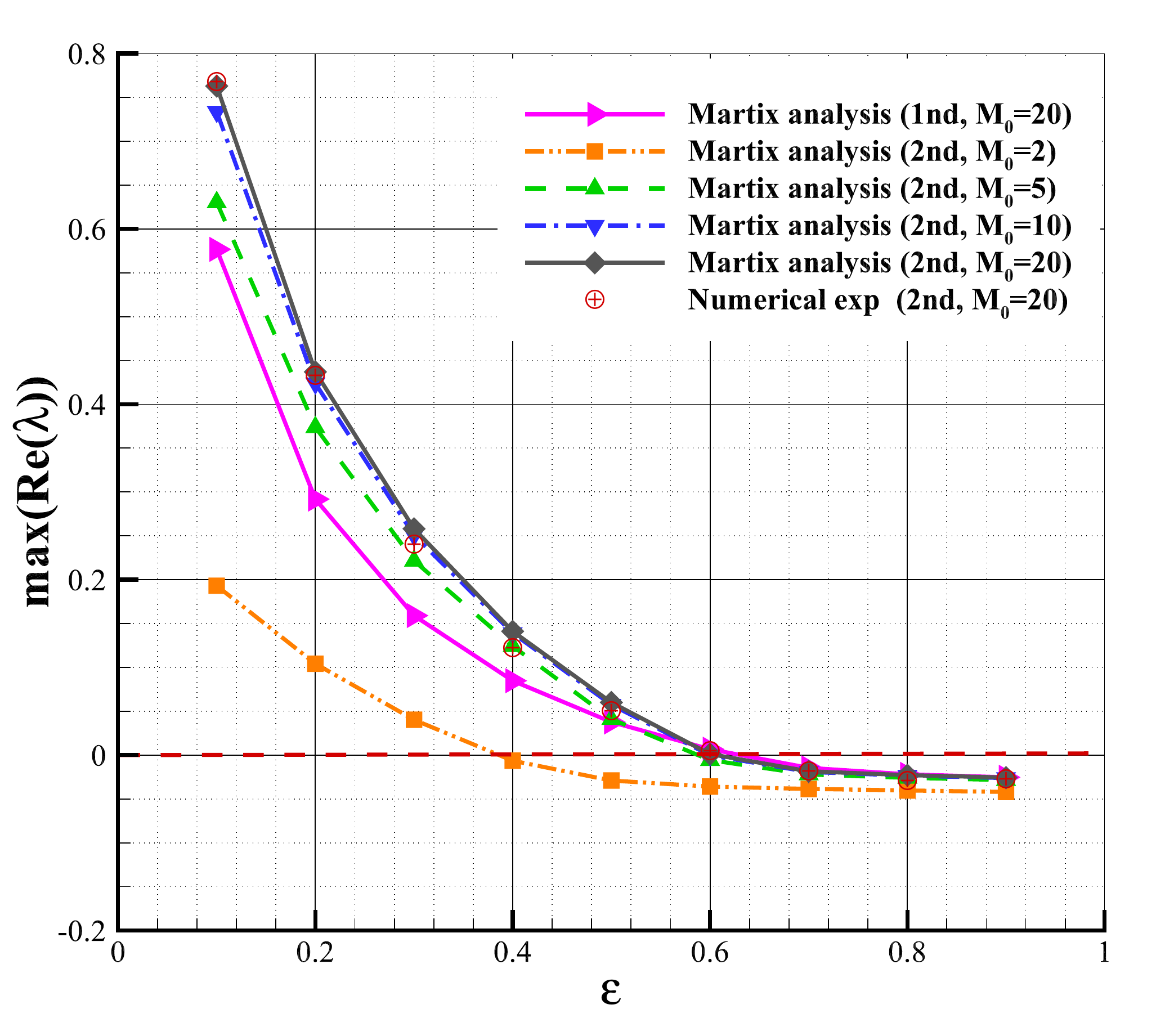}
	\caption{Influence of the numerical shock structure.(Grid with 11$ \times $11 cells, Roe solver, and van Albada limiter.)}
	\label{fig influence of shock structure}
\end{figure}

\begin{figure}[htbp]
	\centering
	\includegraphics[width=0.7\textwidth]{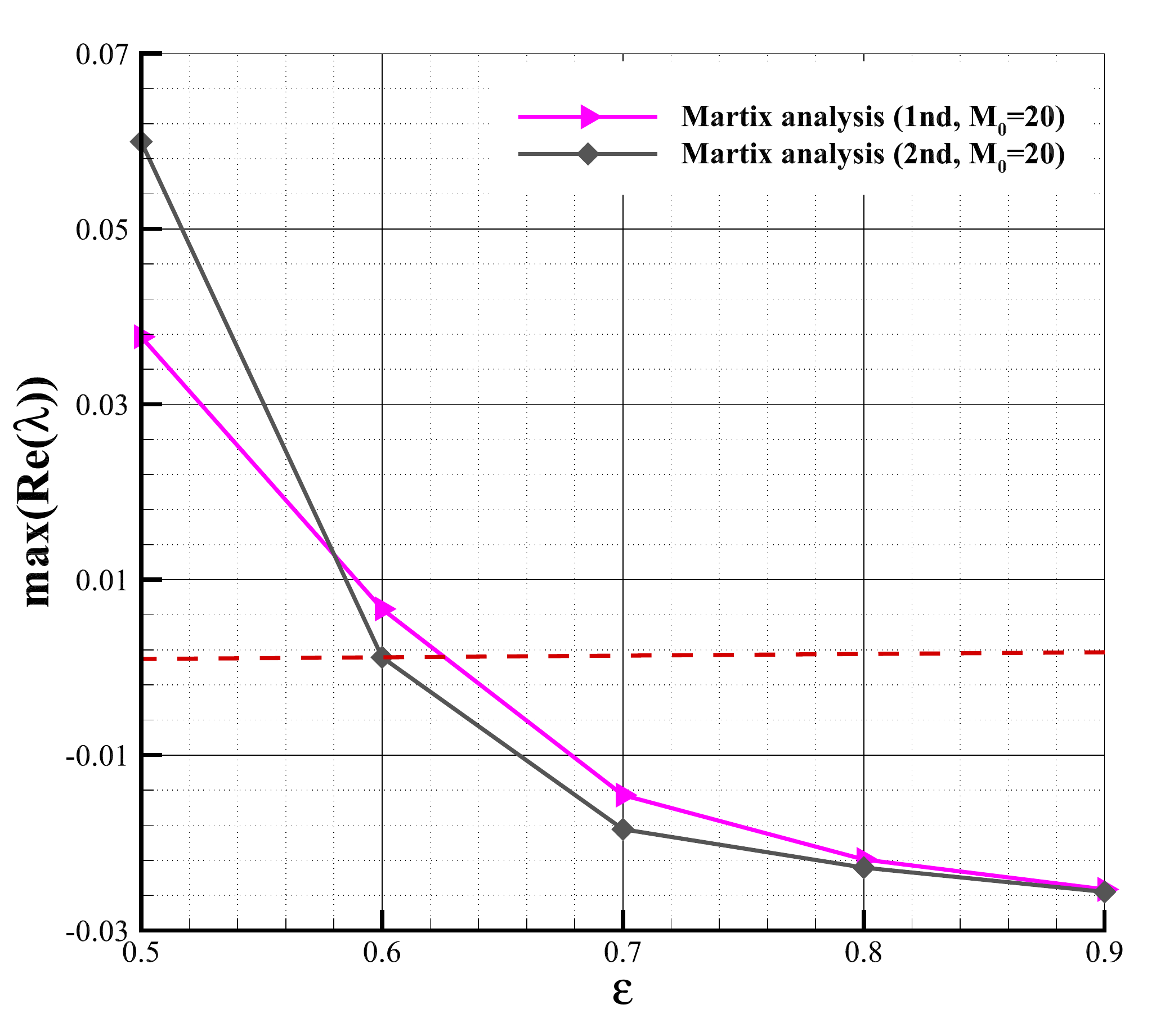}
	\caption{Influence of the numerical shock structure on the first and second-order schemes with $ 0.5\leq \varepsilon \leq0.9 $.(Grid with 11$ \times $11 cells, Roe solver, and van Albada limiter.)}
	\label{fig influence of shock structure(part)}
\end{figure}
\subsection{Analysis of the numerical shock structure}\label{subsection 5.3}

\begin{figure}[htbp]
	\centering
	\includegraphics[width=0.7\textwidth]{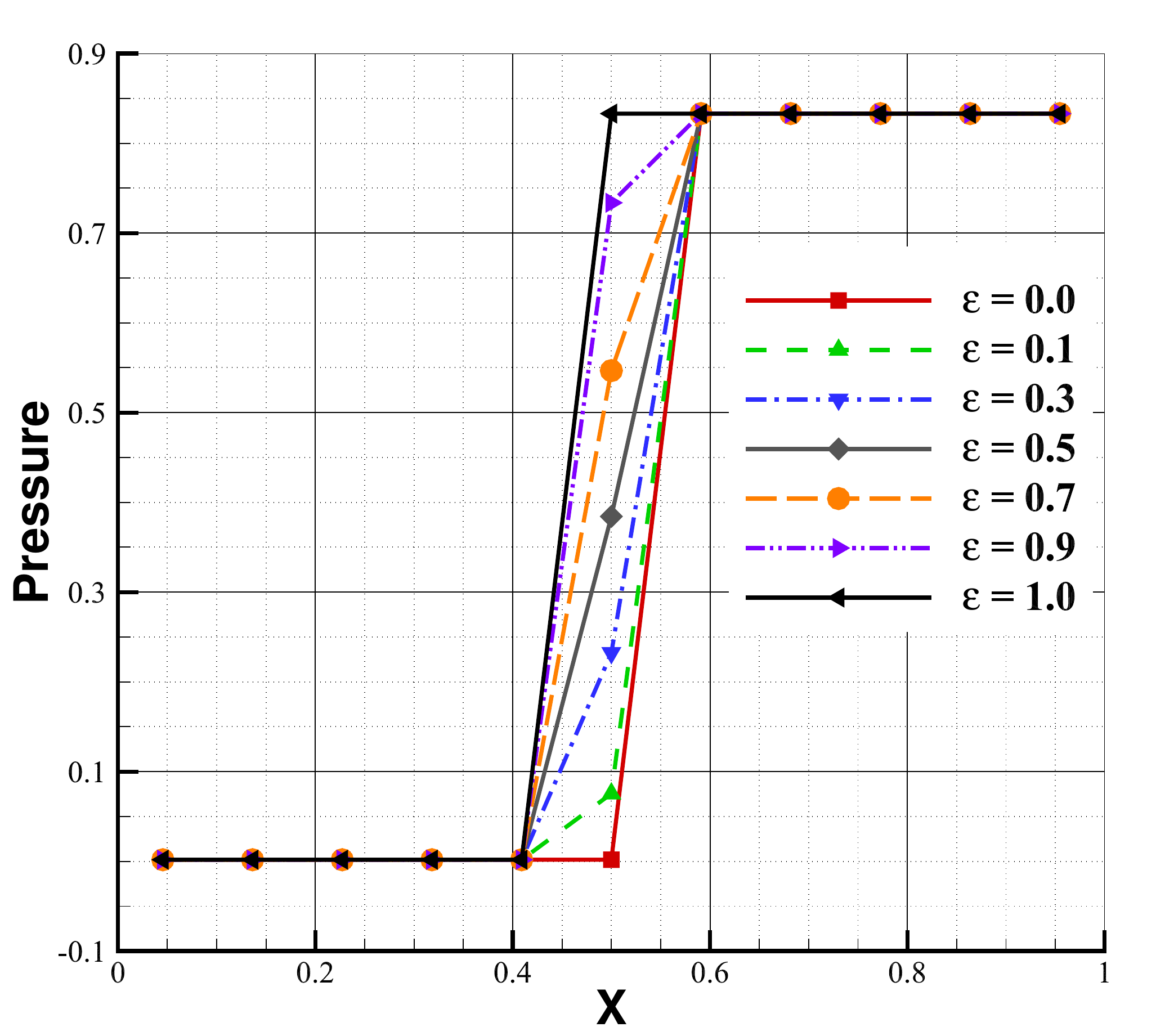}
	\caption{Convergent results of the 1D shock with different shock positions.(Grid with 1$ \times $11 cells, second-order scheme with Roe solver and van Albada limiter, $ M_0=20 $. )}
	\label{fig compare different shock positions}
\end{figure}

\begin{figure}[htbp]
	\centering
	\subfigure[$ \varepsilon =0.5 $]{
	\begin{minipage}[t]{0.46\linewidth}
	\centering
	\includegraphics[width=0.95\textwidth]{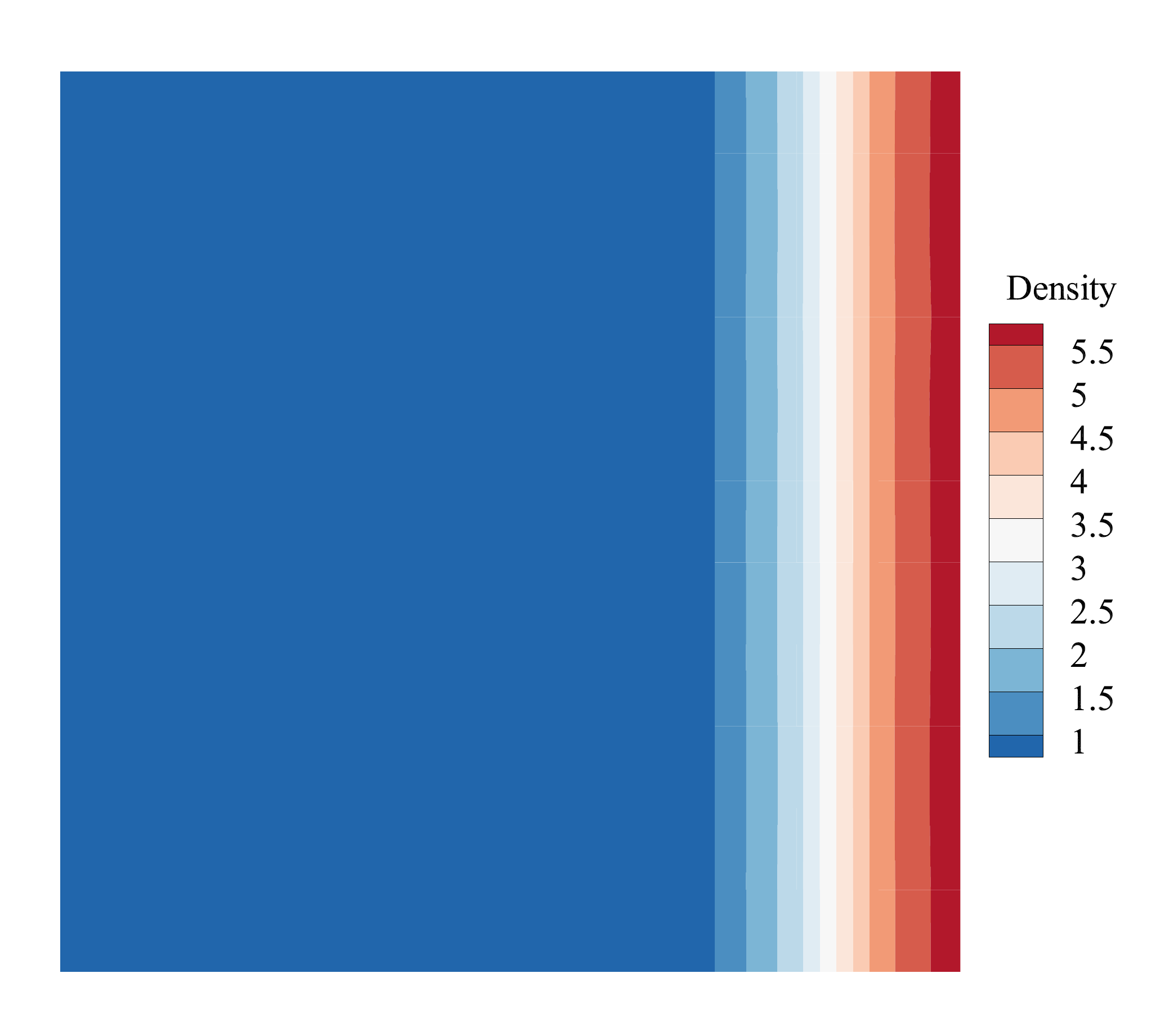}
	\end{minipage}
	}
	\subfigure[$ \varepsilon =0.7 $]{
	\begin{minipage}[t]{0.46\linewidth}
	\centering
	\includegraphics[width=0.95\textwidth]{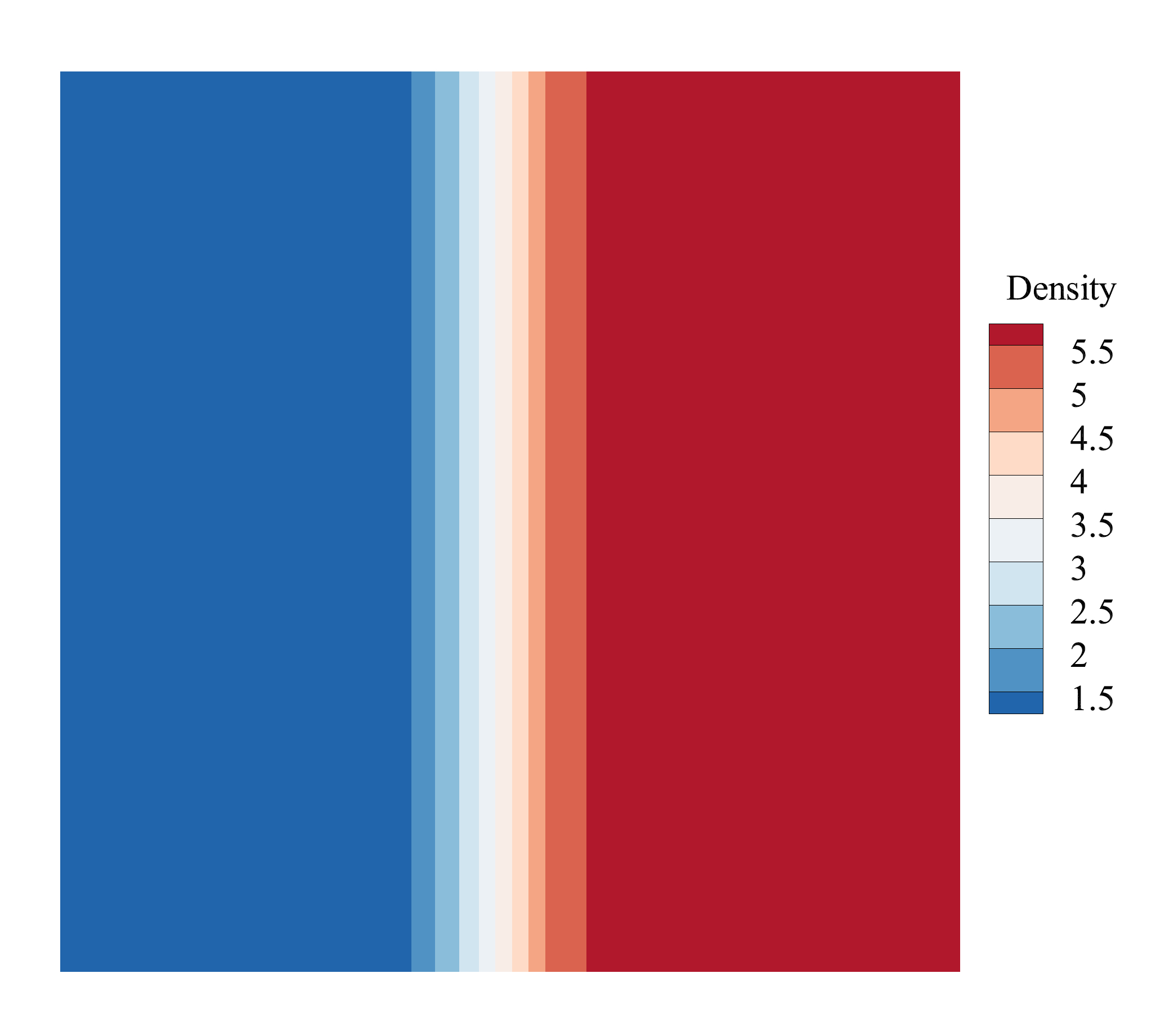}
	\end{minipage}
	}
	\centering
	\caption{Density contours at different numerical shock positions.(Grid with 11$ \times $11 cells, second-order schemes with Roe solver and van Albada limiter, t=1000.)}\label{fig contours at different shock position}
\end{figure}

It has been found that the numerical shock structure plays a vital role in shock instability. In this section, we investigate the impact of the numerical shock structure for the second-order scheme using the matrix stability analysis method. Since we are more concerned about the unstable behavior, Roe solver, which is proved to be more prone to be unstable, is employed to perform the remaining analysis. Fig.\ref{fig influence of shock structure} shows the effect of the numerical shock structure. As shown, with $ \varepsilon $ increasing, the maximal real part of the eigenvalues of the second-order scheme decreases, and when $ \varepsilon $ exceeds a threshold, the maximal real part of the eigenvalues begins to be less than 0. The threshold is related to the Mach number, around 0.4 when $ M_0=2 $ and around 0.6 when $ 5<M_0<20 $. Fig.\ref{fig compare different shock positions} is the convergent results of the 1D shock with different $ \varepsilon $. According to Fig.\ref{fig influence of shock structure} and Fig.\ref{fig compare different shock positions}, we can conclude that the second-order scheme will be more stable when the internal shock conditions are closer to the downstream states. The density contours with different numerical shock positions are shown in Fig.\ref{fig contours at different shock position}, demonstrating the validity of the matrix stability analysis.\par

Fig.\ref{fig influence of shock structure} also indicates that the numerical shock structure has the same impact on both first and second-order schemes. The unstable scheme can be stable by setting the internal shock conditions closer to the downstream states, whether for the first or second-order schemes. Moreover, When $ \varepsilon < 0.6 $, both the first and second-order schemes with Roe Riemann solver are unstable. And the maximal real part of the eigenvalues of the second-order scheme is greater than that of the first-order scheme, which suggests that the stability of the second-order scheme is worse than the first-order scheme. However, as shown in Fig.\ref{fig influence of shock structure} and Fig.\ref{fig influence of shock structure(part)}, when $ \varepsilon \geq  0.6 $ the maximal real parts of the eigenvalues for both first and second-order schemes are less than zero or slightly larger than zero, indicating that the computation is stable or marginally stable. Furthermore, in the region of $ \varepsilon \geq  0.6 $, the maximal real part of the eigenvalues of the first-order scheme is greater than that of the second-order scheme, suggesting that the stability of the second-order scheme is better than the first-order scheme, which is consistent with the conclusion in Section \ref{subsection 5.2}.\par

\begin{figure}[htbp]
	\centering
	\includegraphics[width=0.75\textwidth]{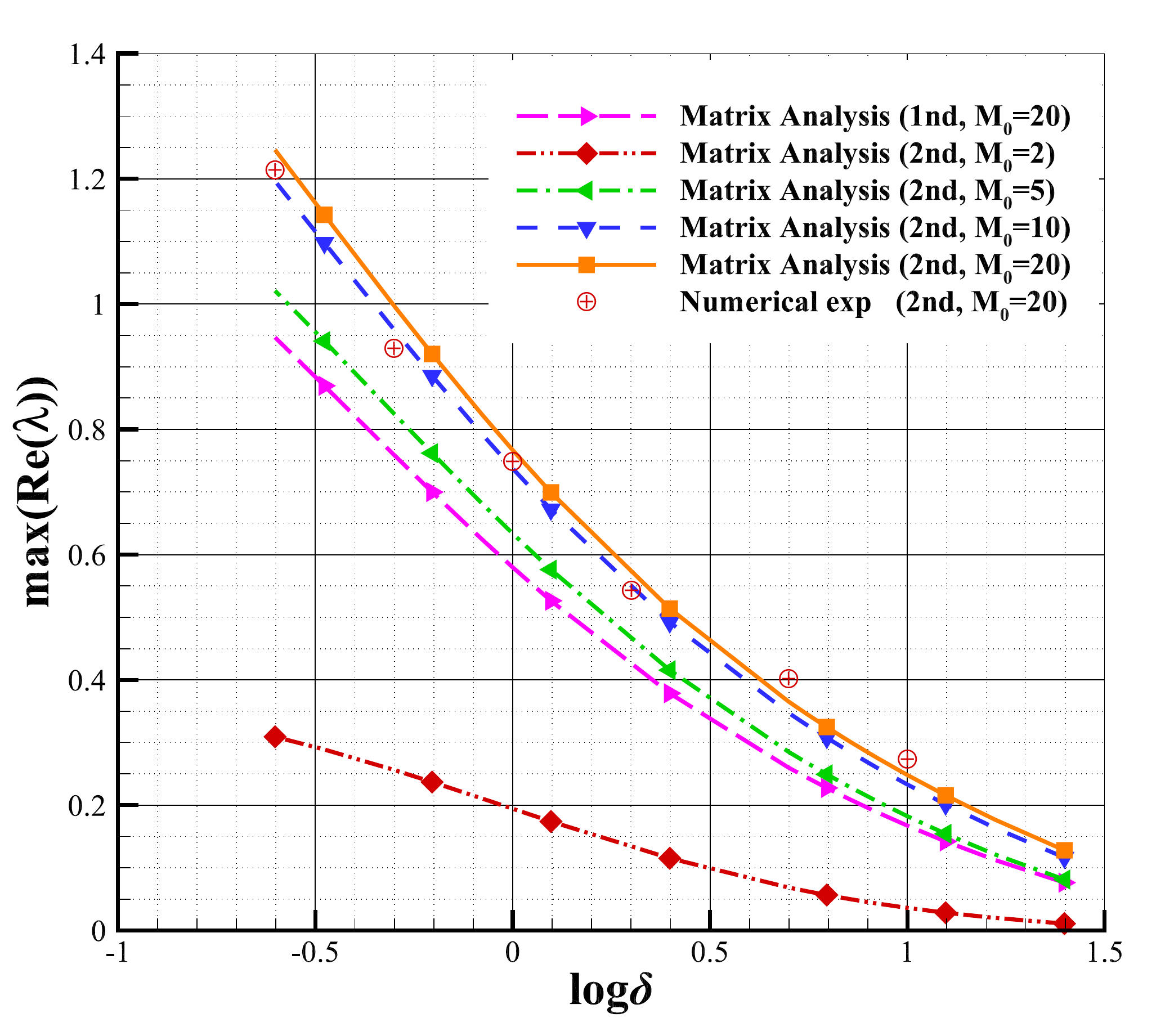}
	\caption{Influence of the aspect ratio.(Computational domain: $ 50\times50 $, 50 cells in x-direction and the cells in y-direction is changed, the initial shock is located in $ i=25 $, $ \varepsilon =0.1 $, Roe solver, and van Albada limiter.)}
	\label{fig influence of the aspect ratio}
\end{figure}

\subsection{Grid dependence of the shock instability}\label{subsection 5.4}
In addition to the influence of the Mach number, spatial accuracy, limiter function, Riemann solver, and the numerical shock structure, the computational grid also plays a significant role in shock instability \cite{Henderson2007,Ohwada2013,Pandolfi2001,Dumbser2004,Chen2018a}. In this section, we study further about how the grid affects the stability of the second-order scheme in capturing strong shocks from the perspectives of aspect ratio and distortion angle.\par

\subsubsection{The aspect ratio}
\begin{figure}[htbp]
	\centering
	\subfigure[$ \delta =0.5 $ (grid: $ 50\times100 $).]{
	\begin{minipage}[t]{0.46\linewidth}
	\centering
	\includegraphics[width=0.95\textwidth]{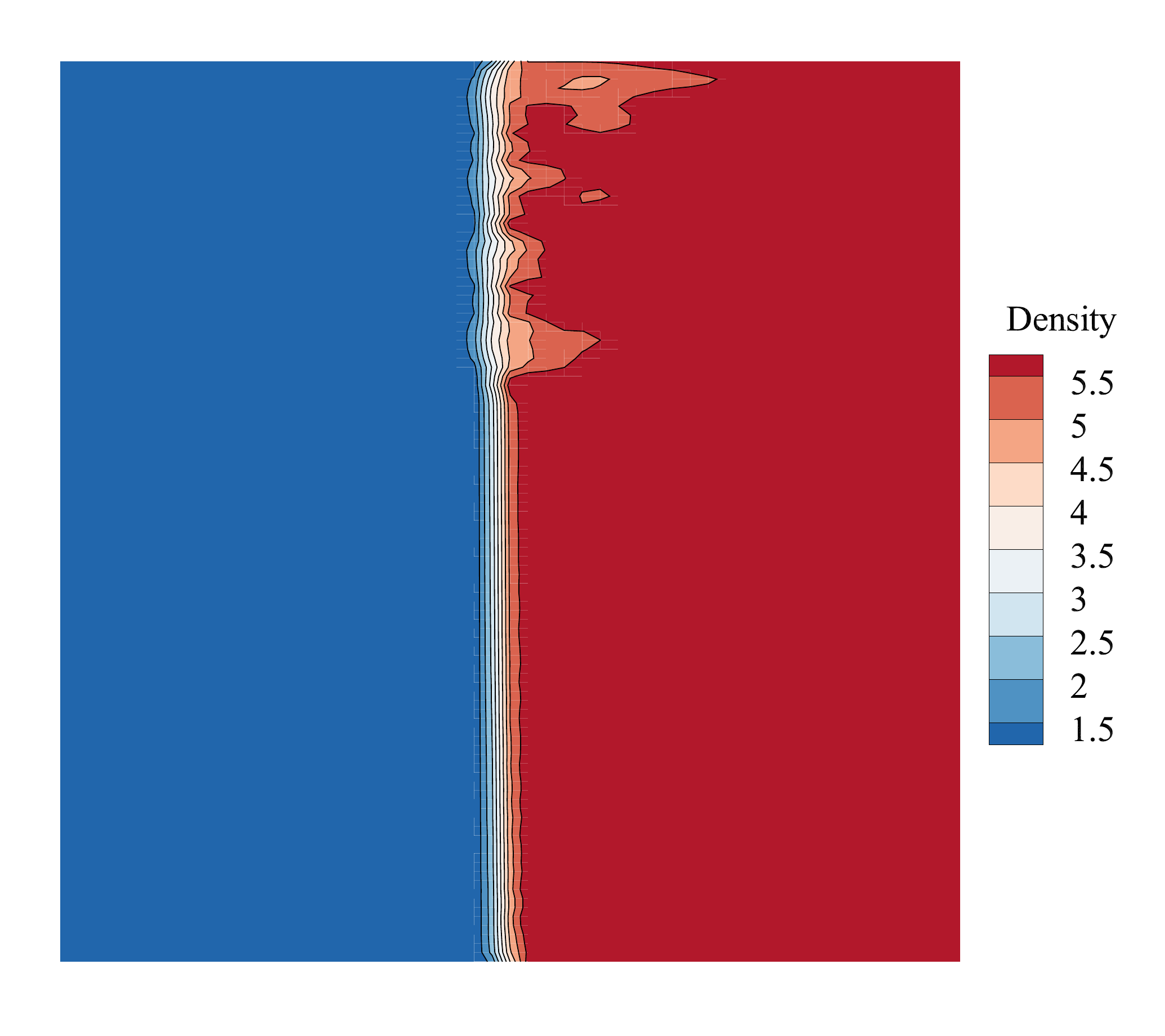}
	\end{minipage}
	}
	\subfigure[$ \delta =1 $ (grid: $ 50\times50 $).]{
	\begin{minipage}[t]{0.46\linewidth}
	\centering
	\includegraphics[width=0.95\textwidth]{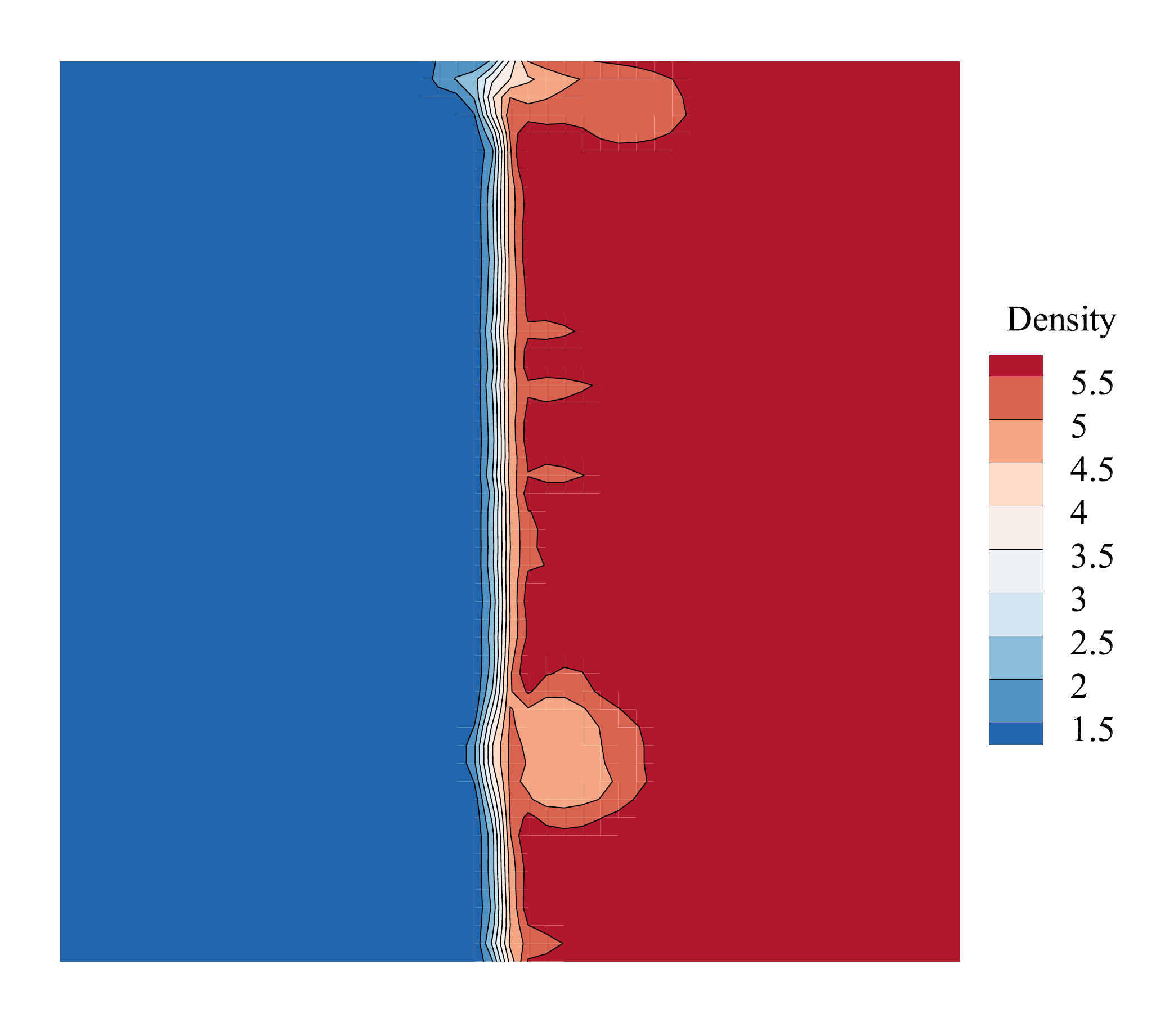}
	\end{minipage}
	}

	\subfigure[$ \delta =5 $ (grid: $ 50\times10 $).]{
	\begin{minipage}[t]{0.46\linewidth}
	\centering
	\includegraphics[width=0.95\textwidth]{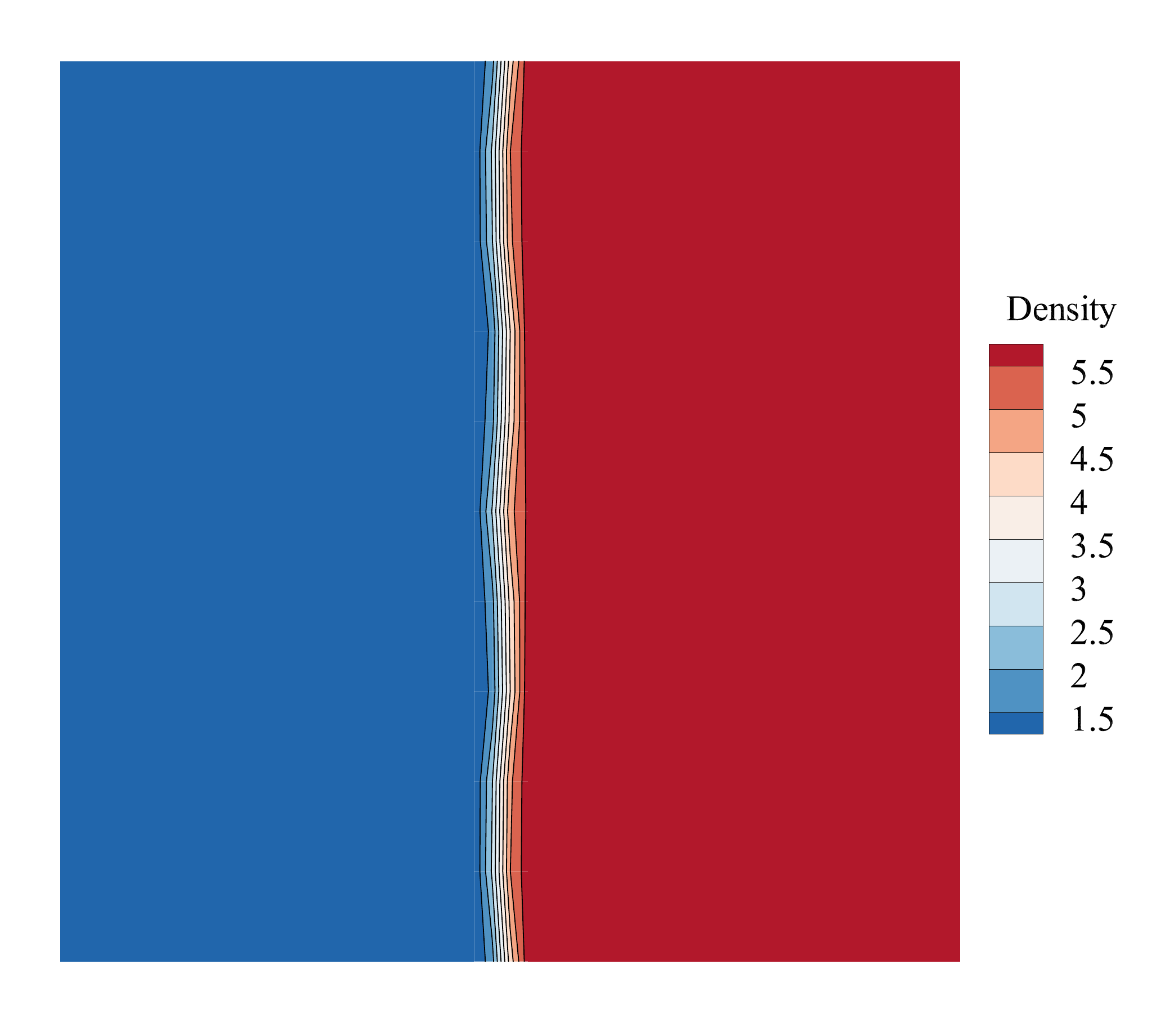}
	\end{minipage}
	}
	\subfigure[$ \delta =10 $ (grid: $ 50\times5 $).]{
	\begin{minipage}[t]{0.46\linewidth}
	\centering
	\includegraphics[width=0.95\textwidth]{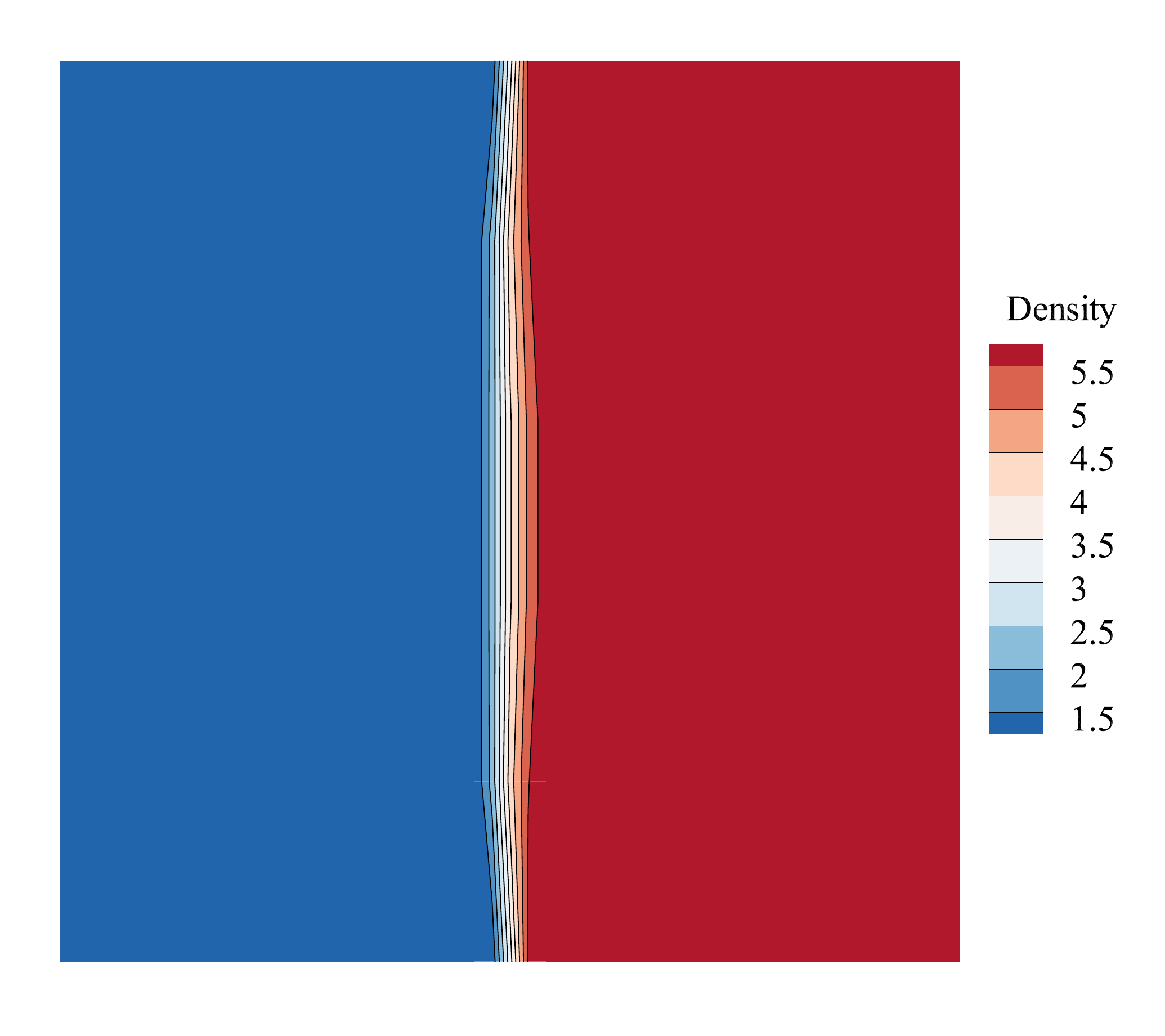}
	\end{minipage}
	}
	\centering
	\caption{The flow field with different aspect ratios.(Computational domain: $ 50\times50 $, 50 cells in x-direction and the cells in y-direction is changed, the initial shock is located in $ i=25 $, $ M_0=20 $ and $ \varepsilon =0.1 $, second-order scheme with Roe solver and van Albada limiter, t=500.)}\label{fig flow field with different aspect ratio}
\end{figure}

In this paper, the aspect ratio can be defined as
\begin{equation}
	\delta = \frac{\varDelta y}{\varDelta x},
\end{equation}
where $ \varDelta x $ represents the x-direction length of the cell, and $ \varDelta y $ is the y-direction length. In the current study, we analyze the case where the computational domain is $ 50 \times 50 $ and $ \varDelta x=1 $. The aspect ratio is then altered by changing the number of cells in the y-direction. Fig.\ref{fig influence of the aspect ratio} shows how the aspect ratio affects the maximal real part of the eigenvalues. According to Fig.\ref{fig influence of the aspect ratio}, it can be found that the maximal real part of eigenvalues decreases when $ \delta $ increases. Consequently, we may conclude that the shock instability will be alleviated when the grid becomes more elongated in the y-direction, while it will worsen when the grid is refined in the y-direction. The conclusion is validated in Fig.\ref{fig flow field with different aspect ratio} by numerical experiments, in which the flow fields become more stable as the aspect ratio increases. Furthermore, it can be inferred from Fig.\ref{fig influence of the aspect ratio} that the aspect ratio has the same effect on the first and second-order schemes: the shock instability can be alleviated by increasing the aspect ratio. And the second-order scheme is always more unstable than the first-order scheme with different aspect ratios under the same Mach number.\par

The analysis of the aspect ratio offers guidance for us during the grid generation. We can alleviate shock instability by increasing the aspect ratio near the shock, although the effect of increasing the aspect ratio is limited (As shown in Fig.\ref{fig influence of the aspect ratio} and \ref{fig flow field with different aspect ratio}, even if the aspect ratio is as large as 10, which may be too large for engineering application, the computation is still unstable.).\par

\subsubsection{The distortion angle}

\begin{figure}[htbp]
	\centering
	\includegraphics[width=0.4\textwidth]{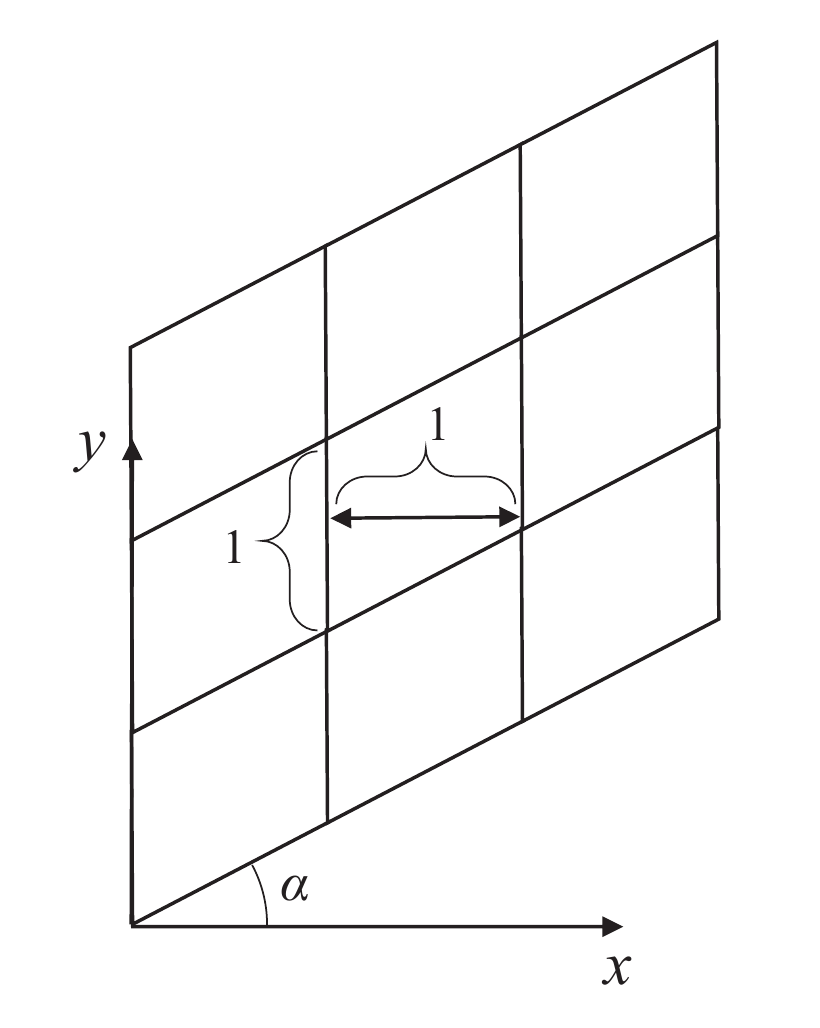}
	\caption{Definition of the distortion angle.}
	\label{fig define the angle}
\end{figure}

\begin{figure}[htbp]
	\centering
	\includegraphics[width=0.75\textwidth]{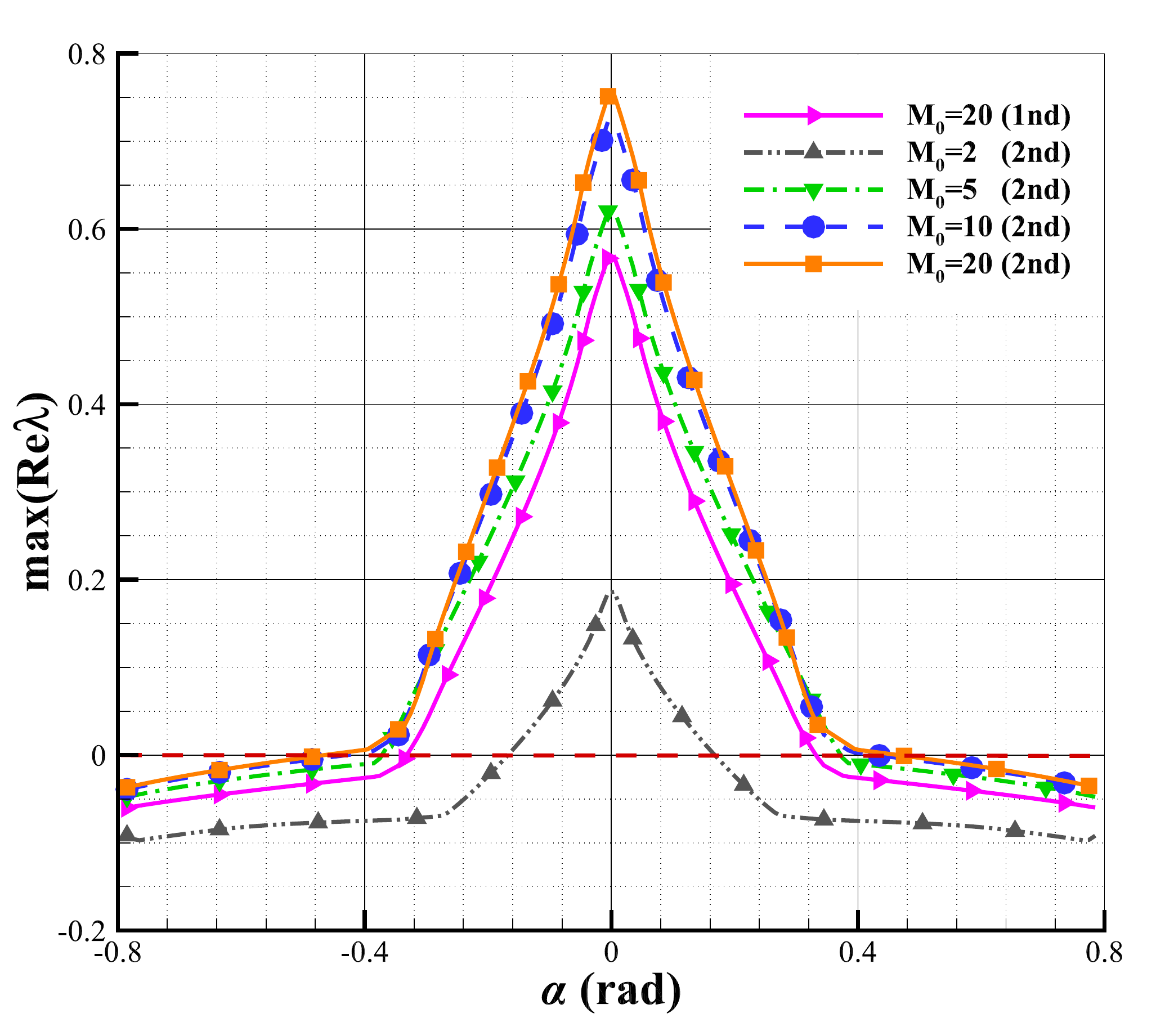}
	\caption{Influence of the distortion angle.($ \varepsilon =0.1 $, Roe solver, and van Albada limiter.)}
	\label{fig influence of angle}
\end{figure}

\begin{figure}[htbp]
	\centering
	\subfigure[$ \alpha  =0^{\circ} $]{
	\begin{minipage}[t]{0.46\linewidth}
	\centering
	\includegraphics[width=0.95\textwidth]{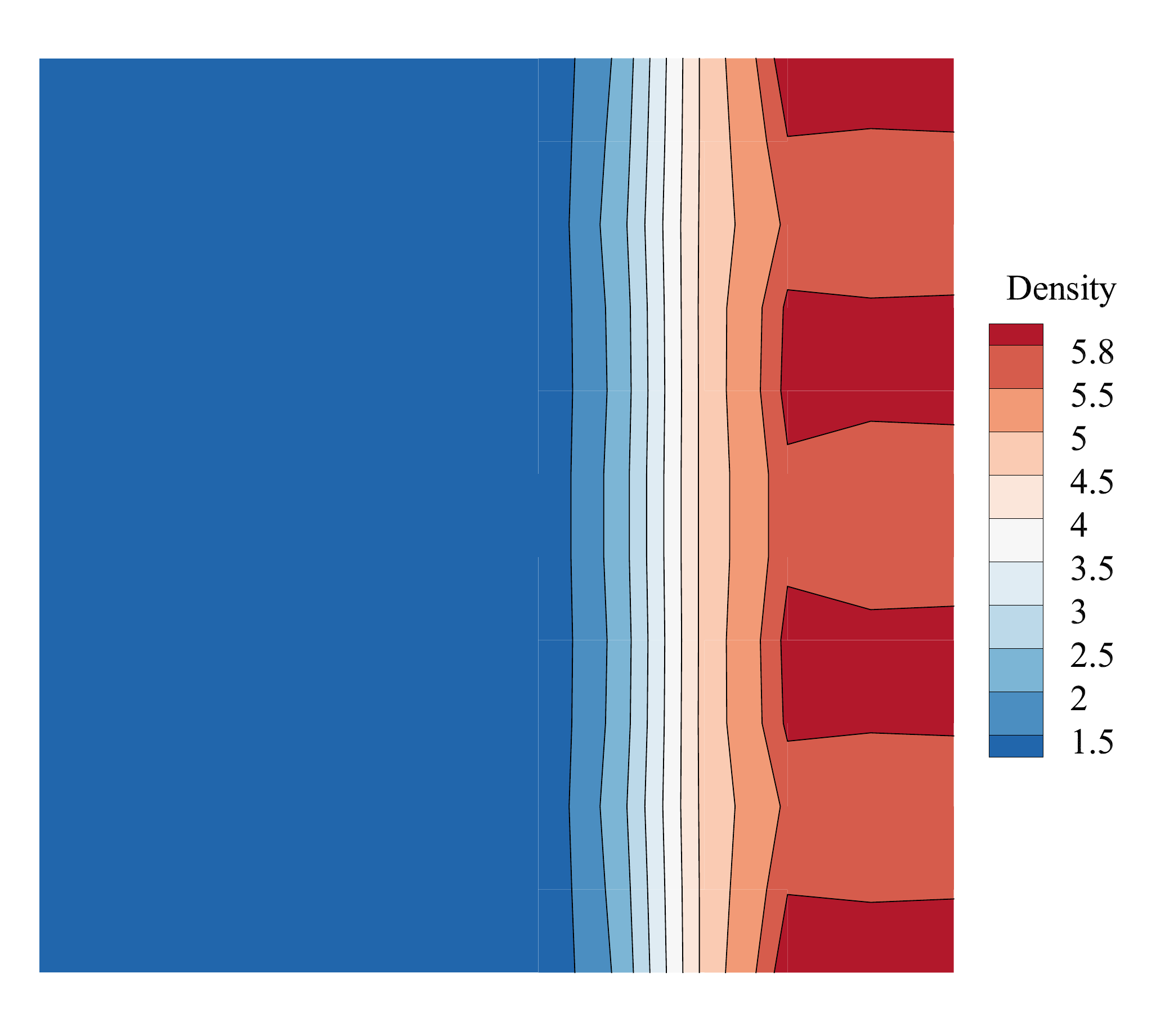}
	\end{minipage}
	}
	\subfigure[$ \alpha  =20^{\circ} $]{
	\begin{minipage}[t]{0.46\linewidth}
	\centering
	\includegraphics[height=1.0\textwidth]{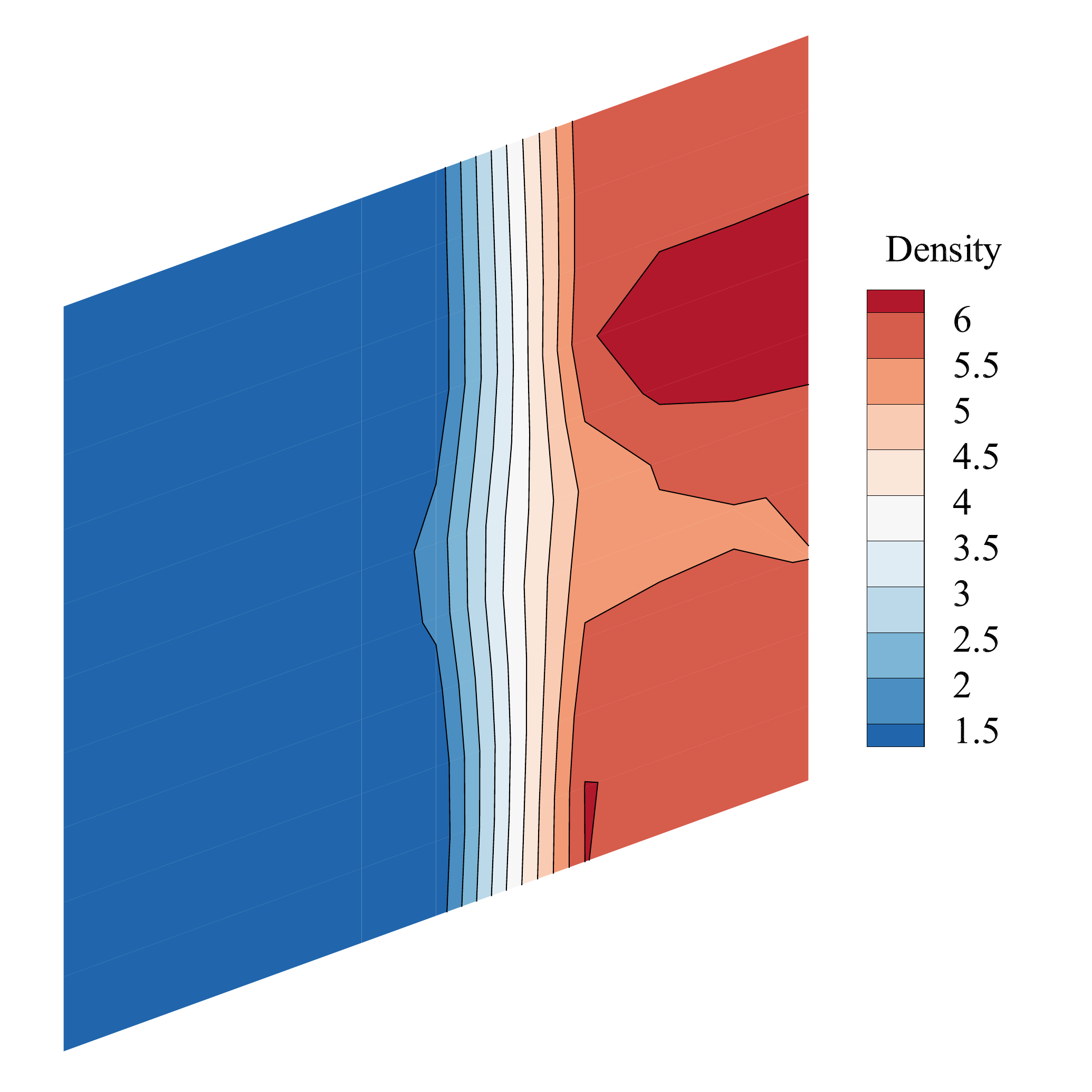}
	\end{minipage}
	}

	\subfigure[$ \alpha  =30^{\circ} $]{
	\begin{minipage}[t]{0.46\linewidth}
	\centering
	\includegraphics[height=1.0\textwidth]{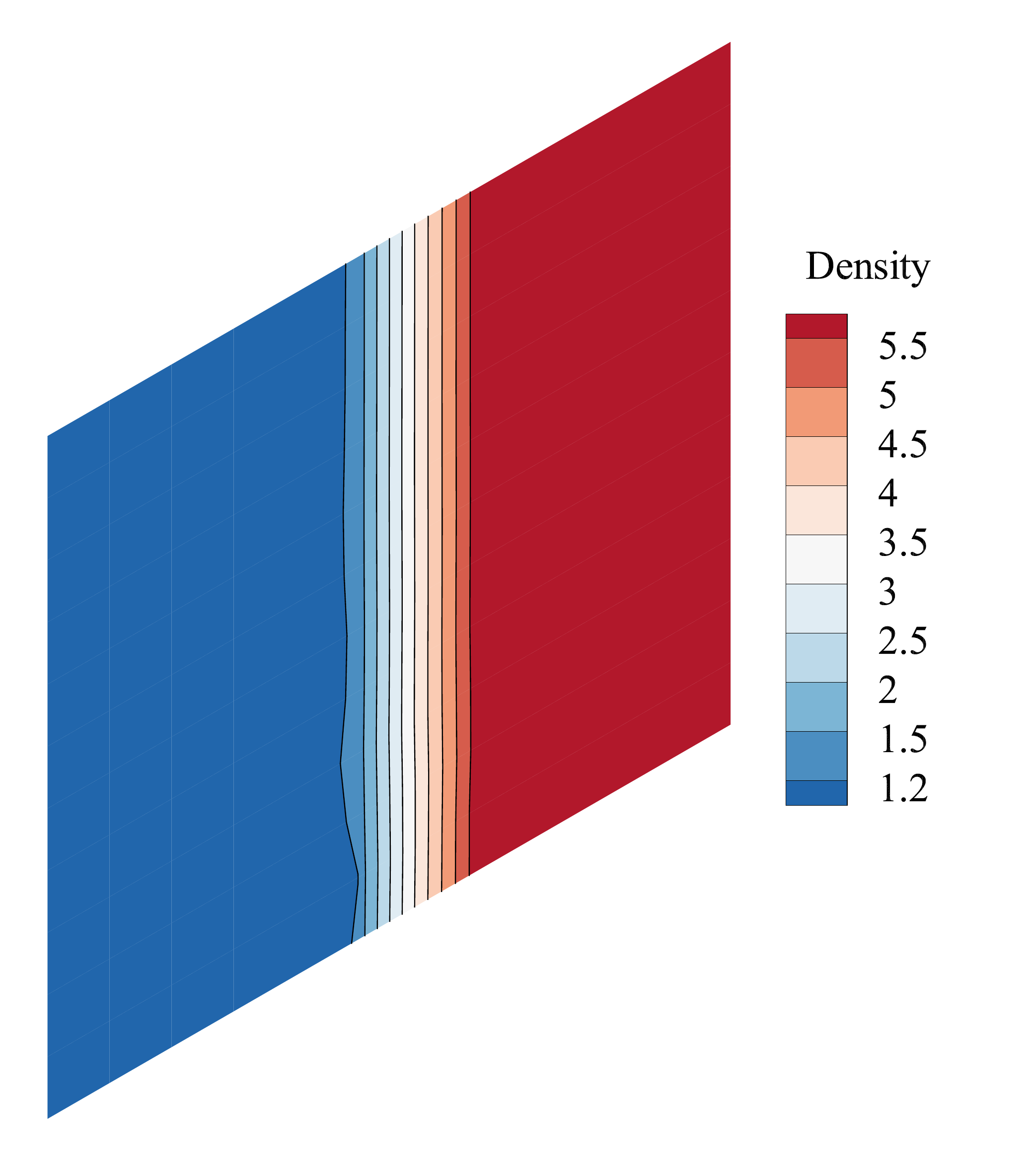}
	\end{minipage}
	}
	\subfigure[$ \alpha  =40^{\circ} $]{
	\begin{minipage}[t]{0.46\linewidth}
	\centering
	\includegraphics[height=1.0\textwidth]{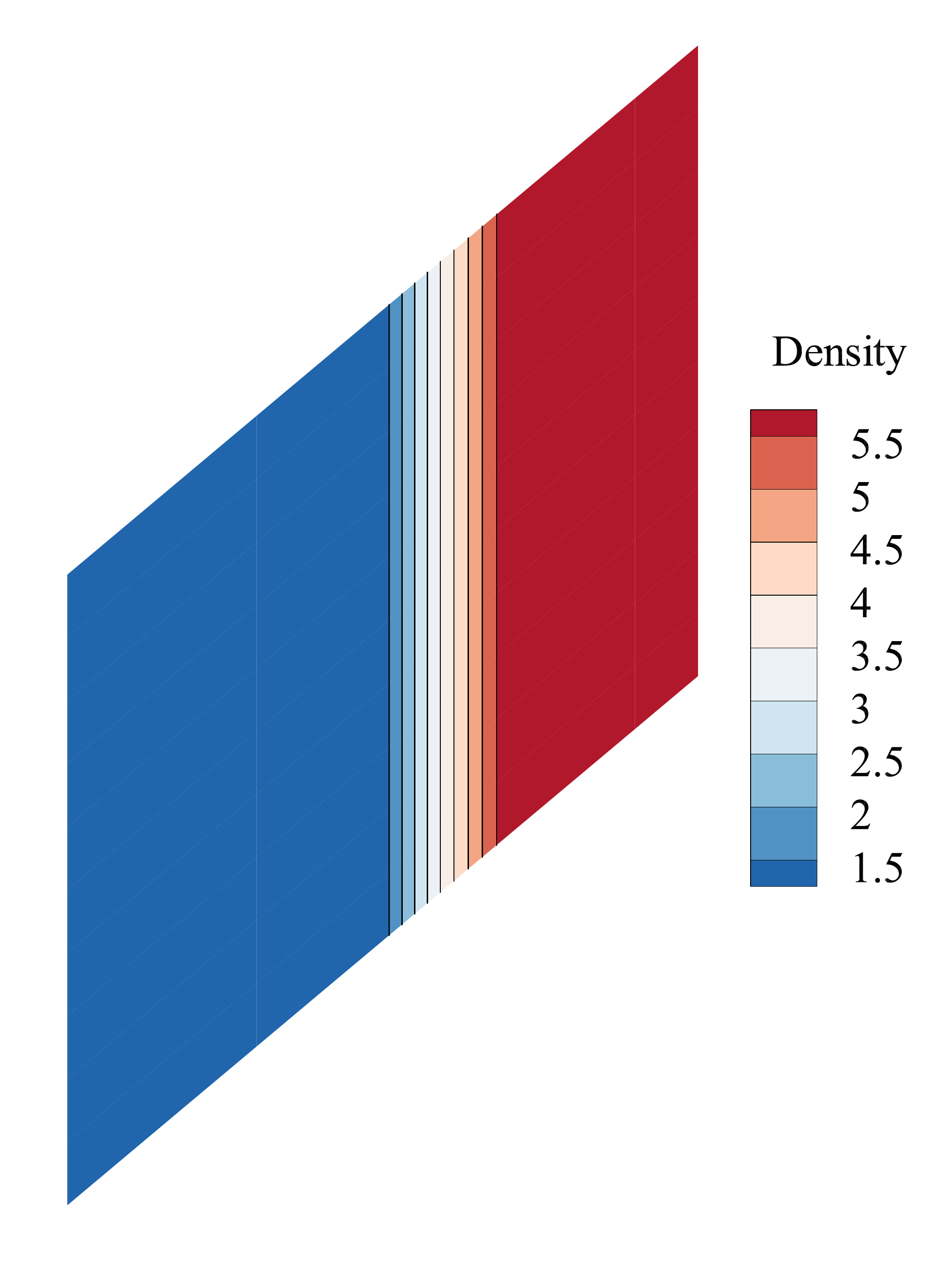}
	\end{minipage}
	}
	\centering
	\caption{The flow field with different distortion angles.(Grid with $ 11\times11 $ cells, $ M_0=20 $ and $ \varepsilon =0.1 $, second-order scheme with Roe solver and van Albada limiter, t=500.)}\label{fig flow field with different distortion angle}
\end{figure}

Apart from the aspect ratio, the distortion angle is another factor that could affect the stability of capturing strong shocks \cite{Chen2018a}. The present work defines the distortion angle as shown in Fig.\ref{fig define the angle}. Be aware that $ \alpha $ is positive when the grid deflects in the counterclockwise direction and negative in the other case. We take $ -45^{\circ}<\alpha<45^{\circ} $ in the current work. Fig.\ref{fig influence of angle} shows that the maximal real part of the eigenvalues of \textbf{S} is a function of $ \alpha $. As shown, for the second-order scheme with the Roe solver, the maximal real part of eigenvalues decreases as $ | \alpha | $ increases, and if $ | \alpha | $ becomes large enough, it may turn negative. So, the shock instability for the second-order scheme can be cured by increasing the distortion angle of the grid, and the threshold is related to the Mach number. When $ M_0=2 $, the threshold is around $ 10^{\circ} $, and it will be around $ 27^{\circ} $ if $ M_0=20 $. According to Fig.\ref{fig flow field with different distortion angle}, which depicts the flow field with different distortion angles, the shock becomes more stable as the distortion angle increases. It confirms the result of the matrix stability analysis. Note that Zhang et al.\cite{Zhang2017b} find that the carbuncle phenomenon will be eliminated when employing the triangular grid instead of the quadrilateral grid. This may be interpreted by the conclusion here. Moreover, If we compare the differences between the first and second-order schemes, we can find that the maximal real part of eigenvalues of the first-order scheme is always smaller than that of the second-order scheme as the distortion angle increases, whether it will be stable or unstable. As a result, it can be inferred that the stability of the second-order scheme is always worse than that of the first-order scheme as the distortion angle increases.\par

\subsection{Spatial localization of the source of instability}\label{subsection 5.5}

In the above sections, we devote our efforts to exploring the primary numerical characteristics of shock instability for the second-order finite-volume scheme. In this section, we will study further about an important problem about the mechanism underlying the shock instability: spatial localization of the source of instability. Dumbser et al.\cite{Dumbser2004} and Chen et al.\cite{Chen2018a} investigate this problem and find that the shock instability origins from the upstream of the shock when the shock is resolved exactly between two cells. However, in actual simulation, the captured shock is rarely exactly between two cells. The effect of numerical shock structure is studied by Xie et al.\cite{Xie2017}. Based on a series of well-designed numerical experiments, they find that the perturbations inside the shock structure are the main factors to trigger the instability. These studies concentrate on first-order schemes. In this section, we will investigate the spatial localization of the source of instability for the second-order scheme using the matrix stability analysis method.\par

\begin{figure}[htbp]
	\centering
	\subfigure[The upstream domain.]{
	\begin{minipage}[t]{0.31\linewidth}
	\centering
	\includegraphics[height=0.95\textwidth]{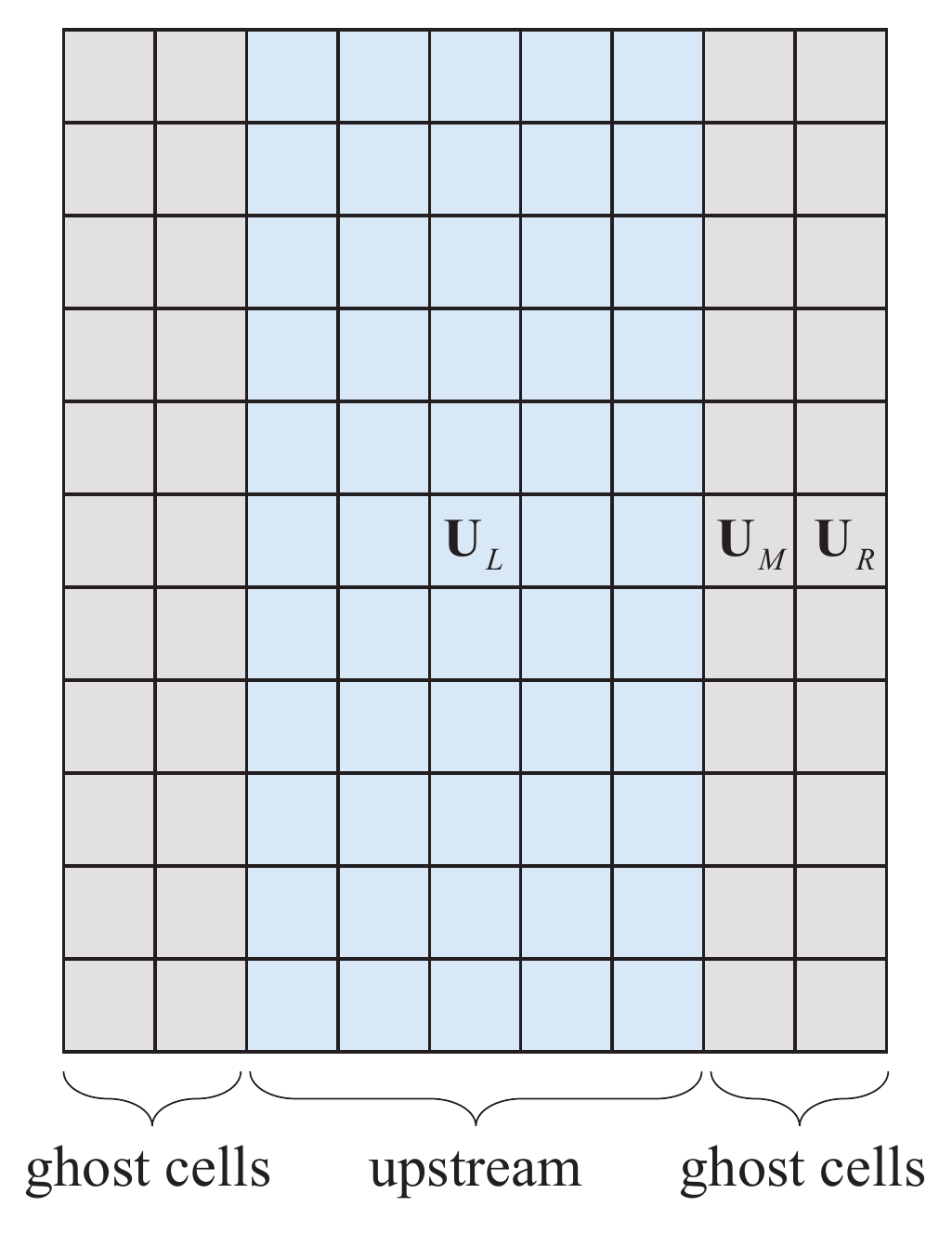}
	\end{minipage}
	}
	\subfigure[The numerical shock structure domain.]{
	\begin{minipage}[t]{0.31\linewidth}
	\centering
	\includegraphics[height=0.95\textwidth]{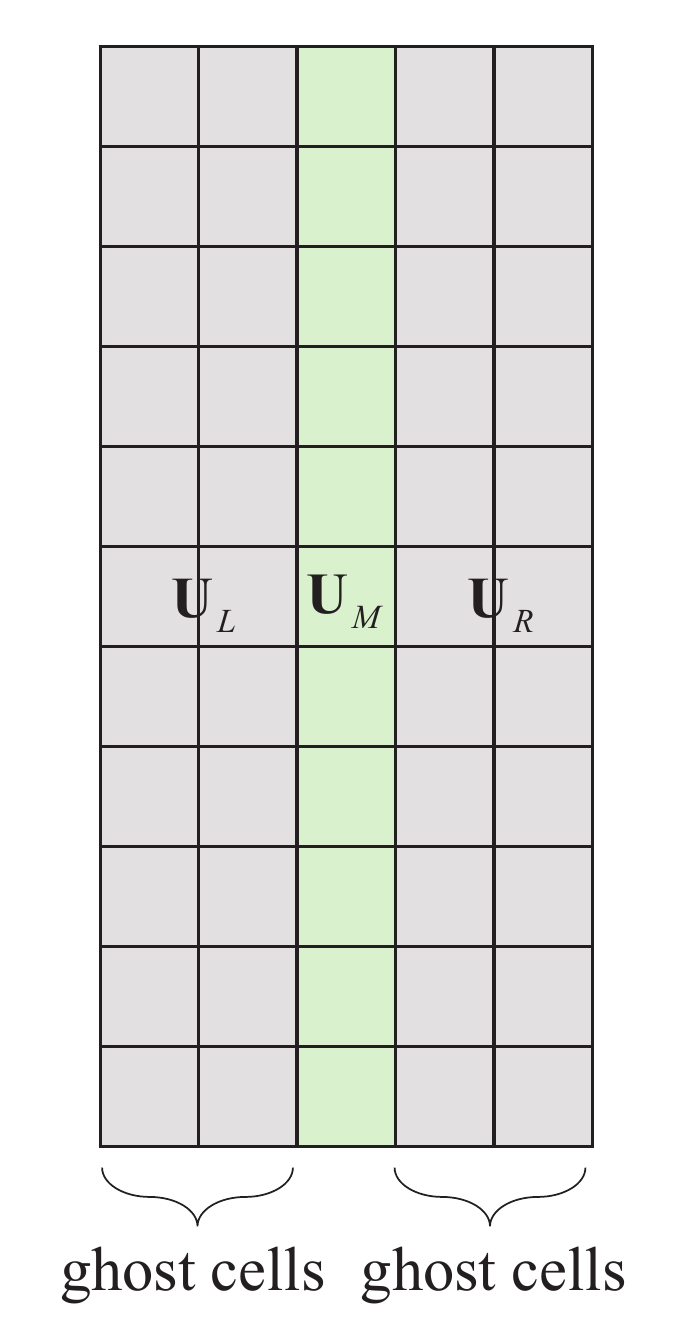}
	\end{minipage}
	}
	\subfigure[The downstream domain.]{
	\begin{minipage}[t]{0.31\linewidth}
	\centering
	\includegraphics[height=0.95\textwidth]{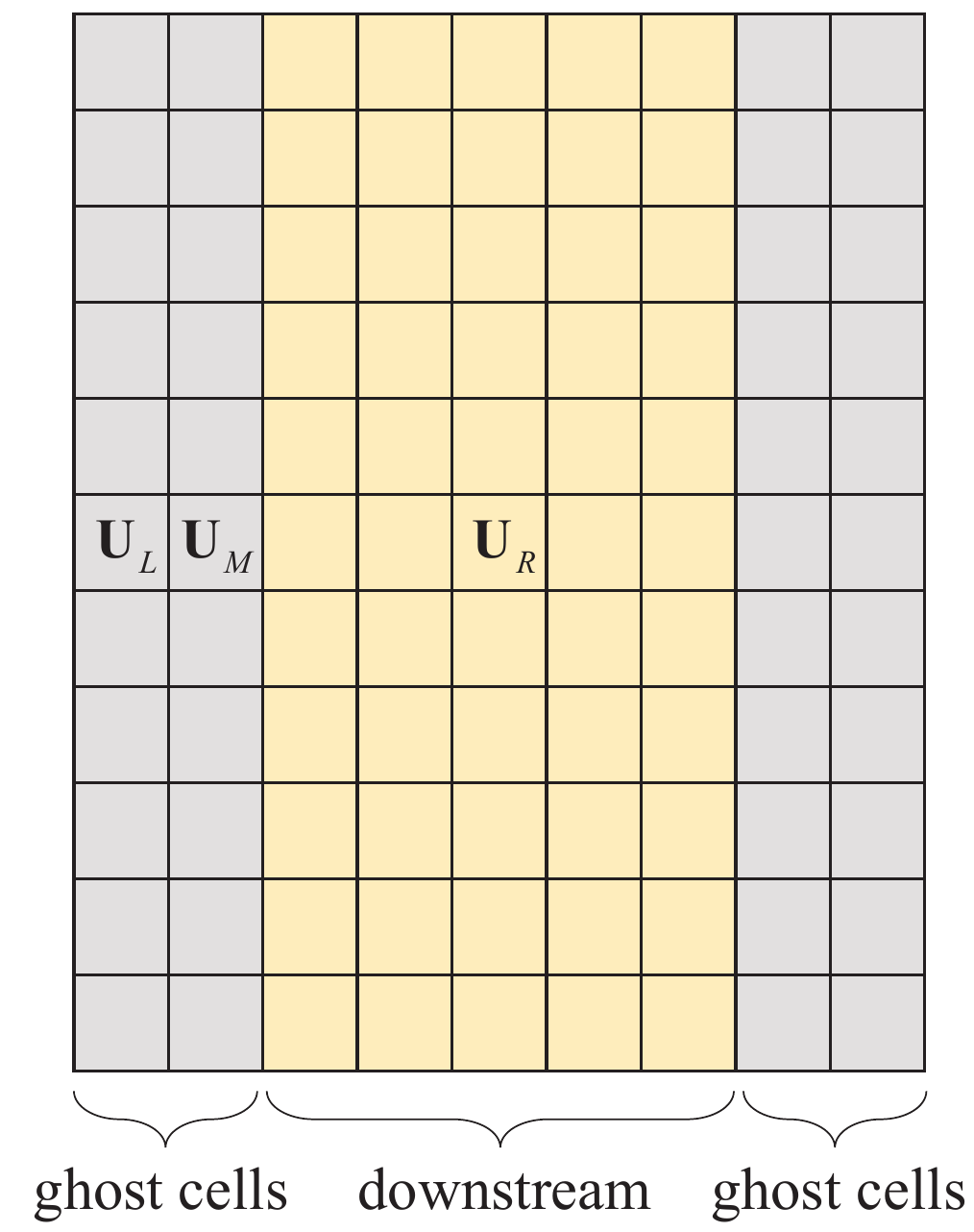}
	\end{minipage}
	}
	\centering
	\caption{Three test cases for investigating the spatial localization of the	source of instability.}\label{fig boundary of different domains}
\end{figure}

\begin{figure}[htbp]
	\centering
	\subfigure[Upstream.]{
	\begin{minipage}[t]{0.46\linewidth}
	\centering
	\includegraphics[width=0.95\textwidth]{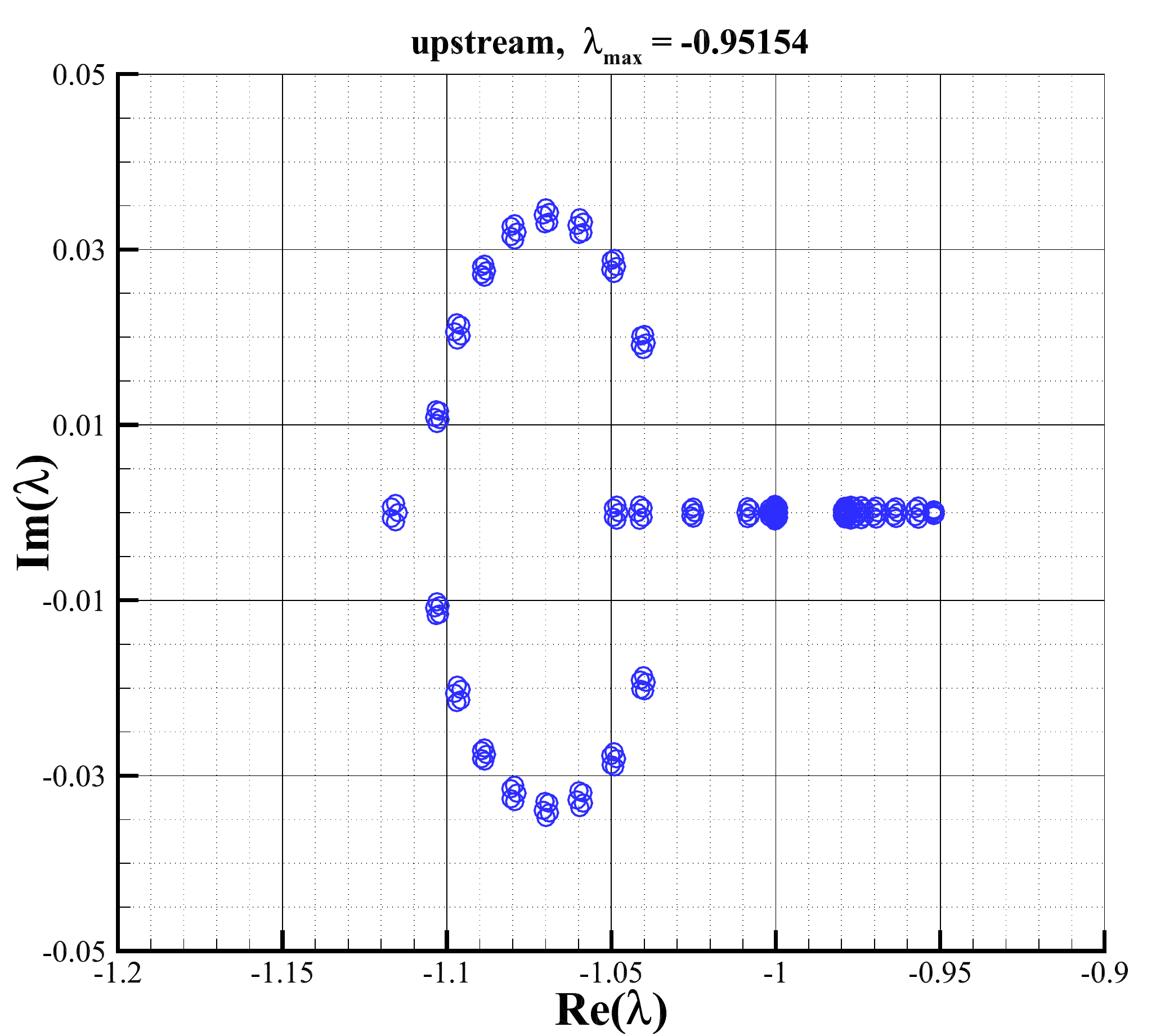}
	\end{minipage}
	}
	\subfigure[Downstream.]{
	\begin{minipage}[t]{0.46\linewidth}
	\centering
	\includegraphics[width=0.95\textwidth]{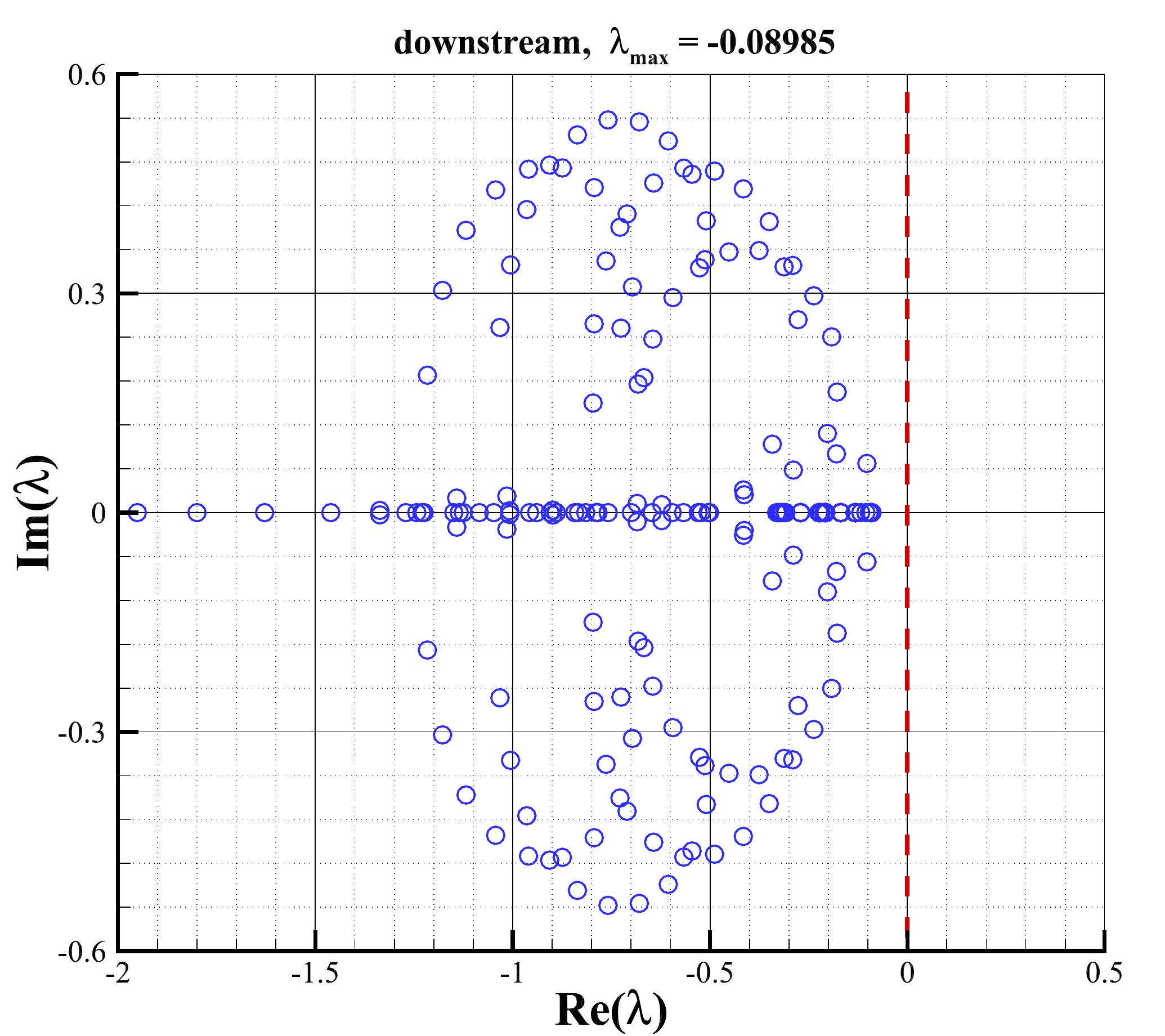}
	\end{minipage}
	}

	\subfigure[Numerical shock structure.]{
	\begin{minipage}[t]{0.46\linewidth}
	\centering
	\includegraphics[width=0.95\textwidth]{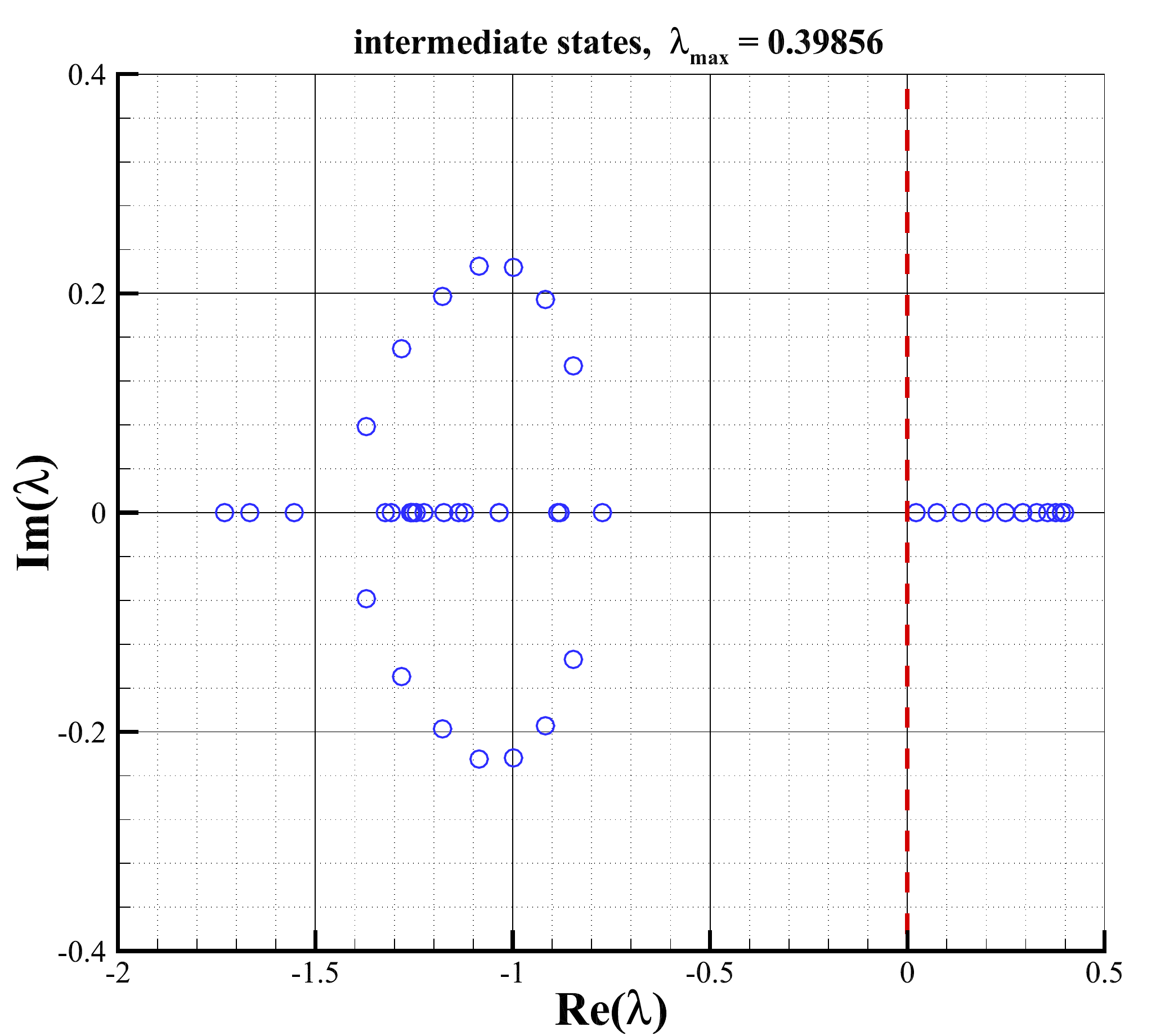}
	\end{minipage}
	}
	\centering
	\caption{Distribution of the eigenvalues of \textbf{S} in the complex plane for three test cases.($ M_0=20 $, second-order scheme with Roe solver and van Albada limiter.)}\label{fig scatters of three domains}
\end{figure}

\begin{figure}[htbp]
	\centering
	\subfigure[The unstable mode for $ \rho $.]{
	\begin{minipage}[t]{0.46\linewidth}
	\centering
	\includegraphics[width=0.95\textwidth]{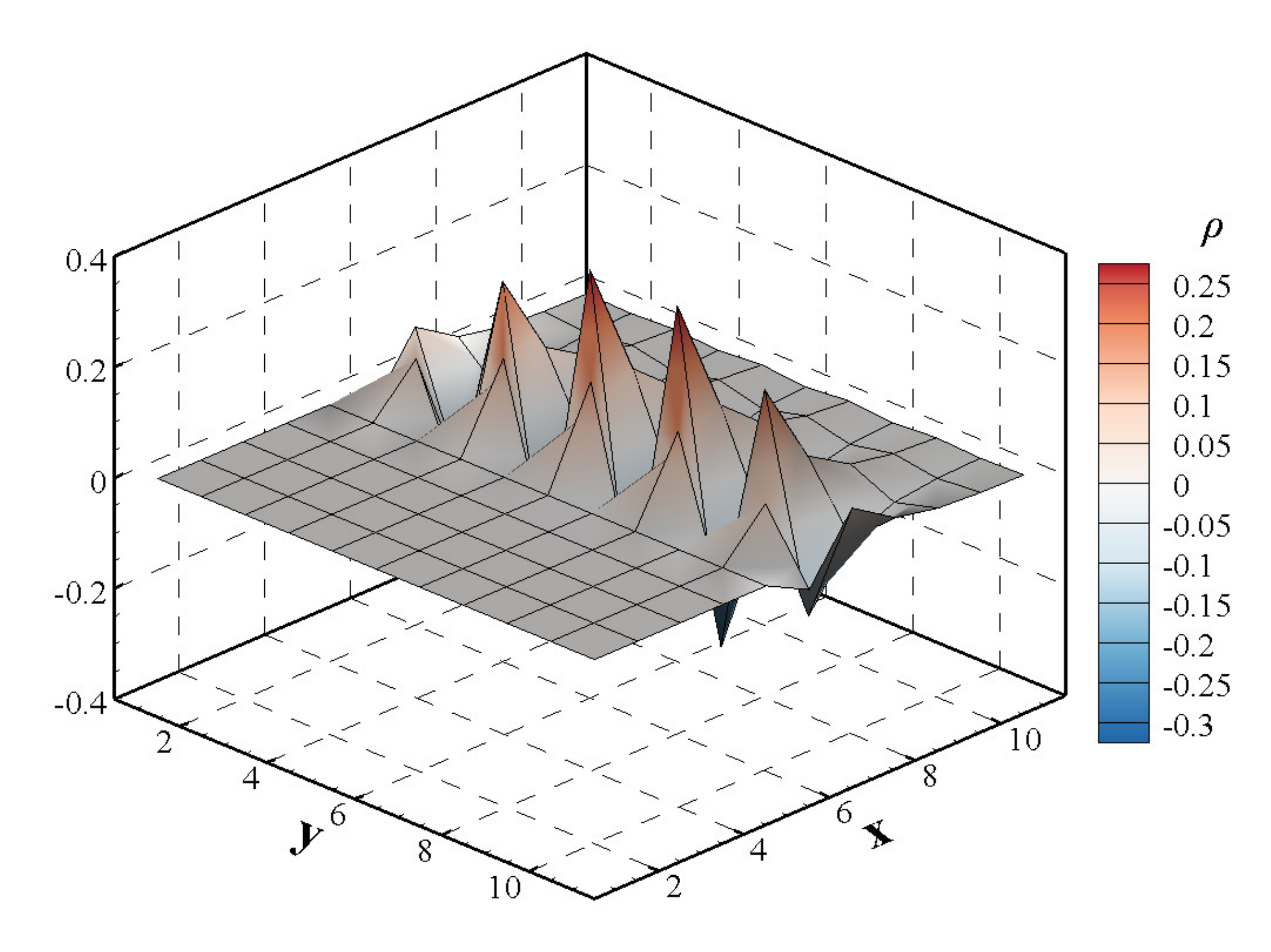}
	\end{minipage}
	}
	\subfigure[The unstable mode for $ u $.]{
	\begin{minipage}[t]{0.46\linewidth}
	\centering
	\includegraphics[width=0.95\textwidth]{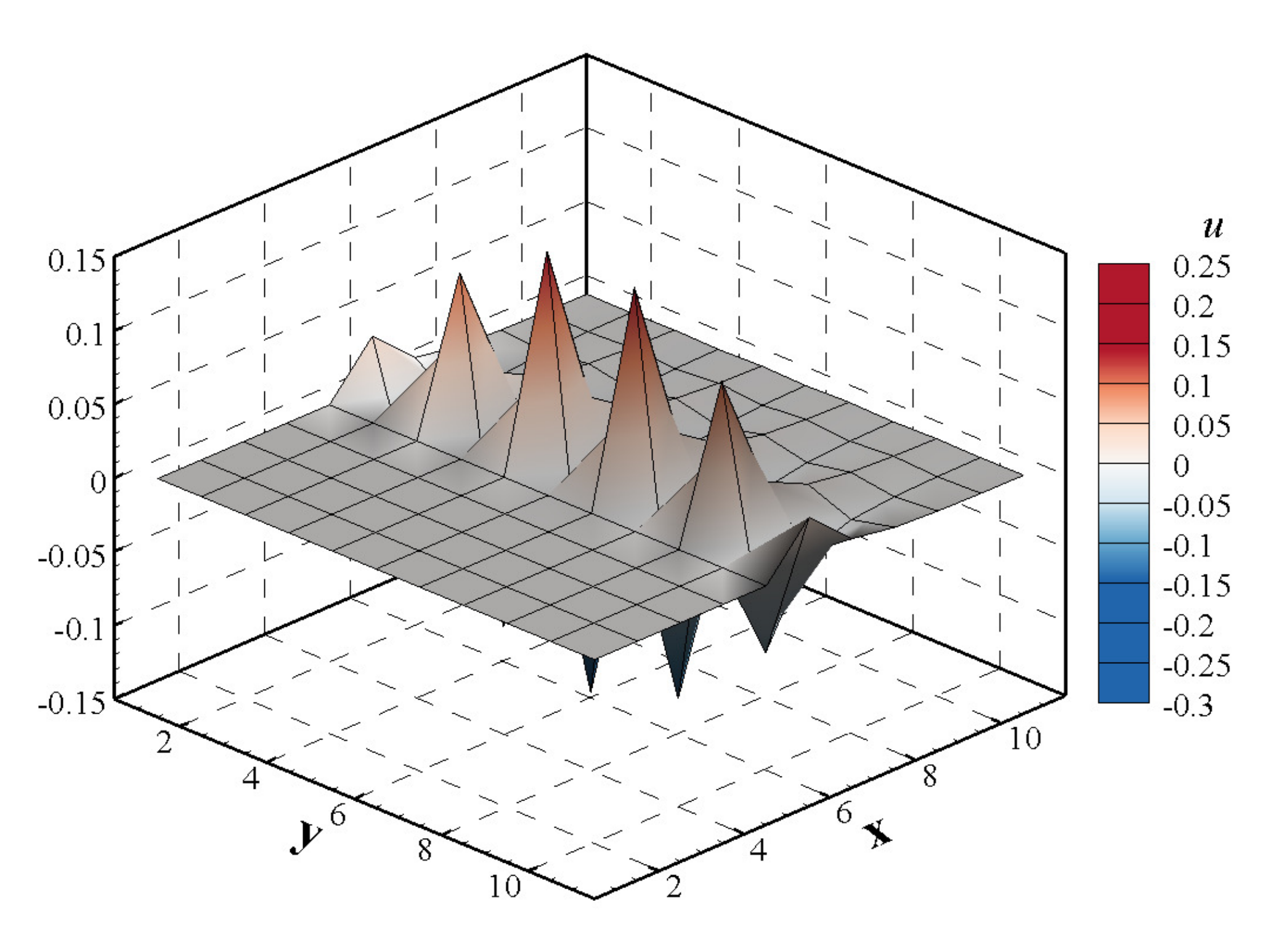}
	\end{minipage}
	}
	
	\subfigure[The unstable mode for $ v $.]{
	\begin{minipage}[t]{0.46\linewidth}
	\centering
	\includegraphics[width=0.95\textwidth]{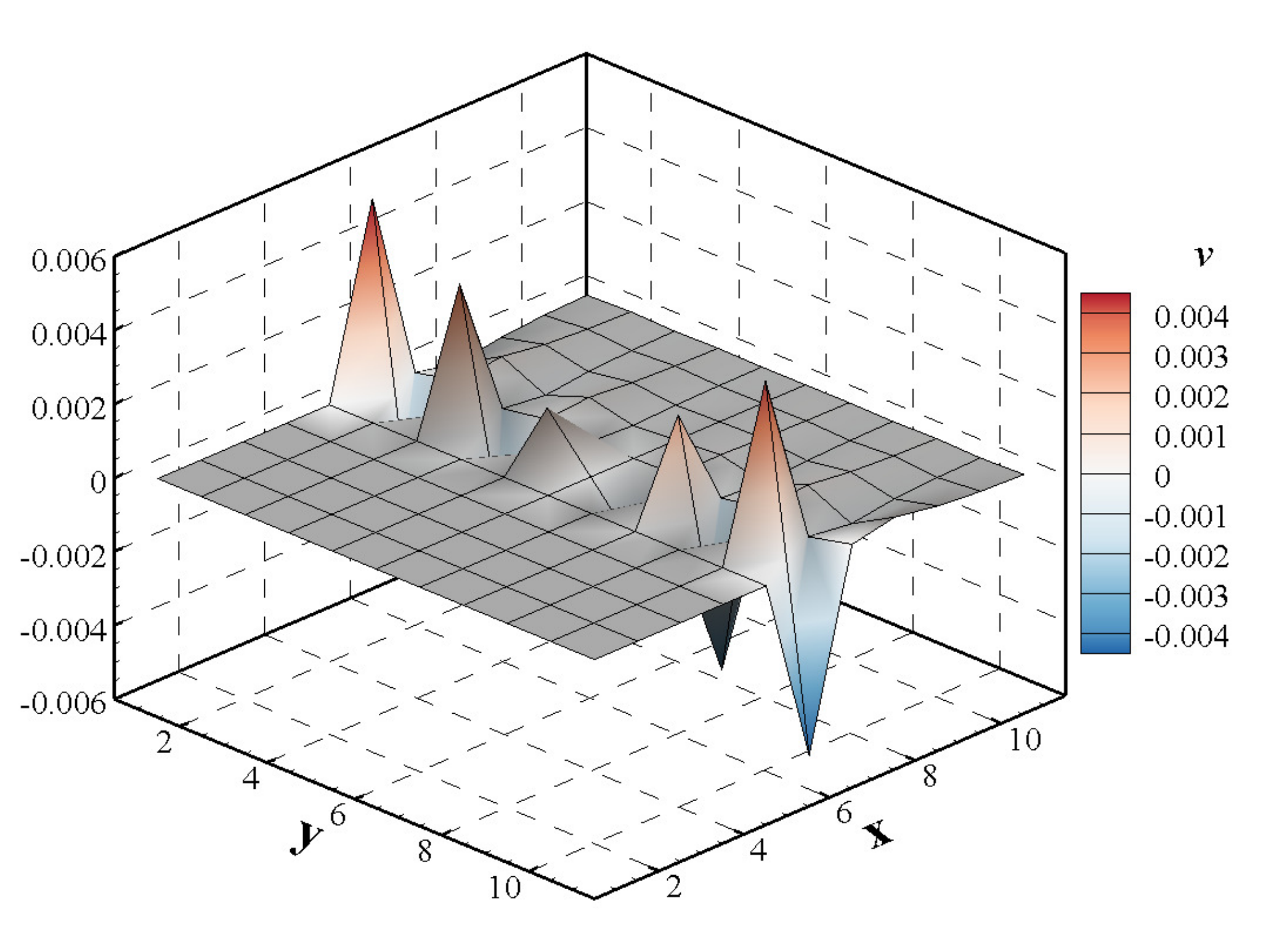}
	\end{minipage}
	}
	\subfigure[The unstable mode for $ p $.]{
	\begin{minipage}[t]{0.46\linewidth}
	\centering
	\includegraphics[width=0.95\textwidth]{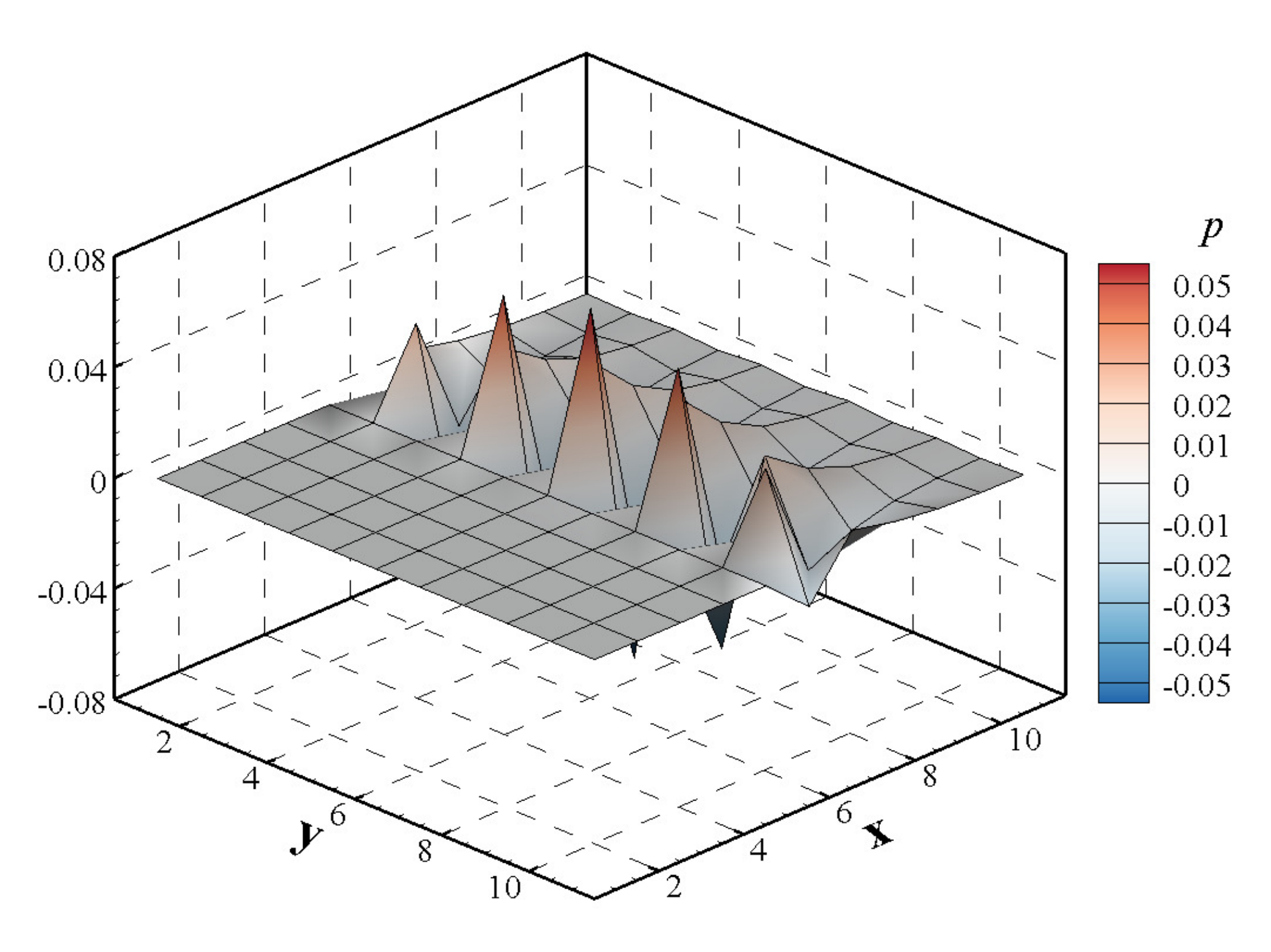}
	\end{minipage}
	}
	\centering
	\caption{The unstable mode for $ \lambda =0.76291+0i $.}\label{fig unstable mode}
\end{figure}

As shown in Fig.\ref{fig boundary of different domains}, in order to explore the location of the source of instability, three test cases are employed. In the first one (Fig.\ref{fig boundary of different domains}(a)), all cells apart from the ghost cells are set to the upstream state $ \mathbf{U}_L $ while the right boundary conditions are set to intermediate state $ \mathbf{U}_M $ and downstream state $ \mathbf{U}_R $. Since the boundaries are supposed to be error-free, so their direct contribution to error evolution vanishes in the stability matrix except from the flux-gradient \cite{Dumbser2004}. This procedure prevents any error from developing in numerical shock structure and downstream region, as a result of which we can analyze the stability of the upstream region separately. Identically, in order to analyze the stability of numerical shock structure and downstream region, similar procedures are shown in Fig.\ref{fig boundary of different domains}(b) and (c).\par

The stabilities of three test cases are investigated using the matrix stability analysis method. Fig.\ref{fig scatters of three domains} shows the distribution of the eigenvalues in the complex plane, from which we can find that the maximal real part of the eigenvalues of the numerical shock structure exceeds 0, while the others are less than 0. As a result, we can infer that in the cases of the shock on the upstream boundary and downstream boundary, the stability analysis predicts a stable behavior, while if the numerical shock structure is contained, the analysis predicts unconditional instability. Moreover, the maximal real part of eigenvalues for the downstream field is larger than that of the upstream field, which may imply that the downstream is more susceptible to being influenced by the unstable information from the numerical shock structure than the upstream. Additionally, the information of the spatial-behavior of the unstable mode is contained in the eigenvectors of the maximal eigenvalue \cite{Dumbser2004}. Fig.\ref{fig unstable mode} depicts the unstable mode for the four primitive variables. As shown in Fig.\ref{fig unstable mode}, the numerical shock structure is the most unstable location, and unstable information will propagate more easily downstream. Thus, according to the analysis in this section, we can conclude that the shock instability originates from the numerical shock structure. And the downstream is more susceptible to being influenced by spurious errors in the shock structure.\par

\section{Conclusion}\label{section 6}
In this paper, we devote our efforts to quantitatively investigating the shock instability of the second-order scheme. To this end, the matrix stability analysis method for the finite-volume MUSCL approach is developed. With the help of the matrix stability analysis, primary numerical characteristics and the underlying mechanism of shock instability for the second-order scheme are investigated in detail. Results reveal that the shock instability is greatly affected by the shock intensity. The computation for shocks will be more unstable as the Mach number increases. Moreover, the Riemann solver is an important factor affecting the shock stability, whose effect on the second-order scheme is consistent with that of the first-order one. One of the most significant findings to emerge from the current study is that the differences in shock instability between the first and second-order scheme are revealed. When employing high dissipative solvers, e.g. HLL, the stability of capturing strong shocks will be better as the spatial accuracy is enhanced to second-order. However, compared with the first-order scheme, the second-order scheme with low dissipative solvers, such as Roe, is prone to be more unstable. For the second-order scheme, the limiter is demonstrated to be a vital factor in triggering the shock instability. Compared with the low dissipative limiter, the limiter with higher dissipation can decelerate the process towards an unstable result. Except for the factors discussed above, the results have also shown that the computational grid is still a critical element that will significantly affect the shock instability, which will be alleviated by increasing the aspect ratio and distortion angle. Furthermore, the underlying mechanism of the shock instability is still investigated in the current study. And it has been demonstrated that the shock instability originates from the numerical shock structure. The computation for shocks will be more stable when the interval shock conditions are closer to the downstream states.\par

The current work offers an effective tool to quantitatively analyze the shock stability of the second-order finite-volume MUSCL approach and provides some clues to better understand and control the occurrence of the carbuncle phenomenon or shock instability. However, numerical evidences have shown that high-order schemes are supposed to be at a higher risk of shock instability, and it is more challenging to analyze it in a quantitative manner. Thus, such a quantitative tool that allows to investigate the shock stability of high-order schemes is urgently needed. Corresponding numerical characteristics and the underlying mechanism of high-order shock-capturing methods are more attractive to be investigated in depth. That is what we are addressing now and will be presented in the follow-up work.

\section*{Acknowledgements}	
This work was supported by National Natural Science Foundation of China (Grant No.12202490), Natural Science Foundation of Hunan Province, China (Grant No. 11472004), the Scientific Research Foundation of NUDT (Grant No. ZK21-10), and Postgraduate Scientific Research Innovation Project of Hunan Province (Grant Nos. CX20220010 and CX20220036).

\bibliography{mybibfile}

\end{document}